 \documentclass[a4paper,12pt]{article}

 \usepackage{amsfonts,amsmath,mathrsfs,amssymb}
 \usepackage{times,helvet,courier,type1cm}

 \usepackage{makeidx}
 \usepackage{harvard}
 \usepackage{color}

 \allowdisplaybreaks

 \makeindex

 \newtheorem{thm0}{Theorem}[section]
 \newtheorem{exa0}{Theorem}[section]

 \newtheorem{def1}[thm0]{Definition}
 \newtheorem{lem1}[thm0]{Lemma}
 \newtheorem{thm1}[thm0]{Theorem}
 \newtheorem{cor1}[thm0]{Corollary}
 \newtheorem{pro1}[thm0]{Proposition}
 \newtheorem{con1}[thm0]{Condition}
 \newtheorem{exa1}[exa0]{\it{Example}}

 \def\bglemma{\begin{lem1}}\def\edlemma{\end{lem1}}
 \def\bgtheorem{\begin{thm1}}\def\edtheorem{\end{thm1}}
 \def\bgcorollary{\begin{cor1}}\def\edcorollary{\end{cor1}}
 \def\bgproposition{\begin{pro1}}\def\edproposition{\end{pro1}}
 \def\bgcondition{\begin{con1}}\def\edcondition{\end{con1}}
 \def\bgexample{\begin{exa1}\rm{}\def\edexample{\end{exa1}}}

 \def\benumerate{\begin{enumerate}}\def\eenumerate{\end{enumerate}}
 \def\bitemize{\begin{itemize}}\def\eitemize{\end{itemize}}
 \def\itm{\item}

 \def\beqlb{\begin{eqnarray}}\def\eeqlb{\end{eqnarray}}
 \def\beqnn{\begin{eqnarray*}}\def\eeqnn{\end{eqnarray*}}

 \def\eqref#1{{\rm(\ref{#1})}}

 \def\nnm{\nonumber}\def\ccr{\nnm\\}

 \def\mbb{\mathbb}\def\mrm{\mathrm}
 \def\mbf{\mathbf}\def\mcr{\mathscr}

 \def\d{\mrm{d}}\def\e{\mrm{e}}\def\b{\mrm{b}}

 \def\bE{\mbf{E}}\def\bP{\mbf{P}}
 \def\bQ{\mbf{Q}}\def\bN{\mbf{N}}

 \def\proof{\noindent{\textit{Proof.~~}}}
 \def\qed{\hfill$\square$\smallskip}

 \def\ulim{\uparrow\!\!\lim}\def\dlim{\downarrow\!\!\lim}

 \def\qqquad{\qquad\qquad}\def\bbigskip{\bigskip\bigskip}

 \def\ar{\!\!&}

\begin{document}

\thispagestyle{empty}

\

\vskip5.0cm

\centerline{\huge\sf{\textbf{CONTINUOUS-STATE BRANCHING}}}

\bigskip

\centerline{\huge\sf{\textbf{PROCESSES WITH IMMIGRATION}}}

\vskip3.0cm

\centerline{\Large\sf{Zenghu Li}}

\vskip5.0cm


\centerline{\large\rm{(Beijing Normal University)}}

\newpage


\

\thispagestyle{empty}

\bbigskip

Zenghu Li, Professor

School of Mathematical Sciences

Beijing Normal University

Beijing 100875, China

E-mail: lizh@bnu.edu.cn

URL: http://math0.bnu.edu.cn/\~{}lizh/

\vskip12.0cm





\newpage

$~$

\pagenumbering{roman}\tableofcontents\newpage$~$


\newpage

\pagenumbering{arabic}

\addcontentsline{toc}{section}{\bf Introduction}

\section*{Introduction}

 \setcounter{equation}{0}

Continuous-state branching processes (CB-processes) and continuous-state branching processes with immigration (CBI-processes) constitute important classes of Markov processes taking values in the positive ($=\,$nonnegative) half line. They were introduced as probabilistic models describing the evolution of large populations with small individuals. The study of CB-processes was initiated by Feller (1951), who noticed that a diffusion process may arise in a limit theorem of Galton--Watson discrete branching processes; see also Aliev and Shchurenkov (1982), Grimvall (1974) and Lamperti (1967a). A characterization of CB-processes by random time changes of L\'{e}vy processes was given by Lamperti (1967b). The convergence of rescaled discrete branching processes with immigration to CBI-processes was studied in Aliev (1985), Kawazu and Watanabe (1971) and Li (2006). {From} a mathematical point of view, the continuous-state processes are usually easier to deal with because both their time and state spaces are smooth, and the distributions that appear are infinitely divisible. For general treatments and backgrounds of CB- and CBI-processes, the reader may refer to Kyprianou (2014) and Li (2011). In the recent work of Pardoux (2016), more complicated probabilistic population models involving competition were studied, which extend the stochastic logistic growth model of Lambert (2005).

A continuous CBI-process with subcritical branching mechanism was used by Cox et al.\ (1985) to describe the evolution of interest rates and has been known in mathematical finance as the \textit{Cox--Ingersoll--Ross model} (CIR-model). Compared with other financial models introduced before, the CIR-model is more appealing as it is positive and mean-reverting. The asymptotic behavior of the estimators of the parameters in this model was studied by Overbeck and Ryd\'en (1997); see also Li and Ma (2015). Applications of stochastic calculus to finance including those of the CIR-model were discussed systematically in Lamberton and Lapeyre (1996). A natural generalization of the CBI-process is the so-called affine Markov process, which has also been used a lot in mathematical finance; see, e.g., Duffie et al.\ (2003) and the references therein.

A strong stochastic equation for general CBI-processes was first established in Dawson and Li (2006). A flow of discontinuous CB-processes was constructed in Bertoin and Le~Gall (2006) by weak solutions to a stochastic equation. Their results were extended to flows of CBI-processes in Dawson and Li (2012) using strong solutions; see also Li (2014) and Li and Ma (2008). For the stable branching CBI-process, a strong stochastic differential equation driven by L\'{e}vy processes was established in Fu and Li (2010). The approach of stochastic equations has played an important role in recent developments of the theory and applications of CB- and CBI-processes.

The purpose of these notes is to provide a brief introduction to CB- and CBI-processes accessible to graduate students with reasonable background in probability theory and stochastic processes. In particular, we give a quick development of the stochastic equations of the processes and some immediate applications. The proofs given here are more elementary than those appearing in the literature before. We have made them readable without requiring too much preliminary knowledge on branching processes and stochastic analysis.

In Section~1, we review some properties of Laplace transforms of finite measures on the positive half line. In Section~2, a construction of CB-processes is given as rescaling limits of Galton--Watson branching processes. This approach also gives the physical interpretation of the CB-processes. Some basic properties of the processes are developed in Section~3. The Laplace transforms of some positive integral functionals are calculated explicitly in Section~4. In Section~5, the CBI-processes are constructed as rescaling limits of Galton--Watson branching processes with immigration. In Section~6, we present reconstructions of the CB- and CBI-processes by Poisson random measures determined by entrance laws, which reveal the structures of the trajectories of the processes. Several equivalent formulations of martingale problems for CBI-processes are given in Section~7. {From} those we derive the stochastic equations of the processes in Section~8. Using the stochastic equations, some characterizations of local and global maximal jumps of the CB- and CBI-processes are given in Section~9. In Section~10, we prove the strong Feller property and the exponential ergodicity of the CBI-process under suitable conditions using a coupling based on one of the stochastic equations.

These lecture notes originated from graduate courses I gave at Beijing Normal University in the past years. They were also used for mini courses at Peking University in 2017 and at the University of Verona in 2018. I would like to thank Professors Ying Jiao and Simone Scotti, who invited me to give the mini courses. I am grateful to the participants of all those courses for their helpful comments. I would also like to thank NSFC for the financial supports. I am indebted to the Laboratory of Mathematics and Complex Systems (Ministry of Education) for providing me the research facilities.

\newpage

\section{Laplace transforms of measures}

 \setcounter{equation}{0}

Let $\mcr{B}[0,\infty)$ be the Borel $\sigma$-algebra on the positive half line $[0,\infty)$. Let $B[0,\infty) = \b\mcr{B}[0,\infty)$ be the set of bounded Borel
functions on $[0,\infty)$. Given a finite measure $\mu$ on $[0,\infty)$, we
define the \index{Laplace transform} \textit{Laplace transform} $L_\mu$
of $\mu$ by
 \beqlb\label{s1.1}
L_\mu(\lambda)
 =
\int_{[0,\infty)} \e^{-\lambda x} \mu(\d x), \qquad \lambda\ge 0.
 \eeqlb

\bgtheorem\label{ts1.1} A finite measure on $[0,\infty)$ is uniquely determined by its Laplace
transform. \edtheorem

\proof Suppose that $\mu_1$ and $\mu_2$ are finite measures on $[0,\infty)$ and
$L_{\mu_1}(\lambda) = L_{\mu_2}(\lambda)$ for all $\lambda\ge 0$. Let
$\mcr{K} = \{x\mapsto \e^{-\lambda x}: \lambda\ge 0\}$ and let $\mcr{L}$
be the class of functions $F\in B[0,\infty)$ so that
 \beqnn
\int_{[0,\infty)} F(x)\mu_1(\d x) = \int_{[0,\infty)} F(x)\mu_2(\d x).
 \eeqnn
Then $\mcr{K}$ is closed under multiplication and $\mcr{L}$ is a monotone
vector space containing $\mcr{K}$. It is easy to see $\sigma(\mcr{K}) =
\mcr{B}[0,\infty)$. Then the monotone class theorem implies
$\mcr{L}\supset \b\sigma(\mcr{K}) = B[0,\infty)$. That proves the
desired result. \qed

\bgtheorem\label{ts1.2} Let $\{\mu_n\}$ be a sequence of finite measures on $[0,\infty)$ and $\lambda\mapsto L(\lambda)$ a continuous function on $[0,\infty)$. If $\lim_{n\to \infty} L_{\mu_n}(\lambda) = L(\lambda)$ for every $\lambda\ge 0$, then there is a finite measure $\mu$ on $[0,\infty)$ such that $L_\mu=L$ and $\lim_{n\to \infty} \mu_n = \mu$ by weak convergence. \edtheorem

\proof We can regard each $\mu_n$ as a finite measure on $[0,\infty]$, the one-point compactification of $[0,\infty)$. Let $F_n$ denote the distribution function of $\mu_n$. By Helly's theorem we infer that $\{F_n\}$ contains a subsequence $\{F_{n_k}\}$ that converges weakly on $[0,\infty]$ to some distribution function $F$. Then the corresponding subsequence $\{\mu_{n_k}\}$ converges weakly on $[0,\infty]$ to the finite measure $\mu$ determined by $F$. It follows that
 \beqlb\label{s1.2}
\mu[0,\infty] = \lim_{k\to \infty} \mu_{n_k}[0,\infty]
 =
\lim_{k\to \infty} \mu_{n_k}[0,\infty) = \lim_{k\to \infty} L_{\mu_{n_k}}(0) = L(0).
 \eeqlb
Moreover, for $\lambda> 0$ we have
 \beqlb\label{s1.3}
\int_{[0,\infty]} \e^{-\lambda x} \mu(\d x)
 \ar=\ar
\lim_{k\to \infty} \int_{[0,\infty]} \e^{-\lambda x} \mu_{n_k}(\d x) \cr
 \ar=\ar
\lim_{k\to \infty} \int_{[0,\infty)} \e^{-\lambda x} \mu_{n_k}(\d x)
 =
L(\lambda),
 \eeqlb
where $\e^{-\lambda\cdot\infty} = 0$ by convention. By letting $\lambda\to 0+$ in \eqref{s1.3} and using the continuity of $L$ at $\lambda=0$ we find $\mu[0,\infty) = L(0)$. {From} this and \eqref{s1.2} we see $\mu$ is supported by $[0,\infty)$. By Theorem~1.7 of Li (2011, p.4) we have $\lim_{n\to \infty} \mu_{n_k} = \mu$ weakly on $[0,\infty)$. It follows that, for $\lambda\ge 0$,
 \beqnn
\int_{[0,\infty)} \e^{-\lambda x} \mu(\d x)
 =
\lim_{k\to \infty} \int_{[0,\infty)} \e^{-\lambda x} \mu_{n_k}(\d x)
 =
L(\lambda).
 \eeqnn
Then $L_\mu=L$. If $\mu_n$ does not converge weakly to $\mu$, then $F_n$ does not converge weakly to $F$, so there is a subsequence $\{F_{n_k^\prime}\}\subset \{F_n\}$ that converges weakly to a limit $G\neq F$. The above arguments show that $G$ corresponds to a finite measure on $[0,\infty)$ with Laplace transform $L=L_\mu$, yielding a contradiction. Then $\lim_{n\to \infty} \mu_n = \mu$ weakly on $[0,\infty)$. \qed

\bgcorollary\label{ts1.3} Let $\mu_1, \mu_2, \ldots$ and $\mu$ be finite measures on $[0,\infty)$. Then $\mu_n \to \mu$ weakly if and only if $L_{\mu_n}(\lambda)\to L_\mu(\lambda)$ for every $\lambda\ge 0$. \edcorollary

\proof If $\mu_n \to \mu$ weakly, we have $\lim_{n\to \infty} L_{\mu_n}(\lambda) = L_\mu(\lambda)$ for every $\lambda\ge 0$ by dominated convergence. The converse assertion is a consequence of Theorem~\ref{ts1.2}. \qed

Given two probability measures $\mu_1$ and $\mu_2$ on $[0,\infty)$, we denote by
$\mu_1\times \mu_2$ their product measure on $[0,\infty)^2$. The image of $\mu_1\times \mu_2$ under the mapping $( x_1, x_2)\mapsto x_1+ x_2$
is called the \index{convolution of probability measures} \textit{convolution} of $\mu_1$ and
$\mu_2$ and is denoted by $\mu_1*\mu_2$, which is a probability measure
on $[0,\infty)$. According to the definition, for any $F\in B[0,\infty)$
we have
 \beqlb\label{s1.4}
\int_{[0,\infty)} F(x) (\mu_1*\mu_2)(\d x)
 =
\int_{[0,\infty)} \mu_1(\d x_1)\int_{[0,\infty)} F(x_1+x_2) \mu_2(\d x_2).
 \eeqlb
Clearly, if $\xi_1$ and $\xi_2$ are independent random variables with
distributions $\mu_1$ and $\mu_2$ on $[0,\infty)$, respectively, then the
random variable $\xi_1+\xi_2$ has distribution $\mu_1*\mu_2$. It is easy
to show that
 \beqlb\label{s1.5}
L_{\mu_1*\mu_2}(\lambda) = L_{\mu_1}(\lambda)L_{\mu_2}(\lambda), \qquad
\lambda\ge 0.
 \eeqlb
Let $\mu^{*0} = \delta_0$ and define $\mu^{*n} = \mu^{*(n-1)}*\mu$ inductively for integers $n\ge 1$.

We say a probability distribution $\mu$ on $[0,\infty)$ is \index{infinitely divisible distribution} \textit{infinitely divisible} if for each integer $n\ge1$, there is a probability $\mu_n$ such that $\mu = \mu_n^{*n}$. In this case, we call $\mu_n$ the \index{the $n$-th root of a probability} \textit{$n$-th root} of $\mu$. A positive random variable $\xi$ is said to be \index{infinitely divisible random variable} \textit{infinitely divisible} if it has infinitely divisible distribution on $[0,\infty)$. Write $\psi\in \mcr{I}$ if $\lambda\mapsto \psi(\lambda)$ is a positive function on $[0,\infty)$ with the \index{L\'{e}vy--Khintchine representation} \textit{L\'{e}vy--Khintchine representation}:
 \beqlb\label{s1.6}
\psi(\lambda) = h\lambda + \int_{(0,\infty)} (1-\e^{-\lambda u}) l(\d u),
 \eeqlb
where $h\ge 0$ and $l(\d u)$ is a $\sigma$-finite measure on $(0,\infty)$ satisfying
 \beqnn
\int_{(0,\infty)} (1\land u) l(\d u)< \infty.
 \eeqnn
The relation $\psi = -\log L_\mu$ establishes a one-to-one correspondence between the functions $\psi\in \mcr{I}$ and infinitely divisible probability measures $\mu$ on $[0,\infty)$; see, e.g., Theorem~1.39 in Li (2011, p.20).

\newpage

\section{Construction of CB-processes}

 \setcounter{equation}{0}

Let $\{p(j): j\in \mbb{N}\}$ be a probability distribution on the space of positive integers $\mbb{N} := \{0,1,2,\ldots\}$. It is well-known that $\{p(j): j\in \mbb{N}\}$ is uniquely determined by its generating function $g$ defined by
 \beqnn
g(z) = \sum_{j=0}^\infty p(j)z^j, \qquad |z|\le 1.
 \eeqnn
Suppose that $\{\xi_{n,i}: n,i=1,2,\ldots\}$ is a family of $\mbb{N}$-valued i.i.d.\ random variables with distribution $\{p(j): j\in \mbb{N}\}$. Given an $\mbb{N}$-valued random variable $x(0)$ independent of $\{\xi_{n,i}\}$, we define inductively
 \beqlb\label{s2.1}
x(n) = \sum_{i=1}^{x(n-1)}\xi_{n,i}, \qquad n=1,2,\ldots.
 \eeqlb
Here we understand $\sum_{i=1}^0=0$. For $i\in \mbb{N}$ let $\{Q(i,j): j\in \mbb{N}\}$ denote the $i$-fold convolution of $\{p(j): j\in \mbb{N}\}$, that is, $Q(i,j) = p^{*i}(j)$ for $i,j\in \mbb{N}$. For any $n\ge 1$ and $\{i_0,i_1, \cdots, i_{n-1}=i,j\}\subset \mbb{N}$ it is easy to see that
 \beqnn
\ar\ar\bP\Big(x(n) = j\big|x(0)=i_0, x(1)=i_1, \cdots, x(n-1)=i_{n-1}\Big) \cr
 \ar\ar\qquad
= \bP\bigg(\sum_{i=1}^{x(n-1)}\xi_{n,i}\Big|x(n-1)=i_{n-1}\bigg) \cr
 \ar\ar\qquad
= \bP\bigg(\sum_{k=1}^{i}\xi_{n,k} = j\bigg) = Q(i,j).
 \eeqnn
Then $\{x(n): n\ge 0\}$ is an $\mbb{N}$-valued Markov chain with one-step transition matrix $Q= (Q(i,j): i, j\in \mbb{N})$. The random variable $x(n)$ can be thought of as the number of individuals in generation $n$ of an evolving population system. After one unit time, each individual in the population splits independently of others into a random number of offspring according to the distribution $\{p(j): j\in \mbb{N}\}$. Clearly, we have, for $i\in \mbb{N}$ and $|z|\le 1$,
 \beqlb\label{s2.2}
\sum^\infty_{j=0} Q(i,j)z^j
 =
\sum^\infty_{j=0} p^{*i}(j)z^j
 =
g(z)^i.
 \eeqlb
Clearly, the transition matrix $Q$ satisfies the \index{branching property} \textit{branching property}:
 \beqlb\label{s2.3}
Q(i_1+i_2,\cdot) = Q(i_1,\cdot)*Q(i_2,\cdot), \qquad i_1,i_2\in \mbb{N}.
 \eeqlb
This means that different individuals in the population propagate independently each other.

A Markov chain in $\mbb{N}$ with one-step transition matrix defined by \eqref{s2.2} is called a \index{Galton--Watson branching process} \textit{Galton--Watson branching process} \index{GW-process} (GW-process) or a \index{Bienaym\'{e}--Galton--Watson branching process} \textit{Bienaym\'{e}--Galton--Watson branching process} \index{BGW-process} (BGW-process) with \index{branching distribution} \textit{branching distribution} given by $g$; see, e.g., Athreya and Ney (1972) and Harris (1963). The study of the model goes back to Bienaym\'{e} (1745) and Galton and Watson (1874).

By a general result in the theory of Markov chains, for any $n\ge 1$ the $n$-step transition matrix of the GW-process is just the $n$-fold product matrix $Q^n= (Q^n(i,j): i, j\in \mbb{N})$.

\bgproposition\label{ts2.1} For any $n\ge 1$ and $i\in \mbb{N}$ we have
 \beqlb\label{s2.4}
\sum^\infty_{j=0} Q^n(i,j)z^j = g^{\circ n}(z)^i,
\qquad |z|\le 1,
 \eeqlb
where $g^{\circ n}(z)$ is defined by $g^{\circ n}(z) = g\circ g^{\circ (n-1)}(z) = g(g^{\circ (n-1)}(z))$ successively with $g^{\circ 0}(z) = z$ by convention.
\edproposition

\proof {From} \eqref{s2.2} we know \eqref{s2.4} holds for $n=1$. Now suppose that \eqref{s2.4} holds for some $n\ge 1$. We have
 \beqnn
\sum^\infty_{j=0} Q^{n+1}(i,j)z^j
 \ar=\ar
\sum^\infty_{j=0} \sum^\infty_{k=0} Q(i,k)Q^n(k,j)z^j \cr
 \ar=\ar
\sum^\infty_{k=0} Q(i,k)g^{\circ n}(z)^k = g^{\circ(n+1)}(z)^i.
 \eeqnn
Then \eqref{s2.4} also holds when $n$ is replaced by $n+1$. That proves the result by induction. \qed

It is easy to see that zero is a trap for the GW-process. If $g^\prime(1-)<\infty$, by differentiating both sides of \eqref{s2.4} we see the first moment of the distribution $\{Q^n(i,j): j\in \mbb{N}\}$ is given by
 \beqlb\label{s2.5}
\sum^\infty_{j=1} jQ^n(i,j) = ig^\prime(1-)^n.
 \eeqlb

\bgexample\label{es2.1} Given a GW-process $\{x(n): n\ge 0\}$, we can define its \index{extinction time} \textit{extinction time} $\tau_0 = \inf\{n\ge 0: x(n)=0\}$. In view of \eqref{s2.1}, we have $x(n)=0$ on the event $\{n\ge \tau_0\}$. Let $q = \bP(\tau_0< \infty|x(0)=1)$ be the \index{extinction probability} \textit{extinction probability}. By the independence of the propagation of different individuals we have $\bP(\tau_0< \infty| x(0) = i) = q^i$ for any $i=0,1,2,\dots$. By the total probability formula,
 \beqnn
q \ar=\ar \sum_{j=0}^\infty\bP(x(1)=j|x(0)=1)\bP(\tau_0< \infty|x(0)=1,x(1)=j) \cr
 \ar=\ar
\sum_{j=0}^\infty\bP(\xi_{1,1}=j)\bP(\tau_0< \infty|x(1)=j)
 =
\sum_{j=0}^\infty p(j)q^j = g(q).
 \eeqnn
Then the extinction probability $q$ is a solution to the equation $z=g(z)$ on $[0,1]$. Clearly, in the case of $p(1)< 1$ we have $q=1$ if and only if $\sum_{j=1}^\infty jp(j)\le 1$. \edexample

Now suppose we have a sequence of GW-processes $\{x_k(n): n\ge 0\}$, $k=1,2,\dots$ with branching distributions given by the probability generating functions $g_k$, $k=1,2,\dots$. Let $z_k(n) = k^{-1}x_k(n)$. Then $\{z_k(n): n\ge 0\}$ is a Markov chain with state space $E_k := \{0,k^{-1},2k^{-1},\ldots\}$ and $n$-step transition probability $Q_k^n(x,\d y)$ determined by
 \beqlb\label{s2.6}
\int_{E_k}\e^{-\lambda y}Q_k^n(x,\d y)
 =
g_k^{\circ n}(\e^{-\lambda/k})^{kx}, \qquad \lambda\ge 0.
 \eeqlb
Suppose that $\{\gamma_k\}$ is a positive sequence so that $\gamma_k\to \infty$ as $k\to \infty$. Let $\lfloor\gamma_kt\rfloor$ denote the integer part of $\gamma_kt$. Clearly, given $z_k(0) = x\in E_k$, for any $t\ge 0$ the random variable $z_k(\lfloor\gamma_kt\rfloor) = k^{-1}x_k(\lfloor\gamma_kt\rfloor)$
has distribution $Q_k^{\lfloor\gamma_kt\rfloor}(x,\cdot)$ on $E_k$ determined by
 \beqlb\label{s2.7}
\int_{E_k}\e^{-\lambda y}Q_k^{\lfloor\gamma_kt\rfloor}(x,\d y)
 =
\exp\{-xv_k(t,\lambda)\},
 \eeqlb
where
 \beqlb\label{s2.8}
v_k(t,\lambda)
 =
-k\log g_k^{\circ\lfloor\gamma_kt\rfloor}(\e^{-\lambda/k}).
 \eeqlb

We are interested in the asymptotic behavior of the sequence of continuous time processes $\{z_k(\lfloor\gamma_kt\rfloor): t\ge 0\}$ as $k\to \infty$. By \eqref{s2.8}, for $\gamma_k^{-1}(i-1)\le t< \gamma_k^{-1}i$ we have
 \beqnn
v_k(t,\lambda) = v_k(\gamma_k^{-1}\lfloor\gamma_k t\rfloor,\lambda)
 =
v_k(\gamma_k^{-1}(i-1),\lambda).
 \eeqnn
It follows that
 \beqlb\label{s2.9}
v_k(t,\lambda)
 \ar=\ar
v_k(0,\lambda) + \sum_{j=1}^{\lfloor\gamma_kt\rfloor} [v_k(\gamma_k^{-1}j,\lambda) - v_k(\gamma_k^{-1}(j-1),\lambda)] \cr
 \ar=\ar
\lambda - k\sum_{j=1}^{\lfloor\gamma_kt\rfloor} [\log g_k^{\circ j}(\e^{-\lambda/k}) - \log g_k^{\circ (j-1)}(\e^{-\lambda/k})] \cr
 \ar=\ar
\lambda - k\sum_{j=1}^{\lfloor\gamma_kt\rfloor} \log \big[g_k(g_k^{\circ (j-1)}(\e^{-\lambda/k}))g_k^{\circ (j-1)}(\e^{-\lambda/k})^{-1}\big] \cr
 \ar=\ar
\lambda - \gamma_k^{-1}\sum_{j=1}^{\lfloor\gamma_kt\rfloor}\bar{\phi}_k(-k\log g_k^{\circ (j-1)}(\e^{-\lambda/k})) \cr
 \ar=\ar
\lambda - \gamma_k^{-1}\sum_{j=1}^{\lfloor\gamma_kt\rfloor}\bar{\phi}_k(v_k(\gamma_k^{-1}(j-1),
\lambda)) \cr
 \ar=\ar
\lambda - \int_0^{\gamma_k^{-1}\lfloor\gamma_kt\rfloor} \bar{\phi}_k(v_k(s,\lambda))\d s,
 \eeqlb
where
 \beqlb\label{s2.10}
\bar{\phi}_k(z)= k\gamma_k\log\big[g_k(\e^{-z/k})\e^{z/k}\big], \qquad z\ge 0.
 \eeqlb
It is easy to see that
 \beqlb\label{s2.11}
\bar{\phi}_k(z)= k\gamma_k\log\big[1+(k\gamma_k)^{-1}\tilde{\phi}_k(z)\e^{z/k}\big],
 \eeqlb
where
 \beqlb\label{s2.12}
\tilde{\phi}_k(z)= k\gamma_k[g_k(\e^{-z/k})-\e^{-z/k}].
 \eeqlb
The sequence $\{\tilde{\phi}_k\}$ is sometimes easier to handle than the original sequence $\{\bar{\phi}_k\}$. The following lemma shows that the two sequences are really not very different.

\bglemma\label{ts2.2} Suppose that the sequence $\{\tilde{\phi}_k\}$ is uniformly bounded on each bounded interval. Then we have: {\rm(i)} $\lim_{k\to\infty} |\bar{\phi}_k(z) - \tilde{\phi}_k(z)| = 0$ uniformly on each bounded interval; {\rm(ii)} $\{\bar{\phi}_k\}$ is uniformly Lipschitz on each bounded interval if and only if so is $\{\tilde{\phi}_k\}$. \edlemma

\proof The first assertion follows immediately from \eqref{s2.11}. By the same relation we have
 \beqnn
\bar{\phi}_k^\prime(z) = {[\tilde{\phi}_k^\prime(z) + k^{-1}\tilde{\phi}_k(z)]\e^{z/k}\over 1+(k\gamma_k)^{-1}\tilde{\phi}_k(z)\e^{z/k}}, \qquad z\ge 0.
 \eeqnn
Then $\{\bar{\phi}_k^\prime\}$ is uniformly bounded on each bounded interval if and only if so is $\{\tilde{\phi}_k^\prime\}$. That gives the second assertion. \qed

By the above lemma, if either $\{\tilde{\phi}_k\}$ or $\{\bar{\phi}_k\}$ is uniformly Lipschitz on each bounded interval, then they converge or diverge simultaneously and in the convergent case they have the same limit. For the convenience of statement of the results, we formulate the following condition:

\bgcondition\label{ts2.3} The sequence $\{\tilde{\phi}_k\}$ is uniformly Lipschitz on $[0,a]$ for every $a\ge 0$ and there is a function $\phi$ on $[0,\infty)$ so that $\tilde{\phi}_k(z)\to \phi(z)$ uniformly on $[0,a]$ for every $a\ge 0$ as $k\to \infty$. \edcondition

\bgproposition\label{ts2.4} Suppose that Condition~\ref{ts2.3} is satisfied. Then the limit function $\phi$ has representation
 \beqlb\label{s2.13}
\phi(z) = bz + cz^2 + \int_{(0,\infty)}\big(\e^{-zu}-1+zu\big)m(\d u), \quad
z\ge 0,
 \eeqlb
where $c\ge 0$ and $b$ are constants and $m(\d u)$ is a $\sigma$-finite measure on $(0,\infty)$ satisfying
 \beqnn
\int_{(0,\infty)} (u\land u^2)m(\d u)< \infty.
 \eeqnn
\edproposition

\proof For each $k\ge 1$ let us define the function $\phi_k$ on $[0,k]$ by
 \beqlb\label{s2.14}
\phi_k(z) = k\gamma_k[g_k(1-z/k) - (1-z/k)].
 \eeqlb
{From} \eqref{s2.12} and \eqref{s2.14} we have
 \beqnn
\tilde{\phi}^\prime_k(z)
 =
\gamma_k\e^{-z/k}[1-g^\prime_k(\e^{-z/k})], \qquad z\ge 0,
 \eeqnn
and
 \beqnn
\phi^\prime_k(z)
 =
\gamma_k[1-g^\prime_k(1-z/k)], \qquad 0\le z\le k.
 \eeqnn
Since $\{\tilde{\phi}_k\}$ is uniformly Lipschitz on each bounded interval, the sequence $\{\tilde{\phi}_k^\prime\}$ is uniformly bounded on each bounded interval. Then $\{\phi_k^\prime\}$ is also uniformly bounded on each bounded interval, and so the sequence $\{\phi_k\}$ is uniformly Lipschitz on each bounded interval. Let $a\ge 0$. By the mean-value theorem, for $k\ge a$ and $0\le z\le a$ we have
 \beqnn
\tilde{\phi}_k(z) - \phi_k(z)
 \ar=\ar
k\gamma_k\big[g_k(\e^{-z/k}) - g_k(1-z/k) - \e^{-z/k} + (1 - z/k)\big]
\quad \cr
 \ar=\ar
k\gamma_k[g^\prime_k(\eta_k)-1](\e^{-z/k}-1+z/k),
 \eeqnn
where
 \beqnn
1-a/k\le 1-z/k\le \eta_k\le \e^{-z/k}\le 1.
 \eeqnn
Choose $k_0\ge a$ so that $\e^{-2a/k_0}\le 1-a/k_0$. Then $\e^{-2a/k}\le
1-a/k$ for $k\ge k_0$ and hence
 \beqnn
\gamma_k |g^\prime_k(\eta_k)-1|
 \le
\sup_{0\le z\le 2a} \gamma_k|g^\prime_k(\e^{-z/k})-1|
 =
\sup_{0\le z\le 2a} \e^{z/k}|\tilde{\phi}^\prime_k(z)|.
 \eeqnn
Since $\{\tilde{\phi}^\prime_k\}$ is uniformly bounded on $[0,2a]$, the sequence $\{\gamma_k |g^\prime_k(\eta_k)-1|: k\ge k_0\}$ is bounded. Then $\lim_{k\to\infty} |\phi_k(z)-\tilde{\phi}_k(z)| = 0$ uniformly on each bounded interval. It follows that $\lim_{k\to \infty}\phi_k(z) = \phi(z)$ uniformly on each bounded interval. Then the result follows by Corollary~1.46 in Li (2011, p.26). \qed

\bgproposition\label{ts2.5} For any function $\phi$ with representation \eqref{s2.13} there is a sequence $\{\tilde{\phi}_k\}$ in the form of \eqref{s2.12} satisfying
Condition~\ref{ts2.3}. \edproposition

\proof By the proof of Proposition~\ref{ts2.4} it suffices to construct a sequence $\{\phi_k\}$ with the expression \eqref{s2.14} that is uniformly Lipschitz on $[0,a]$ and $\phi_k(z)\to \phi(z)$ uniformly on $[0,a]$ for every $a\ge 0$. To simplify the formulations we decompose the function $\phi$ into two parts. Let $\phi_0(z) = \phi(z) - bz$. We first define
 \beqnn
\gamma_{0,k} = (1+2c)k + \int_{(0,\infty)} u(1-\e^{-ku}) m(\d u)
 \eeqnn
and
 \beqnn
g_{0,k}(z) = z + k^{-1}\gamma_{0,k}^{-1} \phi_0(k(1-z)), \qquad |z|\le 1.
 \eeqnn
It is easy to see that $z\mapsto g_{0,k}(z)$ is an analytic function satisfying $g_{0,k}(1) = 1$ and
 \beqnn
\frac{\d^n}{\d z^n}g_{0,k}(0)\ge 0, \qquad n\ge 0.
 \eeqnn
Therefore $g_{0,k}(\cdot)$ is a probability generating function. Let $\phi_{0,k}$ be defined by \eqref{s2.14} with $(\gamma_k,g_k)$ replaced by $(\gamma_{0,k}, g_{0,k})$. Then $\phi_{0,k}(z) = \phi_0(z)$ for $0\le z\le k$. That completes the proof if $b=0$. In the case $b\neq 0$, we set
 \beqnn
g_{1,k}(z) = \frac{1}{2}\bigg(1+\frac{b}{|b|}\bigg) +
\frac{1}{2}\bigg(1-\frac{b}{|b|}\bigg)z^2.
 \eeqnn
Let $\gamma_{1,k} = |b|$ and let $\phi_{1,k}(z)$ be defined by \eqref{s2.14} with $(\gamma_k,g_k)$ replaced by $(\gamma_{1,k}, g_{1,k})$. Then
 \beqnn
\phi_{1,k}(z) = bz + \frac{1}{2k}(|b|-b)z^2.
 \eeqnn
Finally, let $\gamma_k = \gamma_{0,k} + \gamma_{1,k}$ and $g_k = \gamma_k^{-1}(\gamma_{0,k}g_{0,k} + \gamma_{1,k}g_{1,k})$. Then the sequence $\phi_k(z)$ defined by \eqref{s2.14} is equal to $\phi_{0,k}(z) + \phi_{1,k}(z)$ which satisfies the required condition. \qed

\bglemma\label{ts2.6} Suppose that the sequence $\{\tilde{\phi}_k\}$ defined by \eqref{s2.12} is uniformly Lipschitz on $[0,1]$. Then there are constants $B,N\ge 0$ such that $v_k(t,\lambda)\le \lambda\e^{Bt}$ for every $t,\lambda\ge 0$ and
$k\ge N$. \edlemma

\proof Let $b_k := \tilde{\phi}_k^\prime(0+)$ for $k\ge 1$. Since $\{\tilde{\phi}_k\}$ is uniformly Lipschitz on $[0,1]$, the sequence $\{b_k\}$ is bounded. {From} \eqref{s2.12} we have $b_k = \gamma_k[1-g_k^\prime(1-)]$. By \eqref{s2.5} and \eqref{s2.12} it is not hard to obtain
 \beqnn
\int_{E_k} y Q_k^{\lfloor\gamma_kt\rfloor}(x,\d y)
 =
xg_k^\prime(1-)^{\lfloor\gamma_kt\rfloor}
 =
x\bigg(1-\frac{b_k}{\gamma_k}\bigg)^{\lfloor\gamma_kt\rfloor}.
 \eeqnn
Let $B\ge 0$ be a constant such that $2|b_k|\le B$ for all $k\ge 1$. Since $\gamma_k\to \infty$ as $k\to \infty$, there is $N\ge 1$ so that
 \beqnn
0\le \bigg(1-\frac{b_k}{\gamma_k}\bigg)^{\gamma_k/B}
 \le
\bigg(1+\frac{B}{2\gamma_k}\bigg)^{\gamma_k/B}
 \le
\e, \qquad k\ge N.
 \eeqnn
It follows that, for $t\ge 0$ and $k\ge N$,
 \beqnn
\int_{E_k} y Q_k^{\lfloor\gamma_kt\rfloor}(x,\d y)
 \le
x\exp\big\{B\lfloor\gamma_kt\rfloor/\gamma_k\big\}
 \le
x\e^{Bt}.
 \eeqnn
Then the desired estimate follows from \eqref{s2.5} and Jensen's inequality. \qed

\bgtheorem\label{ts2.7} Suppose that Condition~\ref{ts2.3} holds. Then for every $a\ge 0$ we
have $v_k(t,\lambda)\to$ some $v_t(\lambda)$ uniformly on $[0,a]^2$ as
$k\to \infty$ and the limit function solves the integral equation
 \beqlb\label{s2.15}
v_t(\lambda) = \lambda - \int_0^t\phi(v_s(\lambda))\d s, \qquad
\lambda,t\ge 0.
 \eeqlb
\edtheorem

\proof The following argument is a modification of that of Aliev and Shchurenkov (1982) and Aliev (1985). In view of \eqref{s2.9}, we can write
 \beqlb\label{s2.16}
v_k(t,\lambda)
 =
\lambda + \varepsilon_k(t,\lambda) - \int_0^t \bar{\phi}_k(v_k(s,\lambda))\d s,
 \eeqlb
where
 \beqnn
\varepsilon_k(t,\lambda)
 =
\big(t-\gamma_k^{-1}\lfloor\gamma_kt\rfloor\big)
\bar{\phi}_k\big(v_k(\gamma_k^{-1}\lfloor\gamma_kt\rfloor,\lambda)\big).
 \eeqnn
By Lemma~\ref{ts2.2} and Condition~\ref{ts2.3}, for any $0<\varepsilon\le 1$ we can choose $N\ge 1$ so that $|\bar{\phi}_k(z)-\phi(z)|\le \varepsilon$ for $k\ge N$ and $0\le z\le a\e^{Ba}$. It follows that, for $0\le t\le a$ and $0\le \lambda\le a$,
 \beqlb\label{s2.17}
|\varepsilon_k(t,\lambda)|
 \le
\gamma_k^{-1}\big|\bar{\phi}_k\big(v_k(\gamma_k^{-1}\lfloor\gamma_kt\rfloor, \lambda)\big)\big|
 \le
\gamma_k^{-1}M,
 \eeqlb
where
 \beqnn
M = 1+\sup_{0\le z\le a\e^{Ba}}|\phi(z)|.
 \eeqnn
For $n\ge k\ge N$ let
 \beqnn
K_{k,n}(t,\lambda)
 =
\sup_{0\le s\le t}|v_n(s,\lambda)-v_k(s,\lambda)|.
 \eeqnn
By \eqref{s2.16} and \eqref{s2.17} we obtain, for $0\le t\le a$ and $0\le \lambda\le a$,
 \beqnn
K_{k,n}(t,\lambda)
 \ar\le\ar
2\gamma_k^{-1}M + \int_0^t | \bar{\phi}_k(v_k(s,\lambda)) - \bar{\phi}_n(v_n(s,\lambda))| \d s \cr
 \ar\le\ar
2(\gamma_k^{-1}M + \varepsilon a) + \int_0^t | \phi_k(v_k(s,\lambda)) - \phi_n(v_n(s,\lambda))| \d s \cr
 \ar\le\ar
2(\gamma_k^{-1}M + \varepsilon a) + L\int_0^t K_{k,n}(s,\lambda) \d s,
 \eeqnn
where $L = \sup_{0\le z\le a\e^{Ba}}|\phi^\prime(z)|$. By Gronwall's inequality,
 \beqnn
K_{k,n}(t,\lambda)
 \le
2(\gamma_k^{-1}M + \varepsilon a)\exp\{Lt\}, \qquad 0\le t,\lambda\le a.
 \eeqnn
Then $v_k(t,\lambda)\to$ some $v_t(\lambda)$ uniformly on $[0,a]^2$ as
$k\to \infty$ for every $a\ge 0$. {From} \eqref{s2.16} we get
\eqref{s2.15}. \qed

\bgtheorem\label{ts2.8} Suppose that $\phi$ is a function given by \eqref{s2.13}. Then for any $\lambda\ge 0$ there is a unique positive solution $t\mapsto v_t(\lambda)$ to \eqref{s2.15}. Moreover, the solution satisfies the \index{semigroup property} \textit{semigroup property}:
 \beqlb\label{s2.18}
v_{r+t}(\lambda) = v_r\circ v_t(\lambda) = v_r(v_t(\lambda)), \qquad
r,t,\lambda\ge 0.
 \eeqlb
\edtheorem

\proof By Proposition~\ref{ts2.5} there is a sequence $\{\tilde{\phi}_k\}$ in form \eqref{s2.12} satisfying Condition~\ref{ts2.3}. Let $v_k(t,\lambda)$ be given by \eqref{s2.7} and \eqref{s2.8}. By Theorem~\ref{ts2.7} the limit $v_t(\lambda)= \lim_{k\to \infty} v_k(t,\lambda)$ exists and solves \eqref{s2.15}. Clearly, any positive solution $t\mapsto v_t(\lambda)$ to \eqref{s2.15} is locally bounded. The uniqueness of the solution follows by Gronwall's inequality. The relation \eqref{s2.18} is a consequence of the uniqueness of the solution. \qed

\bgtheorem\label{ts2.9} Suppose that $\phi$ is a function given by \eqref{s2.13}. For any $\lambda\ge 0$ let $t\mapsto v_t(\lambda)$ be the unique positive solution to \eqref{s2.15}. Then we can define a transition semigroup $(Q_t)_{t\ge 0}$ on $[0,\infty)$ by
 \beqlb\label{s2.19}
\int_{[0,\infty)} \e^{-\lambda y} Q_t(x,\d y)
 =
\e^{-xv_t(\lambda)}, \qquad \lambda\ge 0,x\ge 0.
 \eeqlb
\edtheorem

\proof By Proposition~\ref{ts2.5}, there is a sequence $\{\tilde{\phi}_k\}$ in form \eqref{s2.12} satisfying Condition~\ref{ts2.3}. By Theorem~\ref{ts2.7} we have $v_k(t,\lambda)\to v_t(\lambda)$ uniformly on $[0,a]^2$ as $k\to \infty$ for every $a\ge 0$. Taking $x_k\in E_k$ satisfying $x_k\to x$ as $k\to \infty$, we see by Theorem~\ref{ts1.2} that \eqref{s2.19} defines a probability measure $Q_t(x,\d y)$ on $[0,\infty)$ and $\lim_{k\to \infty} Q_k^{\lfloor\gamma_k t\rfloor}(x_k,\cdot) = Q_t(x,\cdot)$ by weak convergence. By a monotone class argument one can see that $Q_t(x,\d y)$ is a kernel on $[0,\infty)$. The semigroup property of the family of kernels $(Q_t)_{t\ge 0}$ follows from \eqref{s2.18} and \eqref{s2.19}. \qed

\bgproposition\label{ts2.10} For every $t\ge 0$ the function $\lambda\mapsto v_t(\lambda)$ is strictly increasing on $[0,\infty)$. \edproposition

\proof By the continuity of $t\mapsto v_t(\lambda)$, for any $\lambda_0>0$ there is $t_0>0$ so that $v_t(\lambda_0)>0$ for $0\le t\le t_0$. Then \eqref{s2.19} implies $Q_t(x,\{0\})< 1$ for $x>0$ and $0\le t\le t_0$, and so $\lambda\mapsto v_t(\lambda)$ is strictly increasing for $0\le t\le t_0$. By the semigroup property \eqref{s2.18} we infer $\lambda\mapsto v_t(\lambda)$ is strictly increasing for all $t\ge 0$. \qed

\bgtheorem\label{ts2.11} The transition semigroup $(Q_t)_{t\ge 0}$ defined by \eqref{s2.19} is a Feller semigroup. \edtheorem

\proof For $\lambda>0$ and $x\ge 0$ set $e_\lambda(x) = \e^{-\lambda x}$. We denote by $\mcr{D}_0$ the linear span of $\{e_\lambda: \lambda>0\}$. By Proposition~\ref{ts2.10}, the operator $Q_t$ preserves $\mcr{D}_0$ for every $t\ge 0$. By the continuity of $t\mapsto v_t(\lambda)$ it is easy to show that $t\mapsto Q_te_\lambda(x)$ is continuous for $\lambda> 0$ and $x\ge 0$. Then $t\mapsto Q_tf(x)$ is continuous for every $f\in \mcr{D}_0$ and $x\ge 0$. Let $C_0[0,\infty)$ be the space of continuous functions on $[0,\infty)$ vanishing at infinity. By the Stone--Weierstrass theorem, the set $\mcr{D}_0$ is uniformly dense in $C_0[0,\infty)$; see, e.g., Hewitt and Stromberg (1965, pp.98-99). Then each operator $Q_t$ preserves $C_0[0,\infty)$ and $t\mapsto Q_tf(x)$ is continuous for $x\ge 0$ and $f\in C_0[0,\infty)$. That gives the Feller property of the semigroup $(Q_t)_{t\ge 0}$. \qed

A Markov process in $[0,\infty)$ is called a \index{continuous-state branching process} \textit{continuous-state branching process} \index{CB-process} (CB-process) with \index{branching mechanism} \textit{branching mechanism} $\phi$ if it has transition semigroup $(Q_t)_{t\ge 0}$ defined by \eqref{s2.19}. It is simple to see that $(Q_t)_{t\ge 0}$ satisfies the \index{branching property} \textit{branching property}:
 \beqlb\label{s2.20}
Q_t(x_1+x_2,\cdot)
 =
Q_t(x_1,\cdot)*Q_t(x_2,\cdot), \qquad t,x_1,x_2\ge 0.
 \eeqlb
The family of functions $(v_t)_{t\ge 0}$ is called the \index{cumulant semigroup} \textit{cumulant semigroup} of the CB-process. Since $(Q_t)_{t\ge 0}$ is a Feller semigroup, the process has a Hunt realization; see, e.g., Chung (1982, p.75). Clearly, zero is a trap for the CB-process.

\bgproposition\label{ts2.12} Suppose that $\{(x_1(t),\mcr{F}_t^1): t\ge 0\}$ and $\{(x_2(t),\mcr{F}_t^2): t\ge 0\}$ are two independent CB-processes with branching mechanism $\phi$. Let $x(t) = x_1(t) + x_2(t)$ and $\mcr{F}_t = \sigma(\mcr{F}_t^1\cup \mcr{F}_t^2)$. Then $\{(x(t),\mcr{F}_t): t\ge 0\}$ is also a CB-processes with branching mechanism $\phi$. \edproposition

\proof Let $t\ge r\ge 0$ and for $i=1,2$ let $F_i$ be a bounded $\mcr{F}_r^i$-measurable random variable. For any $\lambda\ge 0$ we have
 \beqnn
\bP\big[F_1F_2\e^{-\lambda x(t)}\big]
 \ar=\ar
\bP\big[F_1\e^{-\lambda x_1(t)}\big] \bP\big[F_2\e^{-\lambda x_2(t)}\big] \cr
 \ar=\ar
\bP\big[F_1\e^{-x_1(r)v_{t-r}(\lambda)}\big] \bP\big[F_2\e^{-x_2(r)v_{t-r}(\lambda)}\big] \cr
 \ar=\ar
\bP\big[F_1F_2\e^{-x(r)v_{t-r}(\lambda)}\big].
 \eeqnn
A monotone class argument shows that
 \beqnn
\bP\big[F\e^{-\lambda x(t)}\big]
 =
\bP\big[F\e^{-x(r)v_{t-r}(\lambda)}\big]
 \eeqnn
for any bounded $\mcr{F}_r$-measurable random variable $F$. Then $\{(x(t),\mcr{F}_t): t\ge 0\}$ is a Markov processes with transition semigroup $(Q_t)_{t\ge 0}$. \qed

Let $D[0,\infty)$ denote the space of positive c\`{a}dl\`{a}g paths on $[0,\infty)$ furnished with the Skorokhod topology. The following theorem is a slight modification of Theorem~2.1 of Li (2006), which gives a physical interpretation of the CB-process as an approximation of the GW-process with small individuals.

\bgtheorem\label{ts2.13} Suppose that Condition~\ref{ts2.3} holds. Let $\{x(t): t\ge 0\}$ be a c\`{a}dl\`{a}g CB-process with transition semigroup $(Q_t)_{t\ge 0}$ defined by \eqref{s2.19}. For $k\ge 1$ let $\{z_k(n): n\ge 0\}$ be a Markov chain with state space $E_k := \{0,k^{-1},2k^{-1},\ldots\}$ and $n$-step transition probability $Q_k^n(x,\d y)$ determined by \eqref{s2.6}. If $z_k(0)$ converges to $x(0)$ in distribution, then $\{z_k(\lfloor\gamma_kt\rfloor): t\ge 0\}$ converges as $k\to \infty$ to $\{x(t): t\ge 0\}$ in distribution on $D[0,\infty)$. \edtheorem

\proof For $\lambda>0$ and $x\ge 0$ set $e_\lambda(x) = \e^{-\lambda x}$. Let $C_0[0,\infty)$ be the space of continuous functions on $[0,\infty)$ vanishing at infinity. By \eqref{s2.7}, \eqref{s2.19} and Theorem~\ref{ts2.7} it is easy to show
 \beqnn
\lim_{k\to\infty} \sup_{x\in E_k}\big|Q_k^{\lfloor\gamma_kt\rfloor}e_\lambda(x) -
Q_te_\lambda(x)\big| = 0, \qquad \lambda>0.
 \eeqnn
Then the Stone--Weierstrass theorem implies
 \beqnn
\lim_{k\to\infty} \sup_{x\in E_k}\big|Q_k^{\lfloor\gamma_kt\rfloor}f(x) - Q_tf(x)\big| = 0, \qquad f\in C_0[0,\infty).
 \eeqnn
By Ethier and Kurtz (1986, p.226 and pp.233--234) we conclude that $\{z_k(\lfloor\gamma_kt\rfloor): t\ge 0\}$ converges to the CB-process $\{x(t): t\ge 0\}$ in distribution on $D[0,\infty)$. \qed

For any $w\in D[0,\infty)$ let $\tau_0(w)= \inf\{s>0: w(s)=0\}$. Let $D_0[0,\infty)$ be the set of paths $w\in D[0,\infty)$ such that $w(t)=0$ for $t\ge \tau_0(w)$. Then $D_0[0,\infty)$ is a Borel subset of $D[0,\infty)$. Clearly, the distributions of the processes $\{z_k(\lfloor\gamma_kt\rfloor): t\ge 0\}$ and $\{x(t): t\ge 0\}$ are all supported by $D_0[0,\infty)$. By Theorem~1.7 of Li (2011, p.4) we have the following:

\bgcorollary\label{ts2.14} Under the conditions of Theorem~\ref{ts2.13}, the sequence $\{z_k(\lfloor\gamma_kt\rfloor): t\ge 0\}$ converges as $k\to \infty$ to $\{x(t): t\ge 0\}$ in distribution on $D_0[0,\infty)$. \edcorollary

The convergence of rescaled GW-processes to diffusion processes was first studied by Feller (1951). Lamperti (1967a) showed that all CB-processes are weak limits of rescaled GW-processes. A characterization of CB-processes by random time changes of L\'{e}vy processes was given by Lamperti (1967b); see also Kyprianou (2014). We have followed Aliev and Shchurenkov (1982) and Li (2006, 2011) in some of the above calculations.

\bgexample\label{es2.2} For any $0\le \alpha\le 1$ the function $\phi(\lambda) = \lambda^{1+\alpha}$ can be represented in the form of \eqref{s2.13}. In particular, for $0< \alpha< 1$ we can use integration by parts to see
 \beqnn
\ar\ar\int_{(0,\infty)} (\e^{-\lambda u}-1+\lambda u) \frac{\d u}{u^{2+\alpha}} \cr
 \ar\ar\qquad
= \lambda^{1+\alpha}\int_{(0,\infty)} (\e^{-v}-1+v) \frac{\d v}{v^{2+\alpha}} \cr
 \ar\ar\qquad
= \lambda^{1+\alpha}\bigg[-\frac{\e^{-v}-1+v}{(1+\alpha)v^{1+\alpha}}\bigg|_0^\infty
+ \int_{(0,\infty)} \frac{(1-\e^{-v})\d v}{(1+\alpha)v^{1+\alpha}}\bigg] \cr
 \ar\ar\qquad
= \frac{\lambda^{1+\alpha}}{1+\alpha}\bigg[-(1-\e^{-v}) \frac{1}{\alpha v^\alpha}\bigg|_0^\infty + \int_{(0,\infty)} \e^{-v} \frac{\d v}{\alpha v^\alpha}\bigg] \cr
 \ar\ar\qquad
= \frac{\Gamma(1-\alpha)}{\alpha(1+\alpha)}\lambda^{1+\alpha}.
 \eeqnn
Thus we have
 \beqlb\label{s2.21}
\lambda^{1+\alpha}
 =
\frac{\alpha(1+\alpha)}{\Gamma(1-\alpha)}\int_{(0,\infty)} (\e^{-\lambda u} -
1 + \lambda u) \frac{\d u}{u^{2+\alpha}}, \qquad \lambda\ge 0.
 \eeqlb
\edexample

\bgexample\label{es2.3} Suppose that there are constants $c>0$, $0<\alpha\le 1$ and $b$ so that $\phi(z) = cz^{1+\alpha} + bz$. Let $q^0_\alpha(t) = \alpha t$ and $q^b_\alpha(t)= b^{-1}(1-\e^{-\alpha bt})$ for $b\neq 0$. By solving the equation
 \beqnn
\frac{\partial}{\partial t}v_t(\lambda)
 =
- cv_t(\lambda)^{1+\alpha} - bv_t(\lambda),
 \qquad
v_0(\lambda) = \lambda
 \eeqnn
we get
 \beqlb\label{s2.22}
v_t(\lambda)
 =
\frac{\e^{-bt}\lambda}{\big[1 + cq^b_\alpha(t) \lambda^\alpha\big]^{1/\alpha}}, \qquad t\ge 0, \lambda\ge 0.
 \eeqlb
\edexample

\newpage

\section{Some basic properties}

 \setcounter{equation}{0}

In this section we prove some basic properties of CB-processes. Most of the results presented here can be found in Grey (1974) and Li (2000). We here use the treatments in Li (2011). Suppose that $\phi$ is a branching mechanism defined by \eqref{s2.13}. This is a convex function on $[0,\infty)$. In fact, it is easy to see that
 \beqlb\label{s3.1}
\phi^\prime(z)
 =
b + 2cz + \int_{(0,\infty)} u\big(1-\e^{-zu}\big)m(\d u), \qquad z\ge 0,
 \eeqlb
which is an increasing function. In particular, we have $\phi^\prime(0) = b$. The limit $\phi(\infty) := \lim_{z\to \infty}\phi(z)$ exists in $[-\infty,\infty]$ and $\phi^\prime(\infty) := \lim_{z\to \infty}\phi^\prime(z)$ exists in $(-\infty,\infty]$. In particular, we have
 \beqlb\label{s3.2}
\phi^\prime(\infty)
 :=
b + 2c\cdot\infty + \int_{(0,\infty)} u\, m(\d u)
 \eeqlb
with $0\cdot\infty = 0$ by convention. Observe also that $-\infty\le \phi(\infty)\le 0$ if and only if $\phi^\prime(\infty)\le 0$, and $\phi(\infty) = \infty$ if and only if $\phi^\prime(\infty)>0$.

The transition semigroup $(Q_t)_{t\ge 0}$ of the CB-process is defined by \eqref{s2.15} and \eqref{s2.19}. {From} the branching property \eqref{s2.20}, we see that the probability measure $Q_t(x,\cdot)$ is infinitely divisible. Then $(v_t)_{t\ge 0}$ has the \index{canonical representation} \textit{canonical representation}:
 \beqlb\label{s3.3}
v_t(\lambda) = h_t\lambda + \int_{(0,\infty)} (1-\e^{-\lambda u}) l_t(\d u),
\qquad t\ge 0,\lambda\ge 0,
 \eeqlb
where $h_t\ge 0$ and $l_t(\d u)$ is a $\sigma$-finite measure on $(0,\infty)$ satisfying
 \beqnn
\int_{(0,\infty)} (1\land u)l_t(\d u)< \infty.
 \eeqnn
The pair $(h_t,l_t)$ is uniquely determined by \eqref{s3.3}; see, e.g., Proposition~1.30 in Li (2011, p.16). By differentiating both sides of the equation and using \eqref{s2.15} it is easy to find
 \beqlb\label{s3.4}
h_t + \int_{(0,\infty)} ul_t(\d u)
 =
\frac{\partial}{\partial \lambda}v_t(0+)
 =
\e^{-bt}, \qquad t\ge 0.
 \eeqlb
Then we infer that $l_t(\d u)$ satisfies
 \beqnn
\int_{(0,\infty)} ul_t(\d u)< \infty.
 \eeqnn
{From} \eqref{s2.19} and \eqref{s3.4} we get
 \beqlb\label{s3.5}
\int_{[0,\infty)} y Q_t(x,\d y)
 =
x\e^{-bt}, \qquad t\ge 0,x\ge 0.
 \eeqlb
We say the branching mechanism $\phi$ is \index{critical branching mechanism} \textit{critical}, \index{subcritical branching mechanism} \textit{subcritical} or \index{supercritical branching mechanism} \textit{supercritical} according as $b=0$, $b\ge 0$ or $b\le 0$, respectively.

{From} \eqref{s2.15} we see that $t\mapsto v_t(\lambda)$ is first continuous and then continuously differentiable. Moreover, it is easy to show that
 \beqnn
\frac{\partial}{\partial t}v_t(\lambda)\Big|_{t=0}
 =
-\phi(\lambda),
 \qquad
\lambda\ge 0.
 \eeqnn
By the semigroup property $v_{t+s} = v_s\circ v_t = v_t\circ v_s$ we get the backward differential equation
 \beqlb\label{s3.6}
\frac{\partial}{\partial t}v_t(\lambda)
 =
-\phi(v_t(\lambda)),
 \qquad
v_0(\lambda) = \lambda,
 \eeqlb
and forward differential equation
 \beqlb\label{s3.7}
\frac{\partial}{\partial t}v_t(\lambda)
 =
-\phi(\lambda)\frac{\partial}{\partial \lambda}v_t(\lambda),
 \quad
v_0(\lambda)=\lambda.
 \eeqlb
The corresponding equations for a branching process with continuous time and discrete state were given in Athreya and Ney (1972, p.106).

\bgproposition\label{ts3.1} Suppose that $\lambda> 0$ and $\phi(\lambda)\neq 0$. Then the equation $\phi(z) = 0$ has no root between $\lambda$ and $v_t(\lambda)$ for every $t\ge 0$. Moreover, we have
 \beqlb\label{s3.8}
\int_{v_t(\lambda)}^\lambda\phi(z)^{-1}\d z =t, \qquad t\ge 0.
 \eeqlb
\edproposition

\proof By \eqref{s2.13} we see $\phi(0)=0$ and $z\mapsto \phi(z)$ is a
convex function. Since $\phi(\lambda)\neq 0$ for some $\lambda> 0$
according to the assumption, the equation $\phi(z) = 0$ has at most one
root in $(0,\infty)$. Suppose that $\lambda_0\ge 0$ is a root of $\phi(z)
= 0$. Then \eqref{s3.7} implies $v_t(\lambda_0) = \lambda_0$ for all
$t\ge 0$. By Proposition~\ref{ts2.10} we have $v_t(\lambda)>\lambda_0$
for $\lambda> \lambda_0$ and $0<v_t(\lambda)<\lambda_0$ for $0<\lambda<
\lambda_0$. Then $\lambda>0$ and $\phi(\lambda)\neq 0$ imply there is no
root of $\phi(z) = 0$ between $\lambda$ and $v_t(\lambda)$. {From}
\eqref{s3.6} we get \eqref{s3.8}. \qed

\bgcorollary\label{ts3.2} Suppose that $\phi(z_0)\neq 0$ for some $z_0> 0$. Let $\theta_0 = \inf\{z>0: \phi(z)\ge 0\}$ with the convention $\inf\emptyset = \infty$. Then $\lim_{t\to \infty}v_t(\lambda) = \theta_0$ increasingly for $0< \lambda< \theta_0$ and decreasingly $\lambda> \theta_0$.
\edcorollary

\proof In the case $\theta_0=\infty$, we have $\phi(z)<0$ for all $z>0$. {From} \eqref{s3.6} we see $\lambda\mapsto v_t(\lambda)$ is increasing. Then \eqref{s3.6} implies $\lim_{t\to \infty}v_t(\lambda) = \infty$ for every $\lambda>0$. In the case $\theta_0<\infty$, we have clearly $\phi(\theta_0)=0$. Furthermore, $\phi(z)<0$ for $0<z<\theta_0$ and $\phi(z)>0$ for $z>\theta_0$. {From} \eqref{s3.7} we see $v_t(\theta_0) = \theta_0$ for all $t\ge 0$. Then \eqref{s3.8} implies that $\lim_{t\to \infty}v_t(\lambda) = \theta_0$ increasingly for $0< \lambda< \theta_0$ and decreasingly $\lambda> \theta_0$. \qed

\bgcorollary\label{ts3.3} Suppose that $\phi(z_0)\neq 0$ for some $z_0> 0$. Then for any $x>0$ we have
 \beqnn
\lim_{t\to \infty}Q_t(x,\cdot) = \e^{-x\theta_0}\delta_0 + (1-\e^{-x\theta_0})\delta_\infty
 \eeqnn
by weak convergence of probability measures on $[0,\infty]$.
\edcorollary

\proof The space of probability measures on $[0,\infty]$ endowed the topology of weak convergence is compact and metrizable; see, e.g., Parthasarathy (1967, p.45). Let $\{t_n\}$ be any positive sequence so that $t_n\to \infty$ and $Q_{t_n}(x,\cdot)\to$ some $Q_{\infty}(x,\cdot)$ weakly as $n\to \infty$. By \eqref{s2.19} and Corollary~\ref{ts3.2}, for every $\lambda>0$ we have
 \beqnn
\int_{[0,\infty]}\e^{-\lambda y}Q_\infty(x,\d y)
 \ar=\ar
\lim_{n\to \infty}\int_{[0,\infty]}\e^{-\lambda y}Q_{t_n}(x,\d y) \ccr
 \ar=\ar
\lim_{n\to \infty}\e^{-xv_{t_n}(\lambda)}
 =
\e^{-x\theta_0}.
 \eeqnn
It follows that
 \beqnn
Q_\infty(x,\{0\})
 =
\lim_{\lambda\to \infty}\int_{[0,\infty]}\e^{-\lambda y}Q_\infty(x,\d y)
 =
\e^{-x\theta_0}
 \eeqnn
and
 \beqnn
Q_\infty(x,\{\infty\})
 =
\lim_{\lambda\to 0}\int_{[0,\infty]}(1-\e^{-\lambda y})Q_\infty(x,\d y)
 =
1-\e^{-x\theta_0}.
 \eeqnn
That shows $Q_\infty(x,\cdot) = \e^{-x\theta_0}\delta_0 + (1-\e^{-x\theta_0}) \delta_\infty$, which is independent of the particular choice of the sequence $\{t_n\}$. Then we have $Q_t(x,\cdot)\to Q_{\infty}(x,\cdot)$ weakly as $t\to \infty$. \qed

A simple asymptotic behavior of the CB-process is described in Corollary~\ref{ts3.3}. Clearly, we have: {\rm(i)} $\theta_0> 0$ if and only if $b< 0$; {\rm(ii)} $\theta_0= \infty$ if and only if $\phi^\prime(\infty)\le 0$. The reader can refer to Grey (1974) and Li (2011, Section~3.2) for more asymptotic results on the CB-process.

Since $(Q_t)_{t\ge 0}$ is a Feller transition semigroup, the CB-process has a Hunt process realization $X = (\Omega, \mcr{F}, \mcr{F}_t, x(t), \bQ_x)$; see, e.g., Chung (1982, p.75). Let $\tau_0 := \inf\{s\ge 0: x(s)=0\}$ denote the \index{extinction time} \textit{extinction time} of the CB-process.

\bgtheorem\label{ts3.4} For every $t\ge 0$ the limit $\bar{v}_t = \ulim_{\lambda\to \infty} v_t(\lambda)$ exists in $(0,\infty]$. Moreover, the mapping $t\mapsto \bar{v}_t$ is decreasing and for any $t\ge 0$ and $x>0$ we have
 \beqlb\label{s3.9}
\bQ_x\{\tau_0\le t\} = \bQ_x\{x(t)=0\} = \exp\{-x\bar{v}_t\}.
 \eeqlb
\edtheorem

\proof By Proposition~\ref{ts2.10} the limit $\bar{v}_t = \ulim_{\lambda\to \infty} v_t(\lambda)$ exists in $(0,\infty]$ for every $t\ge 0$. For
$t\ge r\ge 0$ we have
 \beqlb\label{s3.10}
\bar{v}_t
 =
\ulim_{\lambda\to \infty} v_r(v_{t-r}(\lambda))
 =
v_r(\bar{v}_{t-r})
 \le
\bar{v}_r.
 \eeqlb
Since zero is a trap for the CB-process, we get \eqref{s3.9} by letting $\lambda\to \infty$ in \eqref{s2.19}. \qed

For the convenience of statement of the results in the sequel, we formulate the following condition on the branching mechanism, which is known as \index{Grey's condition} \textit{Grey's condition}:

\bgcondition\label{ts3.5} There is some constant $\theta>0$ so that
 \beqnn
\phi(z)>0 ~\mbox{for}~ z\ge \theta ~\mbox{and}~ \int_\theta^\infty
\phi(z)^{-1}\d z< \infty.
 \eeqnn
\edcondition

\bgtheorem\label{ts3.6} We have $\bar{v}_t< \infty$ for some and hence all $t>0$ if and only if Condition~\ref{ts3.5} holds. \edtheorem

\proof By \eqref{s3.10} it is simple to see that $\bar{v}_t = \ulim_{\lambda\to \infty} v_t(\lambda)< \infty$ for all $t>0$ if and only if this holds for some $t>0$. If Condition~\ref{ts3.5} holds, we can let $\lambda\to \infty$ in \eqref{s3.8} to obtain
 \beqlb\label{s3.11}
\int_{\bar{v}_t}^\infty \phi(z)^{-1}\d z =t
 \eeqlb
and hence $\bar{v}_t< \infty$ for $t>0$. For the converse, suppose that $\bar{v}_t<\infty$ for some $t>0$. By \eqref{s3.6} there exists some $\theta>0$ so that $\phi(\theta)>0$, for otherwise we would have $\bar{v}_t\ge v_t(\lambda)\ge \lambda$ for all $\lambda\ge 0$, yielding a contradiction. Then $\phi(z)>0$ for all $z\ge \theta$ by the convexity of the branching mechanism. As in the
above we see that \eqref{s3.11} still holds, so Condition~\ref{ts3.5} is satisfied. \qed

\bgtheorem\label{ts3.7} Let $\bar{v} = \dlim_{t\to \infty} \bar{v}_t\in [0,\infty]$. Then for any $x>0$ we have
 \beqlb\label{s3.12}
\bQ_x\{\tau_0< \infty\} = \exp\{-x\bar{v}\}.
 \eeqlb
Moreover, we have $\bar{v}< \infty$ if and only if Condition~\ref{ts3.5} holds, and in this case $\bar{v}$ is the largest root of $\phi(z)=0$. \edtheorem

\proof The first assertion follows immediately from Theorem~\ref{ts3.4}. By Theorem~\ref{ts3.6} we have $\bar{v}_t< \infty$ for some and hence all $t>0$ if and only if Condition~\ref{ts3.5} holds. This is clearly equivalent to $\bar{v}< \infty$. {From} \eqref{s3.11} we see $\bar{v}$ is the largest root of $\phi(z)=0$. \qed

\bgcorollary\label{ts3.8} Suppose that Condition~\ref{ts3.5} holds. Then for any $x>0$ we have $\bQ_x\{\tau_0< \infty\} = 1$ if and only if $b\ge 0$. \edcorollary

By Corollary~\ref{ts3.2} and Theorem~\ref{ts3.7} we see $0\le \theta_0\le \bar{v}\le \infty$. In fact, we have $0\le \theta_0 = \bar{v}< \infty$ if Condition~\ref{ts3.5} holds and $0\le \theta_0< \bar{v}= \infty$ if there is $\theta>0$ so that
 \beqnn
\phi(z)>0 ~\mbox{for}~ z\ge \theta ~\mbox{and}~ \int_\theta^\infty \phi(z)^{-1}\d z= \infty.
 \eeqnn

\bgproposition\label{ts3.9} For any $t\ge 0$ and $\lambda\ge 0$ let $v_t^\prime(\lambda) = (\partial/\partial \lambda) v_t(\lambda)$. Then we have
 \beqlb\label{s3.13}
v_t^\prime(\lambda)
 =
\exp\bigg\{- \int_0^t \phi^\prime(v_s(\lambda))\d s\bigg\},
 \eeqlb
where $\phi^\prime$ is given by \eqref{s3.1}.
\edproposition

\proof Based on \eqref{s2.15} and \eqref{s3.6} it is elementary to see that
 \beqnn
\frac{\partial}{\partial t} v_t^\prime(\lambda)
 =
\frac{\partial}{\partial \lambda}\frac{\partial}{\partial t}
v_t(\lambda)
 =
- \phi^\prime(v_t(\lambda))v_t^\prime(\lambda).
 \eeqnn
It follows that
 \beqnn
\frac{\partial}{\partial t}\big[\log v_t^\prime(\lambda)\big]
 =
v_t^\prime(\lambda)^{-1} \frac{\partial}{\partial t} v_t^\prime(\lambda)
 =
- \phi^\prime(v_t(\lambda)).
 \eeqnn
Since $v_0^\prime(\lambda) = 1$, we get \eqref{s3.13}. \qed

\bgtheorem\label{ts3.10} Let $\phi_0^\prime(z) = \phi^\prime(z) - b$ for $z\ge 0$, where $\phi^\prime$ is given by \eqref{s3.1}. We can define a Feller transition semigroup $(Q_t^b)_{t\ge 0}$ on $[0,\infty)$ by
 \beqlb\label{s3.14}
\int_{[0,\infty)}\e^{-\lambda y}Q^b_t(x,\d y)
 =
\exp\bigg\{-xv_t(\lambda) - \int_0^t\phi^\prime_0(v_s(\lambda))\d s\bigg\}.
 \eeqlb
Moreover, we have $Q^b_t(x,\d y)= \e^{bt}x^{-1}y Q_t(x,\d y)$ for $x>0$ and
 \beqlb\label{s3.15}
Q_t^b(0,\d y) = \e^{bt}[h_t\delta_0(\d y) + yl_t(\d y)], \qquad t,y\ge 0.
 \eeqlb
\edtheorem

\proof In view of \eqref{s3.5}, it is simple to check that $Q^b_t(x,\d y) := \e^{bt}x^{-1}y Q_t(x,\d y)$ defines a Markov transition semigroup $(Q_t^b)_{t\ge 0}$ on $(0,\infty)$. Let $q_t(\lambda) = \e^{bt} v_t(\lambda)$ and let $q_t^\prime(\lambda) = (\partial/\partial \lambda) q_t(\lambda)$. By differentiating both sides of \eqref{s2.19} we see
 \beqnn
\int_{(0,\infty)}\e^{-\lambda y}Q^b_t(x,\d y)
 =
\exp\{-xv_t(\lambda)\}q_t^\prime(\lambda), \qquad x>0,\lambda\ge 0.
 \eeqnn
{From} \eqref{s3.3} and \eqref{s3.13} we have
 \beqnn
q_t^\prime(\lambda)
 =
\e^{bt}\bigg[h_t + \int_{(0,\infty)} \e^{-\lambda u}u l_t(\d u)\bigg]
 =
\exp\bigg\{- \int_0^t\phi_0^\prime(v_s(\lambda))\d s\bigg\}.
 \eeqnn
Then we can define $Q_t^b(0,\d y)$ by \eqref{s3.15} and extend $(Q^b_t)_{t\ge 0}$ to a Markov transition semigroup on $[0,\infty)$. The Feller property of the semigroup is immediate by \eqref{s3.14}. \qed

\bgcorollary\label{ts3.11} Let $(Q_t^b)_{t\ge 0}$ be the transition semigroup define by \eqref{s3.14}. Then we have $Q_t^b(0,\{0\}) = \e^{bt}h_t$ and $Q_t^b(x,\{0\})= 0$ for $t\ge 0$ and $x>0$. \edcorollary

\bgtheorem\label{ts3.12} Suppose that $T>0$ and $x>0$. Then $\bP_x^{b,T}(\d\omega)= x^{-1}\e^{bT}x(\omega,T)$ $\bQ_x(\d\omega)$ defines a probability measure on $(\Omega, \mcr{F}_T)$. Moreover, the process $\{(x(t),\mcr{F}_t): 0\le t\le T\}$ under this measure is a Markov process with transition semigroup $(Q_t^b)_{t\ge 0}$ given by \eqref{s3.14}.
\edtheorem

\proof Clearly, the probability measure $\bP_x^{b,T}$ is carried by $\{x(T)>0\}\in \mcr{F}_T$. Then we have $\bP_x^{b,T}\{x(t)>0\}= 1$ for every $0\le t\le T$. Let $0\le r\le t\le T$. Let $F$ be a bounded $\mcr{F}_r$-measurable random variable and $f$ a bounded Borel function on $[0,\infty)$. By \eqref{s3.5} and the Markov property under $\bQ_x$,
 \beqnn
\bP_x^{b,T}[Ff(x(t))]
 \ar=\ar
x^{-1}\e^{bT}\bQ_x\big[Ff(x(t))x(T)\big] \cr
 \ar=\ar
x^{-1}\e^{bt}\bQ_x\big[Ff(x(t))x(t)\big] \cr
 \ar=\ar
x^{-1}\e^{br}\bQ_x\big[Fx(r)Q_{t-r}^bf(x(r))\big] \cr
 \ar=\ar
\bP_x^{b,T}\big[FQ_{t-r}^bf(x(r))\big],
 \eeqnn
where we have used the relation $Q^b_{t-r}(x,\d y)= \e^{bt}x^{-1}y Q_{t-r}(x,\d y)$ for the third equality. Then $\{(x(t),\mcr{F}_t): 0\le t\le T\}$ under $\bP_x^{b,T}$ is a Markov process with transition semigroup $(Q_t^b)_{t\ge 0}$. \qed

Recall that zero is a trap for the CB-process. Let $(Q_t^\circ)_{t\ge 0}$ denote the restriction of its transition semigroup $(Q_t)_{t\ge 0}$ to $(0,\infty)$. For a $\sigma$-finite measure $\mu$ on $(0,\infty)$ write
 \beqnn
\mu Q_t^\circ(\d y) = \int_{(0,\infty)} \mu(\d x) Q_t^\circ(x,\d y), \qquad
t\ge 0, y> 0.
 \eeqnn
A family of $\sigma$-finite measures $(\kappa_t)_{t>0}$ on $(0,\infty)$ is called an \index{entrance rule} \textit{entrance rule} for $(Q^\circ_t)_{t\ge 0}$ if $\kappa_r Q^\circ_{t-r}\le \kappa_t$ for all $t>r>0$ and $\kappa_rQ^\circ_{t-r}\to \kappa_t$ as $r\to t$. We call $(\kappa_t)_{t>0}$ an \index{entrance law} \textit{entrance law} if $\kappa_rQ^\circ_{t-r}= \kappa_t$ for all $t>r>0$.

The special case of the canonical representation \eqref{s3.3} with $h_t=0$ for all $t>0$ is particularly interesting. In this case, we have
 \beqlb\label{s3.16}
v_t(\lambda) = \int_{(0,\infty)} (1-\e^{-\lambda u}) l_t(\d u), \qquad t>0,\lambda\ge 0.
 \eeqlb
{From} this and \eqref{s2.19} we have, for $t>0$ and $\lambda\ge 0$,
 \beqnn
\int_{(0,\infty)} (1-\e^{-y\lambda})l_t(\d y)
 =
\lim_{x\to 0} x^{-1}\int_{(0,\infty)}(1-\e^{-y\lambda})Q_t^\circ(x,\d y).
 \eeqnn
Then, formally,
 \beqlb\label{s3.17}
l_t = \lim_{x\to 0} x^{-1}Q_t(x,\cdot).
 \eeqlb

\bgtheorem\label{ts3.13} The cumulant semigroup $(v_t)_{t\ge 0}$ admits representation \eqref{s3.16} if and only if $\phi^\prime(\infty) = \infty$. In this case, the family $(l_t)_{t>0}$ is an entrance law for $(Q_t^\circ)_{t\ge 0}$. \edtheorem

\proof By differentiating both sides of the general representation \eqref{s3.3} we get
 \beqlb\label{s3.18}
v_t^\prime(\lambda) = h_t + \int_{(0,\infty)} u\e^{-\lambda u} l_t(\d u),
\qquad t\ge 0,\lambda\ge 0.
 \eeqlb
{From} this and \eqref{s3.13} it follows that
 \beqnn
h_t = v_t^\prime(\infty) = \exp\bigg\{- \int_0^t \phi^\prime(\bar{v}_s)\d
s\bigg\}.
 \eeqnn
Then we have $\phi^\prime(\infty)= \infty$ if $h_t = 0$ for any $t>0$. For the converse, assume that $\phi^\prime(\infty)= \infty$. If Condition~\ref{ts3.5} holds, we have $\bar{v}_t< \infty$ for $t>0$ by Theorem~\ref{ts3.6}, so $h_t = 0$ by \eqref{s3.3}. If Condition~\ref{ts3.5} does not hold, we have $\bar{v}_t = \infty$ by Theorem~\ref{ts3.6}. Since $\phi^\prime(\infty) = \infty$, by \eqref{s3.18} and \eqref{s3.13} we see $h_t = v_t^\prime(\infty) = 0$ for $t>0$. If $(v_t)_{t\ge 0}$ admits the representation \eqref{s3.16}, we can use \eqref{s2.18} to see, for $t>r>0$ and $\lambda\ge 0$,
 \beqnn
\int_{(0,\infty)} (1-\e^{-\lambda u}) l_t(\d u)
 \ar=\ar
\int_{(0,\infty)} (1-\e^{- uv_{t-r}(\lambda)}) l_r(\d u) \cr
 \ar=\ar
\int_{(0,\infty)} l_r(\d x)\int_{(0,\infty)} (1-\e^{-\lambda u}) Q_{t-r}^\circ(x,\d
u).
 \eeqnn
Then $(l_t)_{t>0}$ is an entrance law for $(Q_t^\circ)_{t\ge 0}$. \qed

\bgcorollary\label{ts3.14} If Condition~\ref{ts3.5} holds, the cumulant semigroup admits the representation \eqref{s3.16} and $t\mapsto \bar{v}_t = l_t(0,\infty)$ is the unique solution to the differential equation
 \beqlb\label{s3.19}
\frac{\d}{\d t}\bar{v}_t = -\phi(\bar{v}_t), \qquad t>0
 \eeqlb
with singular initial condition $\bar{v}_{0+} = \infty$. \edcorollary

\proof Under Condition~\ref{ts3.5}, for every $t>0$ we have $\bar{v}_t< \infty$ by Theorem~\ref{ts3.6}. Moreover, the condition and the convexity of $z\mapsto \phi(z)$ imply $\phi^\prime(\infty) = \infty$. Then we have the representation \eqref{s3.16} by Theorem~\ref{ts3.13}. The semigroup property of $(v_t)_{t\ge 0}$ implies $\bar{v}_{s+t} = v_s(\bar{v}_t)$ for $s>0$ and $t>0$. Then $t\mapsto \bar{v}_t$ satisfies \eqref{s3.19}. {From} \eqref{s3.11} it is easy to see $\bar{v}_{0+} = \infty$. Suppose that $t\mapsto u_t$ and $t\mapsto v_t$ are two solutions to \eqref{s3.19} with $u_{0+} = v_{0+} = \infty$. For any $\varepsilon>0$ there exits $\delta>0$ so that $u_s\ge v_\varepsilon$ for every $0<s\le \delta$. Since both $t\mapsto u_{s+t}$ and $t\mapsto v_{\varepsilon+t}$ are solutions to \eqref{s3.19}, we have $u_{s+t}\ge v_{\varepsilon+t}$ for $t\ge 0$ and $0<s\le \delta$ by Proposition~\ref{ts2.10}. Then we can let $s\to 0$ and $\varepsilon\to 0$ to see $u_t\ge v_t$ for $t>0$. By symmetry we get the uniqueness of the solution. \qed

\bgtheorem\label{ts3.15} If $\delta:= \phi^\prime(\infty)< \infty$, then we have, for $t\ge 0$ and $\lambda\ge 0$,
 \beqlb\label{s3.20}
v_t(\lambda) = \e^{-\delta t}\lambda + \int_0^t\e^{-\delta s}\d s\int_{(0,\infty)}(1-\e^{-u v_{t-s}(\lambda)}) m(\d u),
 \eeqlb
that is, we have \eqref{s3.3} with
 \beqlb\label{s3.21}
h_t=\e^{-\delta t},
 \quad
l_t = \int_0^t \e^{-\delta s}mQ_{t-s}^\circ\d s, \quad t\ge 0.
 \eeqlb
In this case, the family $(l_t)_{t>0}$ is an entrance rule for $(Q_t^\circ)_{t\ge 0}$.
\edtheorem

\proof If $\delta:= \phi^\prime(\infty)< \infty$, by \eqref{s3.2} we must have $c=0$. In this case, we can write the branching mechanism into
 \beqlb\label{s3.22}
\phi(\lambda) = \delta\lambda + \int_{(0,\infty)} (\e^{-\lambda z}-1)m(\d z), \qquad \lambda\ge 0.
 \eeqlb
By \eqref{s2.15} and integration by parts,
 \beqnn
v_t(\lambda)\e^{\delta t} \ar=\ar \lambda + \int_0^t \delta v_s(\lambda) \e^{\delta s} \d s - \int_0^t \phi(v_s(\lambda))\e^{\delta s} \d s \cr
 \ar=\ar
\lambda + \int_0^t\e^{\delta s}\d s\int_{(0,\infty)}(1-\e^{-u v_s(\lambda)}) m(\d u).
 \eeqnn
That gives \eqref{s3.20} and \eqref{s3.21}. It is easy to see that $(l_t)_{t>0}$ is an entrance rule for $(Q_t^\circ)_{t\ge 0}$. \qed

\bgexample\label{es3.1} Suppose that there are constants $c>0$, $0<\alpha\le 1$ and $b$ so that $\phi(z) = cz^{1+\alpha} + bz$. Then Condition~\ref{ts3.5} is satisfied. Let $q^0_\alpha(t)$ be defined as in Example~\ref{es2.3}. By letting $\lambda\to \infty$ in \eqref{s2.22} we get $\bar{v}_t = c^{-1/\alpha}\e^{-bt} q^b_\alpha(t)^{-1/\alpha}$ for $t>0$. In particular, if $\alpha = 1$, then \eqref{s3.16} holds with
 \beqnn
l_t(\d u)
 =
\frac{\e^{-bt}}{c^2q^b_1(t)^2}\exp\bigg\{-\frac{u}{cq^b_1(t)}\bigg\}\d u,
\qquad t>0, u>0.
 \eeqnn
\edexample

\newpage

\section{Positive integral functionals}

 \setcounter{equation}{0}

In this section, we give characterizations of a class of positive integral functionals of the CB-process in terms of Laplace transforms. The corresponding results in the measure-valued setting can be found in Li (2011). For our purpose, it is more convenient to start the process from an arbitrary initial time $r\ge 0$. Let $X = (\Omega, \mcr{F}, \mcr{F}_{r,t}, x(t), \bQ_{r,x})$ a c\`{a}dl\`{a}g realization of the CB-process with transition semigroup $(Q_t)_{t\ge 0}$ defined by \eqref{s2.15} and \eqref{s2.19}. For any $t\ge r\ge 0$ and $\lambda\ge 0$ we have
 \beqlb\label{s4.1}
\bQ_{r,x}\exp\{-\lambda x(t)\} = \exp\{-xu_r(\lambda)\},
 \eeqlb
where $r\mapsto u_r(\lambda):=v_{t-r}(\lambda)$ is the unique bounded positive solution to
 \beqlb\label{s4.2}
u_r(\lambda) + \int_r^t \phi(u_s(\lambda))\d s
 =
\lambda, \qquad 0\le r\le t.
 \eeqlb

\bgproposition\label{ts4.1} For $\{t_1< \cdots< t_n\}\subset [0,\infty)$ and $\{\lambda_1, \ldots, \lambda_n\}\subset [0,\infty)$ we have
 \beqlb\label{s4.3}
\bQ_{r,x}\exp\bigg\{-\sum_{j=1}^n \lambda_jx(t_j)1_{\{r\le t_j\}}\bigg\}
 =
\exp\{-xu(r)\}, \quad 0\le r\le t_n,
 \eeqlb
where $r\mapsto u(r)$ is a bounded positive function on $[0,t_n]$ solving
 \beqlb\label{s4.4}
u(r) + \int_r^{t_n} \phi(u(s)) \d s
 =
\sum_{j=1}^n \lambda_j1_{\{r\le t_j\}}.
 \eeqlb
\edproposition

\proof We shall give the proof by induction in $n\ge 1$. For $n=1$ the result
follows from \eqref{s4.1} and \eqref{s4.2}. Now supposing
\eqref{s4.3} and \eqref{s4.4} are satisfied when $n$ is replaced
by $n-1$, we prove they are also true for $n$. It is clearly sufficient
to consider the case with $0\le r\le t_1< \cdots< t_n$. By the Markov
property,
 \beqnn
\ar\ar\bQ_{r,x}\exp\bigg\{-\sum_{j=1}^n \lambda_jx(t_j)\bigg\} \cr
 \ar\ar\qquad
= \bQ_{r,x}\bigg[\bQ_{r,x}\bigg(\exp\bigg\{-\sum_{j=1}^n \lambda_jx(t_j)\bigg\}\Big| \mcr{F}_{r,t_1}\bigg)\bigg] \cr
 \ar\ar\qquad
= \bQ_{r,x}\bigg[\e^{-x(t_1)\lambda_1}\bQ_{r,x}\bigg(\exp\bigg\{-\sum_{j=2}^n \lambda_jx(t_j)\bigg\}\Big| \mcr{F}_{r,t_1}\bigg)\bigg] \cr
 \ar\ar\qquad
= \bQ_{r,x}\bigg[\e^{-x(t_1)\lambda_1}\bQ_{t_1,x(t_1)}\bigg(\exp\bigg\{-\sum_{j=2}^n \lambda_jx(t_j)\bigg\}\bigg)\bigg] \cr
 \ar\ar\qquad
= \bQ_{r,x}\exp\Big\{-x(t_1)\lambda_1 - x(t_1)w(t_1)\Big\},
 \eeqnn
where $r\mapsto w(r)$ is a bounded positive Borel function on $[0,t_n]$
satisfying
 \beqlb\label{s4.5}
w(r) + \int_r^{t_n} \phi(w(s)) \d s
 =
\sum_{j=2}^n \lambda_j1_{\{r\le t_j\}}.
 \eeqlb
Then the result for $n=1$ implies that
 \beqnn
\bQ_{r,x}\exp\bigg\{-\sum_{j=1}^n \lambda_jx(t_j)\bigg\}
 =
\exp\{-xu(r)\}
 \eeqnn
with $r\mapsto u(r)$ being a bounded positive Borel function on $[0,t_1]$
satisfying
 \beqlb\label{s4.6}
u(r) + \int_r^{t_1} \phi(u(s))\d s
 =
\lambda_1 + w(t_1).
 \eeqlb
Setting $u(r) = w(r)$ for $t_1<r\le t_n$, from \eqref{s4.5} and
\eqref{s4.6} one checks that $r\mapsto u(r)$ is a bounded positive
solution to \eqref{s4.4} on $[0,t_n]$. \qed

\bgtheorem\label{ts4.2} Suppose that $t\ge 0$ and $\mu$ is a finite measure supported by $[0,t]$. Let $s\mapsto \lambda(s)$ be a bounded positive Borel function on
$[0,t]$. Then we have
 \beqlb\label{s4.7}
\bQ_{r,x}\exp\bigg\{-\int_{[r,t]}\lambda(s)x(s)\mu(\d s)\bigg\}
 =
\exp\{-xu(r)\}, \quad 0\le r\le t,
 \eeqlb
where $r\mapsto u(r)$ is the unique bounded positive solution on $[0,t]$ to
 \beqlb\label{s4.8}
u(r) + \int_r^t \phi(u(s))\d s
 =
\int_{[r,t]}\lambda(s)\mu(\d s).
 \eeqlb
\edtheorem

\proof \textit{Step~1.} We first consider a bounded positive continuous function $s\mapsto \lambda(s)$ on $[0,t]$. To avoid triviality we assume $t>0$. For any integer $n\ge1$ define the finite measure $\mu_n$ on $[0,t]$ by
 \beqnn
\mu_n(\d s)= \sum_{k=1}^{2^n} \mu[(k-1)t/2^n, kt/2^n)\delta_{kt/2^n}(\d s) + \mu(\{t\})\delta_t(\d s).
 \eeqnn
By Proposition~\ref{ts4.1} we see that
 \beqlb\label{s4.9}
\bQ_{r,x}\exp\bigg\{-\int_{[r,t]} \lambda(s)x(s) \mu_n(\d s)\bigg\}
 =
\exp\{-xu_n(r)\},
 \eeqlb
where $r\mapsto u_n(r)$ is a bounded positive solution on $[0,t]$ to
 \beqlb\label{s4.10}
u_n(r) + \int_r^t \phi(u_n(s))\d s
 =
\int_{[r,t]}\lambda(s)\mu_n(\d s).
 \eeqlb
Let $v_n(r)= u_n(t-r)$ for $0\le r\le t$. Observe that $\phi(z)\ge bz\ge -|b|z$ for every $z\ge 0$. {From} \eqref{s4.10} we have
 \beqnn
v_n(r)
 \ar=\ar
\int_{[t-r,t]}\lambda(s)\mu_n(\d s) - \int_0^r \phi(v_n(s))\d s \cr
 \ar\le\ar
\sup_{0\le s\le t}\lambda(s)\mu[0,t] + |b|\int_0^r v_n(s)\d s.
 \eeqnn
By Gronwall's inequality it is easy to show that $\{v_n\}$ and hence $\{u_n\}$ is uniformly bounded on $[0,t]$. Let $q_n(t)= t$. For any $0\le s< t$ let $q_n(s)= (\lfloor 2^ns/t\rfloor+1)t/2^n$, where $\lfloor 2^ns/t\rfloor$ denotes the integer part of $2^ns/t$. Then $s\le q_n(s)\le s+t/2^n$. It is easy to see that
 \beqnn
\int_{[r,t]} f(s)\mu_n(\d s)= \int_{[r,t]} f(q_n(s))\mu(\d s)
 \eeqnn
for any bounded Borel function $f$ on $[0,t]$. By the right-continuity of $s\mapsto \lambda(s)$ and $s\mapsto x(s)$ we have
 \beqnn
\lim_{n\to \infty}\int_{[r,t]}\lambda(s)\mu_n(\d s)
 =
\int_{[r,t]}\lambda(s)\mu(\d s)
 \eeqnn
and
 \beqnn
\lim_{n\to \infty}\int_{[r,t]} \lambda(s)x(s) \mu_n(\d s)
 =
\int_{[r,t]} \lambda(s)x(s) \mu(\d s).
 \eeqnn
{From} \eqref{s4.9} we see the limit $u(r) = \lim_{n\to \infty} u_n(r)$ exists and \eqref{s4.7} holds for $0\le r\le t$. Then we get \eqref{s4.8} by letting $n\to \infty$ in \eqref{s4.10}.

\textit{Step~2.} Let $B_0[0,\infty)$ be the set of bounded Borel functions $s\mapsto \lambda(s)$ for which there exist bounded positive solutions $r\mapsto u(r)$ of \eqref{s4.8} such that \eqref{s4.7} holds. Then $B_0[0,\infty)$ is closed under bounded pointwise convergence. The result of the first step shows that $B_0[0,\infty)$ contains all positive continuous functions on $[0,t]$. By Proposition~1.3 in Li (2011, p.3) we infer that $B_0[0,\infty)$ contains all bounded positive Borel functions on $[0,t]$.

\textit{Step~3.} To show the uniqueness of the solution to \eqref{s4.8}, suppose that $r\mapsto v(r)$ is another bounded positive Borel function on $[0,t]$ satisfying this equation. Since $z\mapsto \phi(z)$ is locally Lipschitz, it is easy to find a constant $K\ge 0$ such that
 \beqnn
|u(r)-v(r)|\ar\le\ar \int_r^t|\phi(u(s))-\phi(v(s))|\d s \cr
 \ar\le\ar
K\int_r^t|u(s)-v(s)|\d s.
 \eeqnn
Let $U(r)= |u(t-r)-v(t-r)|$ for $0\le r\le t$. We have
 \beqnn
U(r)\le K\int_0^rU(s)\d s, \quad 0\le r\le t.
 \eeqnn
Then Gronwall's inequality implies $U(r) =0$ for every $0\le
r\le t$. \qed

Suppose that $\mu(\d s)$ is a locally bounded Borel measure on $[0,\infty)$ and $s\mapsto \lambda(s)$ is a locally bounded positive Borel function on $[0,\infty)$. We define the \index{positive integral functional} \textit{positive integral functional}:
 \beqnn
A[r,t] := \int_{[r,t]} \lambda(s)x(s)\mu(\d s), \qquad t\ge r\ge 0.
 \eeqnn
By replacing $\lambda(s)$ with $\theta \lambda(s)$ in Theorem~\ref{ts4.2} for $\theta\ge 0$ we get a characterization of the Laplace transform of the random variable $A[r,t]$.

\bgtheorem\label{ts4.3} Let $t\ge 0$ be given. Let $\lambda\ge 0$ and let $s\mapsto \theta(s)$ be a bounded positive Borel function on $[0,t]$. Then for $0\le r\le t$ we have
 \beqlb\label{s4.11}
\bQ_{r,x}\exp\bigg\{-\lambda x(t) - \int_r^t \theta(s)x(s)\d s\bigg\}
 =
\exp\{-xu(r)\},
 \eeqlb
where $r\mapsto u(r)$ is the unique bounded positive solution on $[0,t]$ to
 \beqlb\label{s4.12}
u(r) + \int_r^t \phi(u(s))\d s
 =
\lambda + \int_r^t \theta(s)\d s.
 \eeqlb
\edtheorem

\proof This follows by an application of Theorem~\ref{ts4.2} to the measure $\mu(\d s) =\d s + \delta_t(\d s)$ and the function $\lambda(s) = 1_{\{s<t\}} \theta(s) + 1_{\{s=t\}}\lambda$. \qed

\bgcorollary\label{ts4.4} Let $X = (\Omega, \mcr{F}, \mcr{F}_t, x(t), \bQ_x)$ be a Hunt realization of the CB-process started from time zero. Then we have, for $t, \lambda, \theta\ge 0$,
 \beqlb\label{s4.13}
\bQ_x\exp\bigg\{-\lambda x(t) - \theta\int_0^t x(s)\d s\bigg\}
 =
\exp\{-xv(t)\},
 \eeqlb
where $t\mapsto v(t)= v(t,\lambda,\theta)$ is the unique positive solution to
 \beqlb\label{s4.14}
\frac{\partial}{\partial t} v(t) = \theta - \phi(v(t)),
 \quad
v(0) = \lambda.
 \eeqlb
\edcorollary

\proof By Theorem~\ref{ts4.3} we have \eqref{s4.13} with $v(t)= u_t(0)$, where $r\mapsto u_t(r)$ is the unique bounded positive solution on $[0,t]$ to
 \beqnn
u(r) + \int_r^t \phi(u(s))\d s
 =
\lambda + (t-r)\theta.
 \eeqnn
Then $r\mapsto v(r):= u_t(t-r)$ is the unique bounded positive solution on $[0,t]$ of
 \beqlb\label{s4.15}
v(r) + \int_0^r\phi(v(s))\d s
 =
\lambda + r\theta.
 \eeqlb
Clearly, we can extend \eqref{s4.15} to all $r\ge 0$ and the extended equation is equivalent with the differential equation \eqref{s4.14}. The uniqueness of the solution follows by Gronwall's inequality. \qed

\bgcorollary\label{ts4.5} Let $X = (\Omega, \mcr{F}, \mcr{F}_t, x(t), \bQ_x)$ be a Hunt realization of the CB-process started from time zero. Then we have, for $t, \theta\ge 0$,
 \beqlb\label{s4.16}
\bQ_x\exp\bigg\{- \theta\int_0^t x(s)\d s\bigg\}
 =
\exp\{-xv(t)\}.
 \eeqlb
where $t\mapsto v(t)=v(t,\theta)$ is the unique positive solution to
 \beqlb\label{s4.17}
\frac{\partial}{\partial t} v(t) = \theta - \phi(v(t)),
 \quad
v(0) = 0.
 \eeqlb
\edcorollary

\bgcorollary\label{ts4.6} Let $t\ge 0$ be given. Let $\phi_1$ and $\phi_2$ be two branching mechanisms in form \eqref{s2.13} satisfying $\phi_1(z)\ge \phi_2(z)$ for all $z\ge 0$. Let $t\mapsto v_i(t)$ be the solution to \eqref{s4.14} or \eqref{s4.15} with $\phi=\phi_i$. Then $v_1(t)\le v_2(t)$ for all $t\ge 0$.
\edcorollary

\proof Fix $t\ge 0$ and let $u_i(r) = v_i(t-r)$ for $0\le r\le t$. Then $r\mapsto u_1(r)$ is the unique bounded positive solution on $[0,t]$ of
 \beqnn
u(r) + \int_r^t\phi_1(u(s))\d s = \lambda + (t-r)\theta
 \eeqnn
and $r\mapsto u_2(r)$ is the unique bounded positive solution on $[0,t]$ of
 \beqnn
u(r) + \int_r^t\phi_1(u(s))\d s
 =
\lambda + \int_r^t [\theta + g(s)]\d s,
 \eeqnn
where $g(s) = \phi_1(u_2(s)) - \phi_2(u_2(s))\ge 0$. By Theorem~\ref{ts4.3} one can see $u_1(r)\le u_2(r)$ for all $0\le r\le t$. \qed

Recall that $\phi^\prime(\infty)$ is given by \eqref{s3.2}. Under the condition $\phi^\prime(\infty)> 0$, we have $\phi(z)\to \infty$ as $z\to \infty$, so the inverse $\phi^{-1}(\theta):= \inf\{z\ge 0: \phi(z)> \theta\}$ is well-defined for $\theta\ge 0$.

\bgproposition\label{ts4.7} For $\theta>0$ let $t\mapsto v(t,\theta)$ be the unique positive solution to \eqref{s4.17}. Then $\lim_{t\to \infty} v(t,\theta) = \infty$ if $\phi^\prime(\infty)\le 0$, and $\lim_{t\to \infty} v(t,\theta) = \phi^{-1}(\theta)$ if $\phi^\prime(\infty)> 0$. \edproposition

\proof By Proposition~\ref{ts2.10} we have $\bQ_x\{x(t)> 0\}> 0$ for every $x>0$ and $t\ge 0$. {From} \eqref{s4.16} we see $t\mapsto v(t,\theta)$ is strictly increasing, so $(\partial/\partial t) v(t,\theta)> 0$ for all $\theta> 0$. Let $v(\infty,\theta) = \lim_{t\to \infty} v(t,\theta)\in (0,\infty]$. In the case $\phi^\prime(\infty)\le 0$, we clearly have $\phi(z)\le 0$ for all $z\ge 0$. Then $(\partial/\partial t) v(t,\theta)\ge \theta>0$ and $v(\infty,\theta) = \infty$. In the case $\phi^\prime(\infty)> 0$, we note
 \beqnn
\phi(v(t,\theta)) = \theta - {\partial\over \partial t} v(t,\theta)< \theta,
 \eeqnn
and hence $v(t,\theta)< \phi^{-1}(\theta)$, implying $v(\infty,\theta)\le \phi^{-1}(\theta)< \infty$. It follows that
 \beqnn
0 = \lim_{t\to \infty}\frac{\partial}{\partial t} v(t,\theta)
 =
\theta - \lim_{t\to \infty}\phi(v(t,\theta))
 =
\theta - \phi(v(\infty,\theta)).
 \eeqnn
Then we have $v(\infty,\theta) = \phi^{-1}(\theta)$. \qed

\bgtheorem\label{ts4.8} Let $X = (\Omega, \mcr{F}, \mcr{F}_t, x(t), \bQ_x)$ be a Hunt realization of the CB-process started from time zero. If $\phi^\prime(\infty)> 0$, then for $x>0$ and $\theta>0$ we have
 \beqnn
\bQ_x\exp\bigg\{- \theta\int_0^\infty x(s)\d s\bigg\}
 =
\exp\{-x\phi^{-1}(\theta)\}
 \eeqnn
and
 \beqnn
\bQ_x\bigg\{\int_0^\infty x(s)\d s< \infty\bigg\}
 =
\exp\{-x\phi^{-1}(0)\},
 \eeqnn
where $\phi^{-1}(0)= \inf\{z>0: \phi(z)\ge 0\}$. If $\phi^\prime(\infty)\le 0$, then for any $x>0$ we have
 \beqnn
\bQ_x\bigg\{\int_0^\infty x(s)\d s< \infty\bigg\} = 0.
 \eeqnn
\edtheorem

\proof In view of \eqref{s4.16}, we have
 \beqnn
\bQ_x\exp\bigg\{- \theta\int_0^\infty x(s)\d s\bigg\}
 =
\lim_{t\to \infty}\exp\{-xv(t,\theta)\}.
 \eeqnn
Then the result follows from Proposition~\ref{ts4.7}. \qed

\newpage

\section{Construction of CBI-processes}

 \setcounter{equation}{0}

Let $\{p(j): j\in \mbb{N}\}$ and $\{q(j): j\in \mbb{N}\}$ be probability distributions on $\mbb{N}:= \{0,1,2,\ldots\}$ with generating functions $g$ and $h$, respectively. Suppose that $\{\xi_{n,i}: n,i=1,2,\ldots\}$ is a family of $\mbb{N}$-valued i.i.d.\ random variables with distribution $\{p(j): j\in \mbb{N}\}$ and $\{\eta_n: n=1,2,\ldots\}$ is a family of $\mbb{N}$-valued i.i.d.\ random variables with distribution $\{q(j): j\in \mbb{N}\}$. We assume the two families are independent of each other. Given an $\mbb{N}$-valued random variable $y(0)$ independent of $\{\xi_{n,i}\}$ and $\{\eta_n\}$, we define inductively
 \beqlb\label{s5.1}
y(n)= \sum_{i=1}^{y(n-1)}\xi_{n,i} + \eta_n, \qquad n=1,2,\ldots.
 \eeqlb
This is clearly a generalization of \eqref{s2.1}. For $i\in \mbb{N}$ let $\{Q(i,j): j\in \mbb{N}\}$ denote the $i$-fold convolution of $\{p(j): j\in \mbb{N}\}$. Let
 \beqnn
P(i,j) = (Q(i,\cdot)*q)(j) = (p^{*i}*q)(j), \qquad i,j\in \mbb{N}.
 \eeqnn
For any $n\ge 1$ and $\{i_0, \cdots, i_{n-1}=i,j\}\subset \mbb{N}$ we have
 \beqnn
\ar\ar\bP\Big(y(n) = j\big|y(0)=i_0, y(1)=i_1, \cdots, y(n-1)=i_{n-1}\Big) \cr
 \ar\ar\qquad
= \bP\bigg(\sum_{k=1}^{y(n-1)}\xi_{n,k} + \eta_n\Big|y(n-1)=i\bigg) \cr
 \ar\ar\qquad
= \bP\bigg(\sum_{k=1}^i\xi_{n,k} + \eta_n = j\bigg) = P(i,j).
 \eeqnn
Then $\{y(n): n\ge 0\}$ is a Markov chain with one-step transition matrix $P= (P(i,j): i, j\in \mbb{N})$. The random variable $y(n)$ can be thought of as the number of individuals in generation $n$ of a population system with immigration. After one unit time, each of the $y(n)$ individuals splits independently of others into a random number of offspring according to the distribution $\{p(j): j\in \mbb{N}\}$ and a random number of immigrants are added to the system according to the distribution $\{q(j): j\in \mbb{N}\}$. It is easy to see that
 \beqlb\label{s5.2}
\sum^\infty_{j=0} P(i,j)z^j = g(z)^ih(z), \qquad |z|\le 1.
 \eeqlb

A Markov chain in $\mbb{N}$ with one-step transition matrix defined by \eqref{s5.2} is called a \index{Galton--Watson branching process with immigration} \textit{Galton--Watson branching process with immigration} \index{GWI-process} (GWI-process) or a \index{Bienaym\'{e}--Galton--Watson branching process with immigration} \textit{Bienaym\'{e}--Galton--Watson branching process with immigration} \index{BGWI-process} (BGWI-process) with \index{branching distribution} \textit{branching distribution} given by $g$ and \index{immigration distribution} \textit{immigration distribution} given by $h$. When $h\equiv 1$, this reduces to the GW-process defined before. For any $n\ge 1$ the $n$-step transition matrix of the GWI-process is just the $n$-fold product $P^n= (P^n(i,j): i, j\in \mbb{N})$.

\bgproposition\label{ts5.1} For any $n\ge 1$ and $i\in \mbb{N}$ we have
 \beqlb\label{s5.3}
\sum^\infty_{j=0} P^n(i,j)z^j = g^{\circ n}(z)^i\prod_{j=1}^nh(g^{\circ(j-1)}(z)),
\qquad |z|\le 1.
 \eeqlb
\edproposition

\proof {From} \eqref{s5.2} we see \eqref{s5.3} holds for $n=1$. Now suppose that \eqref{s5.3} holds for some $n\ge 1$. We have
 \beqnn
\sum^\infty_{j=0} P^{n+1}(i,j)z^j
 \ar=\ar
\sum^\infty_{j=0} \sum^\infty_{k=0} P(i,k)P^n(k,j)z^j \cr
 \ar=\ar
\sum^\infty_{k=0} P(i,k)g^{\circ n}(z)^k\prod_{j=1}^nh(g^{\circ{j-1}}(z)) \cr
 \ar=\ar
g(g^{\circ n}(z))^ih(g^{\circ n}(z))\prod_{j=1}^nh(g^{\circ{j-1}}(z)) \cr
 \ar=\ar
g^{\circ(n+1)}(z)^i\prod_{j=1}^{n+1}h(g^{\circ{j-1}}(z)).
 \eeqnn
Then \eqref{s5.3} also holds when $n$ is replaced by $n+1$. That gives the result by induction.
\qed

Suppose that for each integer $k\ge 1$ we have a GWI-process $\{y_k(n): n\ge 0\}$ with branching distribution given by the probability generating function $g_k$ and immigration distribution given by the probability generating function $h_k$. Let $z_k(n) = y_k(n)/k$. Then $\{z_k(n): n\ge 0\}$ is a Markov chain with state space $E_k := \{0,1/k,2/k,\ldots\}$ and $n$-step transition probability $P_k^n(x,\d y)$
determined by
 \beqlb\label{s5.4}
\int_{E_k}\e^{-\lambda y}P_k^n(x,\d y)
 =
g_k^{\circ n}(\e^{-\lambda/k})^{kx}\prod_{j=1}^nh_k(g_k^{\circ(j-1)}(\e^{-\lambda/k})).
 \eeqlb
Suppose that $\{\gamma_k\}$ is a positive real sequence so that $\gamma_k\to \infty$ increasingly as $k\to \infty$. Let $\lfloor\gamma_kt\rfloor$ denote the integer part of $\gamma_kt$. In view of \eqref{s5.4}, given $z_k(0) = x\in E_k$, the random variable $z_k(\lfloor\gamma_kt\rfloor) = k^{-1}y_k(\lfloor\gamma_kt\rfloor)$ has distribution $P_k^{\lfloor\gamma_kt\rfloor}(x,\cdot)$ on $E_k$ determined
by
 \beqlb\label{s5.5}
\ar\ar\int_{E_k}\e^{-\lambda y}P_k^{\lfloor\gamma_kt\rfloor}(x,\d y) \cr
 \ar\ar\qquad
= g_k^{\circ\lfloor\gamma_kt\rfloor}(\e^{-\lambda/k})^{kx} \prod_{j=1}^{\lfloor\gamma_kt\rfloor} h_k(g_k^{j-1}(\e^{-\lambda/k})) \cr
 \ar\ar\qquad
= \exp\big\{x k\log g_k^{\circ\lfloor\gamma_kt\rfloor}(\e^{-\lambda/k})\big\} \exp\bigg\{\sum_{j=1}^{\lfloor\gamma_kt\rfloor} \log h_k(g_k^{j-1}(\e^{-\lambda/k}))\bigg\} \cr
 \ar\ar\qquad
= \exp\bigg\{-xv_k(t,\lambda) - \int_0^{\gamma_k^{-1}\lfloor\gamma_kt\rfloor}
\bar{\psi}_k(v_k(s,\lambda))\d s\bigg\},
 \eeqlb
where $v_k(t,\lambda)$ is given by \eqref{s2.8} and
 \beqlb\label{s5.6}
\bar{\psi}_k(z)= -\gamma_k\log h_k(\e^{-z/k}).
 \eeqlb
For any $z\ge 0$ we have
 \beqlb\label{s5.7}
\bar{\psi}_k(z)
 =
-\gamma_k\log\big[1-\gamma_k^{-1}\tilde{\psi}_k(z)\big],
 \eeqlb
where
 \beqlb\label{s5.8}
\tilde{\psi}_k(z) = \gamma_k[1-h_k(\e^{-z/k})].
 \eeqlb

\bglemma\label{ts5.2} Suppose that the sequence $\{\tilde{\psi}_k\}$ is uniformly bounded on each bounded interval. Then we have $\lim_{k\to\infty} |\bar{\psi}_k(z)-\tilde{\psi}_k(z)| = 0$ uniformly on each bounded interval. \edlemma

\proof This is immediate by the relation \eqref{s5.7}. \qed

\bgcondition\label{ts5.3} There is a function $\psi$ on $[0,\infty)$ such that $\tilde{\psi}_k(z)\to \psi(z)$ uniformly on $[0,a]$ for every $a\ge 0$ as $k\to \infty$. \edcondition

\bgproposition\label{ts5.4} Suppose that Condition~\ref{ts5.3} is satisfied. Then the limit function $\psi$ has representation
 \beqlb\label{s5.9}
\psi(z) = \beta z + \int_{(0,\infty)}\big(1-\e^{-zu}\big)\nu(\d u), \qquad z\ge
0,
 \eeqlb
where $\beta\ge 0$ is a constant and $\nu(\d u)$ is a $\sigma$-finite measure on $(0,\infty)$ satisfying
 \beqnn
\int_{(0,\infty)} (1\land u)\nu(\d u)< \infty.
 \eeqnn
\edproposition

\proof It is well-known that $\psi$ has representation \eqref{s5.9} if and only if $\e^{-\psi} = L_\mu$ is the Laplace transform of an infinitely divisible probability distribution $\mu$ on $[0,\infty)$; see, e.g., Theorem~1.39 in Li (2011, p.20). In view of \eqref{s5.8}, the function $\tilde{\psi}_k$ can be represented by a special form \eqref{s5.9}, so $\e^{-\tilde{\psi}_k} = L_{\mu_k}$ is the Laplace transform of an infinitely divisible distribution $\mu_k$ on $[0,\infty)$. By Lemma~\ref{ts5.2} and Condition~\ref{ts5.3} we have $\tilde{\psi}_k(z)\to \psi(z)$ uniformly on $[0,a]$ for every $a\ge 0$ as $k\to \infty$. By Theorem~\ref{ts1.2} there is a probability distribution $\mu$ on $[0,\infty)$ so that $\mu= \lim_{k\to \infty}\mu_k$ weakly and $\e^{-\psi} = L_\mu$. Clearly $\mu$ is also infinitely divisible, so $\psi$ has representation \eqref{s5.9}. \qed

\bgproposition\label{ts5.5} For any function $\psi$ with representation \eqref{s5.9} there is a sequence $\{\tilde{\psi}_k\}$ in the form of \eqref{s5.6} satisfying Condition~\ref{ts5.3}. \edproposition

\proof This is similar to the proof of Proposition~\ref{ts2.5} and is left to the reader as an exercise. \qed

\bgtheorem\label{ts5.6} Suppose that $\phi$ and $\psi$ are given by \eqref{s2.13} and \eqref{s5.9}, respectively. For any $\lambda\ge 0$ let $t\mapsto v_t(\lambda)$ be the unique positive solution to \eqref{s2.15}. Then there is a Feller transition semigroup $(P_t)_{t\ge 0}$ on $[0,\infty)$ defined by
 \beqlb\label{s5.10}
\int_{[0,\infty)} \e^{-\lambda y} P_t(x,\d y)
 =
\exp\bigg\{-xv_t(\lambda)-\int_0^t\psi(v_s(\lambda))\d s\bigg\}.
 \eeqlb
\edtheorem

\proof
This follows by arguments similar to those in Section~2.
\qed

If a Markov process in $[0,\infty)$ has transition semigroup $(P_t)_{t\ge 0}$ defined by \eqref{s5.10}, we call it a \index{continuous-state branching process with immigration} \textit{continuous-state branching process with immigration} \index{CBI-process} (CBI-process) with \index{branching mechanism} \textit{branching mechanism} $\phi$ and \index{immigration mechanism} \textit{immigration mechanism} $\psi$. In particular, if
 \beqlb\label{s5.11}
\int_{(0,\infty)} u\nu(\d u)< \infty,
 \eeqlb
one can differentiate both sides of \eqref{s5.10} and use \eqref{s3.4} to see
 \beqlb\label{s5.12}
\int_{[0,\infty)} y P_t(x,\d y)
 =
x\e^{-bt} + \psi^\prime(0)\int_0^t\e^{-bs}\d s,
 \eeqlb
where
 \beqlb\label{s5.13}
\psi^\prime(0) = \beta + \int_{(0,\infty)} u \nu(\d u).
 \eeqlb

\bgproposition\label{ts5.7} Suppose that $\{(y_1(t),\mcr{G}_t^1): t\ge 0\}$ and $\{(y_2(t),\mcr{G}_t^2): t\ge 0\}$ are two independent CBI-processes with branching mechanism $\phi$ and immigration mechanisms $\psi_1$ and $\psi_2$, respectively. Let $y(t) = y_1(t)+y_2(t)$ and $\mcr{G}_t = \sigma(\mcr{G}_t^1\cup \mcr{G}_t^2)$. Then $\{(y(t),\mcr{G}_t): t\ge 0\}$ is a CBI-processes with branching mechanism $\phi$ and immigration mechanism $\psi=\psi_1+\psi_2$. \edproposition

\proof Let $t\ge r\ge 0$ and for $i=1,2$ let $F_i$ be a bounded positive $\mcr{G}_r^i$-measurable random variable. For any $\lambda\ge 0$ we have
 \beqnn
\bP\big[F_1F_2\e^{-\lambda y(t)}\big]
 \ar=\ar
\bP\big[F_1\e^{-\lambda y_1(t)}\big] \bP\big[F_2\e^{-\lambda y_2(t)}\big] \cr
 \ar=\ar
\bP\bigg[F_1\exp\bigg\{-y_1(r)v_{t-r}(\lambda)-\int_0^{t-r}\psi_1(v_s(\lambda))\d s\bigg\}\bigg] \cr
 \ar\ar\quad
\cdot\,\bP\bigg[F_2\exp\bigg\{-y_2(r)v_{t-r}(\lambda)-\int_0^{t-r}\psi_2(v_s(\lambda))\d s\bigg\}\bigg] \cr
 \ar=\ar
\bP\bigg[F_1F_2\exp\bigg\{-y(r)v_{t-r}(\lambda)-\int_0^{t-r}\psi(v_s(\lambda))\d s\bigg\}\bigg].
 \eeqnn
As in the proof of Proposition~\ref{ts2.12}, one can see $\{(y(t),\mcr{G}_t): t\ge 0\}$ is a CBI-processes with branching mechanism $\phi$ and immigration mechanism $\psi$. \qed

The next theorem follows by a modification of the proof of Theorem~\ref{ts2.13}.

\bgtheorem\label{ts5.8} Suppose that Conditions~\ref{ts2.3} and~\ref{ts5.3} are satisfied. Let $\{y(t): t\ge 0\}$ be a CBI-process with transition semigroup $(P_t)_{t\ge 0}$ defined by \eqref{s5.10}. For $k\ge 1$ let $\{z_k(n): n\ge 0\}$ be a Markov chain with state space $E_k := \{0,k^{-1},2k^{-1},\ldots\}$ and $n$-step transition probability $P_k^n(x,\d y)$ determined by \eqref{s5.4}. If $z_k(0)$ converges to $y(0)$ in distribution, then $\{z_k(\lfloor\gamma_kt\rfloor): t\ge 0\}$ converges to $\{y(t): t\ge 0\}$ in distribution on $D[0,\infty)$. \edtheorem

The convergence of rescaled GWI-processes to CBI-processes have been studied by many authors; see, e.g., Aliev (1985), Kawazu and Watanabe (1971) and Li (2006, 2011).

\bgexample\label{es5.1} The transition semigroup $(Q_t^b)_{t\ge 0}$ defined by \eqref{s3.14} corresponds to a CBI-process with branching mechanism $\phi$ and immigration mechanism $\phi_0^\prime$. \edexample

\bgexample\label{es5.2} Suppose that $c>0$, $0<\alpha\le 1$ and $b$ are constants and let $\phi(z) = cz^{1+\alpha} + bz$. Let $v_t(\lambda)$ and $q^0_\alpha(t)$ be defined as in Example~\ref{es2.3}. Let $\beta\ge 0$ and let $\psi(z) = \beta z^\alpha$. We can use \eqref{s5.10} to define the transition semigroup $(P_t)_{t\ge 0}$. It is easy to show that
 \beqnn
\int_{[0,\infty)} \e^{-\lambda y} P_t(x,\d y)
 =
\frac{1}{\big[1 + cq^b_\alpha(t)\lambda^\alpha\big]^{\beta/c\alpha}}\,
\e^{-xv_t(\lambda)}, \qquad \lambda\ge 0.
 \eeqnn
\edexample

\newpage

\section{Structures of sample paths}

 \setcounter{equation}{0}

In this section, we give some reconstructions of the type of Pitman and Yor (1982) for the CB- and CBI-processes, which reveal the structures of their sample paths. Let $(Q_t)_{t\ge 0}$ be the transition semigroup of the CB-process with branching mechanism $\phi$ given by \eqref{s2.13}. Let $(Q_t^\circ)_{t\ge 0}$ be the restriction of $(Q_t)_{t\ge 0}$ on $(0,\infty)$. Let $D[0,\infty)$ denote the space of positive c\`{a}dl\`{a}g paths on $[0,\infty)$. On this space, we define the $\sigma$-algebras $\mcr{A} = \sigma(\{w(s): 0\le s< \infty\})$ and $\mcr{A}_t = \sigma(\{w(s): 0\le s\le t\})$ for $t\ge 0$. For any $w\in D[0,\infty)$ let $\tau_0(w)= \inf\{s>0: w(s)=0\}$. Let $D_0[0,\infty)$ be the set of paths $w\in D[0,\infty)$ such that $w(t)=0$ for $t\ge \tau_0(w)$. Let $D_1[0,\infty)$ be the set of paths $w\in D_0[0,\infty)$ satisfying $w(0)= 0$. Then both $D_0[0,\infty)$ and $D_1[0,\infty)$ are $\mcr{A}$-measurable subsets of $D[0,\infty)$.

\bgtheorem\label{ts6.1} Suppose that $\phi^\prime(\infty) = \infty$ and let $(l_t)_{t>0}$ be the entrance law for $(Q_t^\circ)_{t\ge 0}$ determined by \eqref{s3.16}. Then there is a unique $\sigma$-finite measure $\bN_0$ on $(D[0,\infty),$ $\mcr{A})$ supported by $D_1[0,\infty)$ such that, for $0<t_1< t_2< \cdots< t_n$ and $x_1,x_2,\ldots,x_n$ $\in (0,\infty)$,
 \beqlb\label{s6.1}
\ar\ar\bN_0(w(t_1)\in \d x_1, w(t_2)\in \d x_2, \ldots, w(t_n)\in \d x_n) \ccr
 \ar\ar\qquad
=\, l_{t_1}(\d x_1) Q_{t_2-t_1}^\circ(x_1,\d x_2) \cdots Q_{t_n-t_{n-1}}^\circ(x_{n-1},\d x_n).
 \eeqlb
\edtheorem

\proof
Recall that $(Q_t^b)_{t\ge 0}$ is the transition semigroup on $[0,\infty)$ given by \eqref{s3.14}. Let $X = (D[0,\infty), \mcr{A}, \mcr{A}_t, w(t), \bQ_x)$ be the canonical realization of $(Q_t)_{t\ge 0}$ and $Y = (D[0,\infty),$ $\mcr{A}, \mcr{A}_t, w(t),$ $\bQ_x^b)$ the canonical realization of $(Q_t^b)_{t\ge 0}$. For any $T>0$ there is a probability measure $\bP_0^{b,T}$ on $(D[0,\infty),\mcr{A})$ so that
 \beqnn
\ar\ar\bP_0^{b,T}[F(\{w(s): s\ge 0\})G(\{w(T+s): s\ge 0\})] \ccr
 \ar\ar\qquad
=\, \bQ_0^b[F(\{w(s): s\ge 0\})\bQ_{w(T)}G(\{w(s): s\ge 0\})],
 \eeqnn
where $F$ is a bounded $\mcr{A}_T$-measurable function and $G$ is a bounded $\mcr{A}$-measurable function. The formula above means that under $\bP_0^{b,T}$ the random path $\{w(s): s\in [0,T]\}$ is a Markov process with initial state $w(0)=0$ and transition semigroup $(Q_t^b)_{t\ge 0}$ and $\{w(s): s\in [T,\infty)\}$ is a Markov process with transition semigroup $(Q_t)_{t\ge 0}$. Then Corollary~\ref{ts3.11} implies $\bP_0^{b,T}(w(s)=0)= Q_s^b(0,\{0\})= 0$ for every $0< s\le T$. Let $\bN_0^{b,T}(\d w) = \e^{-bT}w(T)^{-1}1_{\{w(T)>0\}} \bP_0^{b,T}(\d w)$. We have
 \beqlb\label{s6.2}
\ar\ar\bN_0^{b,T}(w(s_1)\in \d x_1, w(s_2)\in \d x_2, \ldots, w(s_m)\in
\d x_m, w(T)\in \d z, \ccr
 \ar\ar\qqquad\qqquad\qqquad
w(t_1)\in \d y_1, w(t_2)\in \d y_2, \ldots, w(t_n)\in \d y_n) \ccr
 \ar\ar\qquad
=\, Q^b_{s_1}(0,\d x_1) Q_{s_2-s_1}^b(x_1,\d x_2) \cdots Q_{s_m-s_{m-1}}^b(x_{m-1},\d x_m)Q_{T-s_m}^b(x_m,\d z) \ccr
 \ar\ar\qquad\qquad
\e^{-bT}z^{-1}Q_{t_1-T}(z,\d y_1) Q_{t_2-t_1}(y_1,\d y_2) \cdots Q_{t_n-t_{n-1}} (y_{n-1},\d y_n) \ccr
 \ar\ar\qquad
=\,Q^b_{s_1}(0,\d x_1) x_1^{-1}x_1 \cdots Q_{s_m-s_{m-1}}^b(x_{m-1},\d x_m)x_m^{-1}x_mQ_{T-s_m}^b(x_m,\d z) \ccr
 \ar\ar\qquad\qquad
\e^{-bT}z^{-1}Q_{t_1-T}(z,\d y_1) Q_{t_2-t_1}(y_1,\d y_2) \cdots Q_{t_n-t_{n-1}}(y_{n-1},\d y_n) \ccr
 \ar\ar\qquad
=\, l_{s_1}(\d x_1) Q_{s_2-s_1}^\circ(x_1,\d x_2) \cdots Q_{s_m-s_{m-1}}^\circ(x_{m-1},\d x_m)Q_{T-s_m}^\circ(x_m,\d z) \ccr
 \ar\ar\qquad\qquad
Q_{t_1-T}(z,\d y_1) Q_{t_2-t_1}(y_1,\d y_2) \cdots Q_{t_n-t_{n-1}}(y_{n-1},\d y_n),
 \eeqlb
where $0< s_1< \cdots< s_m< T< t_1< \cdots< t_n$ and $x_1,\ldots,x_m,z,y_1,\ldots,y_n\in (0,\infty)$. Then for any $T_1\ge T_2> 0$ the two measures $\bN_0^{b,T_1}$ and $\bN_0^{b,T_2}$ coincide on $\{w\in D_1[0,\infty): \tau_0(w)> T_1\}$, so the increasing limit $\bN_0:= \lim_{T\to 0} \bN_0^{b,T}$ exists and defines a $\sigma$-finite measure supported by $D_1[0,\infty)$. {From} \eqref{s6.2} we get \eqref{s6.1}. The uniqueness of the measure $\bN_0$ satisfying \eqref{s6.1} follows by the measure extension theorem. \qed

The elements of $D_1[0,\infty)$ are called \index{excursion} \textit{excursions} and the measure $\bN_0$ is referred to as the \index{excursion law} \textit{excursion law} for the CB-process. In view of \eqref{s3.17} and \eqref{s6.1}, we have formally, for $0<t_1< t_2< \cdots< t_n$ and $x_1,x_2,\ldots,x_n\in (0,\infty)$,
 \beqlb\label{s6.3}
\ar\ar\bN_0(w(t_1)\in \d x_1, w(t_2)\in \d x_2, \ldots, w(t_n)\in \d x_n) \ccr
 \ar\ar\qquad
= \lim_{x\to 0} x^{-1} Q_{t_1}^\circ(x,\d x_1) Q_{t_2-t_1}^\circ(x_1,\d x_2) \cdots Q_{t_n-t_{n-1}}^\circ(x_{n-1},\d x_n) \ccr
 \ar\ar\qquad
= \lim_{x\to 0} x^{-1} \bQ_x(w(t_1)\in \d x_1,w(t_2)\in \d x_2, \ldots, w(t_n)\in \d x_n),
 \eeqlb
which explains why $\bN_0$ is supported by $D_1[0,\infty)$.

{From} \eqref{s6.1} we see that the excursion law is Markovian, namely, the path $\{w(t): t>0\}$ behaves under this law as a Markov process with transition semigroup $(Q_t^\circ)_{t\ge 0}$. Based on the excursion law, a reconstruction of the CB-process is given in the following theorem:

\bgtheorem\label{ts6.2} Suppose that $\phi^\prime(\infty) = \infty$. Let $z\ge 0$ and let $N_z= \sum_{i=1}^\infty \delta_{w_i}$ be a Poisson random measure on $D[0,\infty)$ with intensity $z\bN_0(\d w)$. Let $X_0=z$ and for $t>0$ let
 \beqlb\label{s6.4}
X_t^z = \int_{D[0,\infty)} w(t) N_z(\d w)
 =
\sum_{i=1}^\infty w_i(t).
 \eeqlb
For $t\ge 0$ let $\mcr{G}_t^z$ be the $\sigma$-algebra generated by the collection of random variables $\{N_z(A): A\in \mcr{A}_t\}$. Then $\{(X_t^z,\mcr{G}_t^z): t\ge 0\}$ is a CB-process with branching mechanism $\phi$.
\edtheorem

\proof It is easy to see that $\{X_t^z: t\ge 0\}$ is adapted relative to the filtration $\{\mcr{G}_t^z: t\ge 0\}$. We claim that the random variable $X_t^z$ has distribution $Q_t(z,\cdot)$ on $[0,\infty)$. Indeed, for $t=0$ this is immediate. For any $t>0$ and $\lambda\ge 0$ we have
 \beqnn
\bP\big[\exp\{-\lambda X_t^z\}\big]
 \ar=\ar
\exp\bigg\{-z\int_{D[0,\infty)} (1-\e^{-\lambda w(t)}) \bN_0(\d w)\bigg\} \cr
 \ar=\ar
\exp\bigg\{-z\int_{(0,\infty)} (1 - \e^{-\lambda u})l_t(\d u)\bigg\}
 =
\exp\{-zv_t(\lambda)\}.
 \eeqnn
By the Markov property \eqref{s6.1}, for any $t\ge r>0$ and any bounded $\mcr{A}_r$-measurable function $H$ on $D[0,\infty)$ we have
 \beqnn
 \ar\ar\int_{D[0,\infty)} H(w)(1-\e^{-\lambda w(t)})\bN_0(\d w) \cr
 \ar\ar\qquad
=\, \int_{D[0,\infty)} H(w)(1-\e^{-v_{t-r}(\lambda)w(r)})
\bN_0(\d w).
 \eeqnn
It follows that, for any bounded positive $\mcr{A}_r$-measurable function $F$ on $D[0,\infty)$,
 \beqnn
 \ar\ar\bP\bigg[\exp\bigg\{-\int_{D[0,\infty)}F(w)N_z(\d w)\bigg\}
\cdot\exp\Big\{-\lambda X_t^z\Big\}\bigg] \cr
 \ar\ar\qquad
=\, \bP\bigg[\exp\bigg\{-\int_{D[0,\infty)}\big[F(w)+\lambda w(t)\big]N_z(\d w)\bigg\}\bigg] \cr
 \ar\ar\qquad
=\, \exp\bigg\{-z\int_{D[0,\infty)}\big(1 - \e^{-F(w)-\lambda w(t)}\big)
\bN_0(\d w)\bigg\} \cr
 \ar\ar\qquad
=\, \exp\bigg\{-z\int_{D[0,\infty)} \big(1 - \e^{-F(w)}\big)
\bN_0(\d w)\bigg\} \cr
 \ar\ar\qquad\qquad\qquad
\cdot\exp\bigg\{-z\int_{D[0,\infty)} \e^{-F(w)}\big(1 - \e^{- \lambda
w(t)}\big)\bN_0(\d w)\bigg\} \cr
 \ar\ar\qquad
=\, \exp\bigg\{-z\int_{D[0,\infty)} \big(1 - \e^{-F(w)}\big)
\bN_0(\d w)\bigg\} \cr
 \ar\ar\qquad\qquad\qquad
\cdot\exp\bigg\{-z\int_{D[0,\infty)} \e^{-F(w)}\big(1 - \e^{-w(r)v_{t-r}(\lambda)}\big)\bN_0(\d w)\bigg\} \cr
 \ar\ar\qquad
=\, \exp\bigg\{-z\int_{D[0,\infty)} \big(1 - \e^{-F(w)}\e^{-w(r)v_{t-r}(\lambda)}\big)
\bN_0(\d w)\bigg\} \cr
 \ar\ar\qquad
=\, \bP\bigg[\exp\bigg\{-\int_{D[0,\infty)}\big[F(w)+v_{t-r}(\lambda)w(r)\big]N_z(\d w)\bigg\}\bigg] \cr
 \ar\ar\qquad
=\, \bP\bigg[\exp\bigg\{-\int_{D[0,\infty)}F(w)N_z(\d w)\bigg\}
\cdot\exp\Big\{-v_{t-r}(\lambda)X_r^z\Big\}\bigg].
 \eeqnn
Clearly, the $\sigma$-algebra $\mcr{G}^z_r$ is generated by the collection of random variables
 \beqnn
\exp\bigg\{-\int_{D[0,\infty)} F(w) N_z(\d w)\bigg\},
 \eeqnn
where $F$ runs over all bounded positive $\mcr{A}_r$-measurable functions on $D[0,\infty)$. Then $\{(X_t^z,\mcr{G}_t^z): t\ge 0\}$ is a Markov process with transition semigroup $(Q_t)_{t\ge 0}$. \qed

The above theorem gives a description of the structures of the population represented by the CB-process. {From} \eqref{s6.4} we see that the population at any time $t>0$ consists of at most countably many families, which evolve as the excursions $\{w_i: i=1,2,\cdots\}$ selected by the Poisson random measure $N_z(\d w)$. Unfortunately, this reconstruction is only available under the condition $\phi^\prime(\infty)= \infty$. To give reconstructions of the CB- and CBI-processes when this condition is not necessarily satisfied, we need to consider some inhomogeneous immigration structures and more general Markovian measures on the path space. As a consequence of Theorem~\ref{ts6.1} we obtain the following:

\bgproposition\label{ts6.3} Let $\beta\ge 0$ be a constant and $\nu$ a $\sigma$-finite measure on $(0,\infty)$ such that $\int_{(0,\infty)} u\nu(\d u)< \infty$. If $\beta> 0$, assume in addition $\phi^\prime(\infty)= \infty$. Then we can define a $\sigma$-finite measure $\bN$ on $(D[0,\infty),\mcr{A})$ by
 \beqlb\label{s6.5}
\bN(\d w) = \beta\bN_0(\d w) + \int_{(0,\infty)}\nu(\d x)\bQ_x(\d w), \qquad w\in D[0,\infty).
 \eeqlb
Moreover, for $0<t_1< t_2< \cdots< t_n$ and $x_1,x_2,\ldots,x_n\in (0,\infty)$, we have
 \beqlb\label{s6.6}
\ar\ar\bN(w(t_1)\in \d x_1, w(t_2)\in \d x_2, \ldots, w(t_n)\in \d x_n) \ccr
 \ar\ar\qquad
=\, H_{t_1}(\d x_1) Q_{t_2-t_1}^\circ(x_1,\d x_2) \cdots Q_{t_n-t_{n-1}}^\circ(x_{n-1},\d x_n),
 \eeqlb
where
 \beqnn
H_t(\d y) = \beta l_t(\d y) + \int_{(0,\infty)}\nu(\d x)Q_t^\circ(x,\d y), \qquad y>0.
 \eeqnn
\edproposition

Clearly, the measure $\bN$ defined by \eqref{s6.5} is actually supported by $D_0[0,\infty)$. Under this law, the path $\{w(t): t>0\}$ behaves as a Markov process with transition semigroup $(Q_t^\circ)_{t\ge 0}$ and one-dimensional distributions $(H_t)_{t> 0}$.

\bgtheorem\label{ts6.4}
Suppose that the conditions of Proposition~\ref{ts6.3} are satisfied and let $\bN$ be defined by \eqref{s6.5}. Let $\rho$ be a Borel measure on $(0,\infty)$ such that $\rho(0,t]< \infty$ for each $0<t< \infty$. Suppose that $N= \sum_{i=1}^\infty \delta_{(s_i,w_i)}$ is a Poisson random measure on $(0,\infty)\times D[0,\infty)$ with intensity $\rho(\d s)\bN(\d w)$. For $t\ge 0$ let
 \beqlb\label{s6.7}
Y_t = \int_{(0,t]} \int_{D[0,\infty)} w(t-s) N(\d s,\d w)
 =
\sum_{0< s_i\le t}w_i(t-s_i)
 \eeqlb
and let $\mcr{G}_t$ be the $\sigma$-algebra generated by the random variables $\{N((0,u]\times A): A\in \mcr{A}_{t-u}, 0\le u\le t\}$. Then $\{(Y_t,\mcr{G}_t): t\ge 0\}$ is a Markov process with inhomogeneous transition semigroup $(P_{r,t})_{t\ge r\ge 0}$ given by
 \beqlb\label{s6.8}
\int_{[0,\infty)} \e^{-\lambda y} P_{r,t}(x,\d y)
 =
\exp\bigg\{-xv_{t-r}(\lambda)-\int_{(r,t]}\psi(v_{t-s}(\lambda))\rho(\d s)\bigg\}, ~~
 \eeqlb
where the function $\psi$ is given by \eqref{s5.9}.
\edtheorem

\proof {From} \eqref{s6.7} we see that $\{Y_t: t\ge 0\}$ is adapted to the filtration $\{\mcr{G}_t: t\ge 0\}$. Let $t\ge r\ge u\ge 0$ and let $F$ be a bounded positive function on $D[0,\infty)$ measurable relative to $\mcr{A}_{r-u}$. For $\lambda\ge 0$, we can use the Markov property \eqref{s6.6} to see
 \beqnn
\ar\ar\bP\bigg[\exp\bigg\{-\int_{(0,u]}\int_{D[0,\infty)} F(w)N(\d s,\d w) - \lambda Y_t\bigg\}\bigg] \cr
 \ar\ar\qquad
=\, \bP\bigg[\exp\bigg\{-\int_{(0,t]}\int_{D[0,\infty)} \big[F(w)1_{\{s\le u\}} + \lambda w(t-s)\big] N(\d s,\d w)\bigg\}\bigg] \cr
 \ar\ar\qquad
=\, \exp\bigg\{-\int_{(0,t]}\rho(\d s)\int_{D[0,\infty)}\big(1-\e^{-F(w) 1_{\{s\le u\}}}\e^{-\lambda w(t-s)}\big)\bN(\d w)\bigg\} \cr
 \ar\ar\qquad
=\, \exp\bigg\{-\int_{(0,r]}\rho(\d s)\int_{D[0,\infty)} \big(1-\e^{-F(w)1_{\{s\le u\}}} \e^{-w(t-s)}\big) \bN(\d w)\bigg\} \cr
 \ar\ar\qquad\qquad
\cdot\,\exp\bigg\{-\int_{(r,t]}\rho(\d s)\int_{D[0,\infty)}\big(1-\e^{-\lambda w(t-s)}\big) \bN(\d w)\bigg\} \cr
 \ar\ar\qquad
=\, \exp\bigg\{-\int_{(0,r]}\rho(\d s)\int_{D[0,\infty)} \big(1-\e^{-F(w)1_{\{s\le u\}}}\big) \bN(\d w)\bigg\} \cr
 \ar\ar\qquad\qquad
\cdot\,\exp\bigg\{-\int_{(0,r]}\rho(\d s)\int_{D[0,\infty)} \e^{-F(w)1_{\{s\le u\}}}\big(1- \e^{-w(t-s)}\big) \bN(\d w)\bigg\} \cr
 \ar\ar\qquad\qquad
\cdot\,\exp\bigg\{-\int_{(r,t]}\rho(\d s)\int_{D[0,\infty)}\big(1-\e^{-\lambda w(t-s)}\big) \bN(\d w)\bigg\} \cr
 \ar\ar\qquad
=\, \exp\bigg\{-\int_{(0,r]}\rho(\d s)\int_{D[0,\infty)} \big(1-\e^{-F(w)1_{\{s\le u\}}}\big) \bN(\d w)\bigg\} \cr
 \ar\ar\qquad\qquad
\cdot\,\exp\bigg\{-\int_{(0,r]}\rho(\d s)\int_{D[0,\infty)} \e^{-F(w)1_{\{s\le u\}}}\big(1- \e^{-v_{t-r}(\lambda)w(r-s)}\big) \bN(\d w)\bigg\} \cr
 \ar\ar\qquad\qquad
\cdot\,\exp\bigg\{-\int_{(r,t]}\rho(\d s)\int_{(0,\infty)}\big(1-\e^{-\lambda y}\big) H_{t-s}(\d w)\bigg\} \cr
 \ar\ar\qquad
=\, \exp\bigg\{-\int_{(0,r]}\rho(\d s)\int_{D[0,\infty)} \big(1-\e^{-F(w)1_{\{s\le u\}}} \e^{-v_{t-r}(\lambda)w(r-s)}\big) \bN(\d w)\bigg\} \cr
 \ar\ar\qquad\qquad
\cdot\exp\bigg\{-\beta\int_{(r,t]}\rho(\d s)\int_{(0,\infty)}(1-\e^{-\lambda y})l_{t-s}(\d y) \cr
 \ar\ar\qqquad\qqquad\qqquad
- \int_{(r,t]}\rho(\d s)\int_{(0,\infty)}(1-\e^{-\lambda y})\nu Q^\circ_{t-s}(\d y)\bigg\} \cr
 \ar\ar\qquad
=\, \bP\bigg[\exp\bigg\{-\int_{(0,r]}\int_{D[0,\infty)} \big[F(w)1_{\{s\le u\}} + v_{t-r}(\lambda)w(r-s)\big] N(\d s,\d w)\bigg\}\bigg] \cr
 \ar\ar\qquad\qquad
\cdot\exp\bigg\{-\int_{(r,t]}\bigg[\beta v_{t-s}(\lambda) + \int_{(0,\infty)}(1-\e^{-y v_{t-s}(\lambda)})\nu(\d y)\bigg]\rho(\d s)\bigg\} \cr
 \ar\ar\qquad
=\, \bP\bigg[\exp\bigg\{-\int_{(0,u]}\int_{D[0,\infty)} F(w) N(\d s,\d w)\bigg\} \cr
 \ar\ar\qqquad\qqquad\qqquad
\cdot\,\exp\bigg\{-v_{t-r}(\lambda)Y_r - \int_{(r,t]} \psi(v_{t-s}(\lambda))\rho(\d s)\bigg\}\bigg].
 \eeqnn
Then $\{(Y_t,\mcr{G}_t): t\ge 0\}$ is a Markov process in $[0,\infty)$ with inhomogeneous transition semigroup $(P_{r,t})_{t\ge r\ge 0}$ given by \eqref{s6.8}. \qed

\bgcorollary\label{ts6.5} Suppose that $\phi^\prime(\infty)= \infty$. Let $\beta> 0$ and let $N_\beta= \sum_{i=1}^\infty \delta_{(s_i,w_i)}$ be a Poisson random measure on $(0,\infty)\times D[0,\infty)$ with intensity $\beta\d s\bN_0(\d w)$. For $t\ge 0$ let
 \beqnn
Y_t^\beta = \int_{(0,t]} \int_{D[0,\infty)} w(t-s) N_\beta(\d s,\d w)
 =
\sum_{0<s_i\le t}w_i(t-s_i)
 \eeqnn
and let $\mcr{G}_t^\beta$ be the $\sigma$-algebra generated by the random variables $\{N_\beta((0,u]\times A): A\in \mcr{A}_{t-u}, 0\le u\le t\}$. Then $\{(Y_t^\beta,\mcr{G}_t^\beta): t\ge 0\}$ is a CBI-process with branching mechanism $\phi$ and immigration mechanism $\psi_\beta$ defined by $\psi_\beta(\lambda)= \beta\lambda$, $\lambda\ge 0$.
\edcorollary

\bgcorollary\label{ts6.6} Let $\nu$ be a $\sigma$-finite measure on $(0,\infty)$ such that $\int_{(0,\infty)} u\nu(\d u)< \infty$. Let $N_\nu= \sum_{i=1}^\infty \delta_{(s_i,w_i)}$ be a Poisson random measure on $(0,\infty)\times D[0,\infty)$ with intensity $\d s\bN_\nu(\d w)$, where
 \beqlb\label{s6.9}
\bN_\nu(\d w)= \int_{(0,\infty)}\nu(\d x)\bQ_x(\d w).
 \eeqlb
For $t\ge 0$ let
 \beqnn
Y_t^\nu = \int_{(0,t]} \int_{D[0,\infty)} w(t-s) N_\nu(\d s,\d w)
 =
\sum_{0<s_i\le t}w_i(t-s_i)
 \eeqnn
and let $\mcr{G}_t^\nu$ be the $\sigma$-algebra generated by the random variables $\{N_\nu((0,u]\times A): A\in \mcr{A}_{t-u}, 0\le u\le t\}$. Then $\{(Y_t^\nu,\mcr{G}_t^\nu): t\ge 0\}$ is a CBI-process with branching mechanism $\phi$ and immigration mechanism $\psi_\nu$ defined by
 \beqlb\label{s6.10}
\psi_\nu(\lambda) = \int_{(0,\infty)}(1-\e^{-u\lambda})\nu(\d u), \qquad \lambda\ge 0.
 \eeqlb
\edcorollary

The transition semigroup $(P_{r,t})_{t\ge r\ge 0}$ defined by \eqref{s6.8} is a generalization of the one given by \eqref{s5.10}; see also Li (1996, 2003) and Li (2011, p.224). A Markov process with transition semigroup $(P_{r,t})_{t\ge r\ge 0}$ is naturally called a CBI-process with \index{inhomogeneous immigration rate} \textit{inhomogeneous immigration rate} $\rho$. In view of \eqref{s6.7}, the population $\{Y_t: t\ge 0\}$ consists of a countable families of immigrants, whose immigration times $\{s_i: i=1,2,\cdots\}$ and evolution trajectories $\{w_i: i=1,2,\cdots\}$ are both selected by the Poisson random measure $N(\d s,\d w)$. The processes constructed in Corollaries~\ref{ts6.5} and~\ref{ts6.6} can be interpreted similarly.

\bgtheorem\label{ts6.7} Suppose that $\delta:= \phi^\prime(\infty)< \infty$. Let $z> 0$ and let $N_z = \sum_{i=1}^\infty \delta_{(s_i,w_i)}$ be a Poisson random measure on $(0,\infty)\times D[0,\infty)$ with intensity $z\e^{-\delta s}\d s\bN_m(\d w)$, where $\bN_m$ is defined by \eqref{s6.9} with $\nu=m$. For $t\ge 0$ let
 \beqlb\label{s6.11}
X_t^z= z\e^{-\delta t} + \int_{(0,t]} \int_{D[0,\infty)} w(t-s) N_z(\d s,\d w)
 \eeqlb
and let $\mcr{G}_t^z$ be the $\sigma$-algebra generated by the random variables $\{N_z((0,u]\times A): A\in \mcr{A}_{t-u}, 0\le u\le t\}$. Then $\{(X_t^z,\mcr{G}_t^z): t\ge 0\}$ is a CB-process with branching mechanism $\phi$.
\edtheorem

\proof Let $Z_t= X_t^z- z\e^{-\delta t}$ denote the second term on the right-hand side of \eqref{s6.11}. By Theorem~\ref{ts6.4} we infer that $\{(Z_t,\mcr{G}_t^z): t\ge 0\}$ is a Markov process with inhomogeneous transition semigroup $(P_{r,t}^z)_{t\ge r\ge 0}$ given by
 \beqnn
\int_{[0,\infty)} \e^{-\lambda y} P_{r,t}^z(x,\d y)
 =
\exp\bigg\{-xv_{t-r}(\lambda)-z\int_r^t\psi_m(v_{t-s}(\lambda))\e^{-\delta s}\d s\bigg\},
 \eeqnn
where $\psi_m$ is defined by \eqref{s6.10} with $\nu=m$. Let $t\ge r\ge u\ge 0$ and let $F$ be a bounded positive $\mcr{A}_{r-u}$-measurable function on $D[0,\infty)$. For $\lambda\ge 0$ we have
 \beqnn
\ar\ar\bP\bigg[\exp\bigg\{-\int_{(0,u]}\int_{D[0,\infty)} F(w)N_z(\d s,\d w) - \lambda X_t^z\bigg\}\bigg] \cr
 \ar\ar\qquad
=\, \bP\bigg[\exp\bigg\{-\int_{(0,u]}\int_{D[0,\infty)} F(w)N_z(\d s,\d w) - \lambda z\e^{-\delta t} - \lambda Z_t\bigg\}\bigg] \cr
 \ar\ar\qquad
=\, \bP\bigg[\exp\bigg\{-\int_{(0,u]}\int_{D[0,\infty)} F(w)N_z(\d s,\d w) - \lambda z\e^{-\delta t}\bigg\}\bigg] \cr
 \ar\ar\qquad\qquad
\cdot\exp\bigg\{- v_{t-r}(\lambda)Z_r - z\int_r^t\psi_m(v_{t-s}(\lambda))\e^{-\delta s}\d s\bigg\} \cr
 \ar\ar\qquad
=\, \bP\bigg[\exp\bigg\{-\int_{(0,u]}\int_{D[0,\infty)} F(w)N_z(\d s,\d w) - \lambda z\e^{-\delta t}\bigg\}\bigg] \cr
 \ar\ar\qquad\qquad
\cdot\exp\bigg\{-v_{t-r}(\lambda)Z_r - \e^{-\delta r}z\int_0^{t-r} \psi_m(v_{t-r-s}(\lambda)) \e^{-\delta s}\d s\bigg\} \cr
 \ar\ar\qquad
=\, \bP\bigg[\exp\bigg\{-\int_{(0,u]}\int_{D[0,\infty)} F(w)N_z(\d s,\d w) - v_{t-r}(\lambda)Z_r - z\e^{-\delta r}v_{t-r}(\lambda)\bigg\}\bigg] \cr
 \ar\ar\qquad
=\, \bP\bigg[\exp\bigg\{-\int_{(0,u]}\int_{D[0,\infty)} F(w)N_z(\d s,\d w) - v_{t-r}(\lambda)X_r^z\bigg\}\bigg],
 \eeqnn
where we have used \eqref{s3.20}. Then $\{(X_t^z,\mcr{G}_t^z): t\ge 0\}$ is a CB-process with transition semigroup $(Q_t)_{t\ge 0}$ defined by \eqref{s2.19}. \qed

\bgtheorem\label{ts6.8} Suppose that $\delta:= \phi^\prime(\infty)< \infty$. Let $\beta>0$ and let $N_\beta = \sum_{i=1}^\infty \delta_{(s_i,w_i)}$ be a Poisson random measure on $(0,\infty)\times D[0,\infty)$ with intensity $\beta\delta^{-1}(1-\e^{-\delta s})$ $\d s\bN_m(\d w)$, where $\bN_m$ is defined by \eqref{s6.9} with $\nu=m$. For $t\ge 0$ let
 \beqlb\label{s6.12}
Y_t^\beta= \beta\delta^{-1}(1-\e^{-\delta t}) + \int_{(0,t]} \int_{D[0,\infty)} w(t-s) N_\beta(\d s,\d w)
 \eeqlb
and let $\mcr{G}_t^\beta$ be the $\sigma$-algebra generated by the random variables $\{N_\beta((0,u]\times A): A\in \mcr{A}_{t-u}, 0\le u\le t\}$. Then $\{(X_t^\beta,\mcr{G}_t^\beta): t\ge 0\}$ is a CBI-process with branching mechanism $\phi$ and immigration mechanism $\psi_\beta$ defined by $\psi_\beta(\lambda)= \beta\lambda$, $\lambda\ge 0$.
\edtheorem

\proof Let $Z_t$ denote the second term on the right-hand side of \eqref{s6.12}. By Theorem~\ref{ts6.4}, the process $\{(Z_t,\mcr{G}_t^\beta): t\ge 0\}$ is a Markov process with inhomogeneous transition semigroup $(P_{r,t}^\beta)_{t\ge r\ge 0}$ given by
 \beqnn
\int_{[0,\infty)} \e^{-\lambda y} P_{r,t}^\beta(x,\d y)
 =
\exp\bigg\{-xv_{t-r}(\lambda) - \beta\delta^{-1} \int_r^t \psi_m(v_{t-s}(\lambda)) (1-\e^{-\delta s})\d s\bigg\},
 \eeqnn
where $\psi_m$ is defined by \eqref{s6.10} with $\nu=m$. For $t\ge 0$ and $\lambda\ge 0$, we can use Theorem~\ref{ts3.15} to see
 \beqlb\label{s6.13}
\int_0^tv_s(\lambda)\d s
 \ar=\ar
\lambda\int_0^t\e^{-\delta s}\d s + \int_0^t\d s\int_0^s\e^{-\delta u} \psi_m(v_{s-u}(\lambda))\d u \cr
 \ar=\ar
\lambda\int_0^t\e^{-\delta s}\d s + \int_0^t\d s\int_0^{t-s}\e^{-\delta u} \psi_m(v_{t-s-u}(\lambda))\d u \cr
 \ar=\ar
\lambda\delta^{-1}(1-\e^{-\delta t}) + \int_0^t\d s\int_s^t \e^{-\delta(u-s)} \psi_m(v_{t-u}(\lambda)) \d u \cr
 \ar=\ar
\lambda\delta^{-1}(1-\e^{-\delta t}) + \int_0^t\d u\int_0^u \e^{-\delta(u-s)} \psi_m(v_{t-u}(\lambda)) \d s \cr
 \ar=\ar
\lambda\delta^{-1}(1-\e^{-\delta t}) + \delta^{-1}\int_0^t(1-\e^{-\delta u}) \psi_m(v_{t-u}(\lambda))\d u.\quad
 \eeqlb
Let $t\ge r\ge u\ge 0$ and let $F$ be a bounded positive $\mcr{A}_{r-u}$-measurable function on $D[0,\infty)$. Then
 \beqnn
\ar\ar\bP\bigg[\exp\bigg\{-\int_{(0,u]}\int_{D[0,\infty)} F(w)N_\beta(\d s,\d w) - \lambda Y_t^\beta\bigg\}\bigg] \cr
 \ar\ar\qquad
=\, \bP\bigg[\exp\bigg\{-\int_{(0,u]}\int_{D[0,\infty)} F(w)N_\beta(\d s,\d w) - \lambda\beta\delta^{-1}(1-\e^{-\delta t}) - \lambda Z_t\bigg\}\bigg] \cr
 \ar\ar\qquad
=\, \bP\bigg[\exp\bigg\{-\int_{(0,u]}\int_{D[0,\infty)} F(w)N_\beta(\d s,\d w) - \lambda\beta\delta^{-1}(1-\e^{-\delta t})\bigg\}\bigg] \cr
 \ar\ar\qquad\qquad
\cdot\exp\bigg\{- v_{t-r}(\lambda)Z_r - \beta\delta^{-1} \int_r^t \psi_m(v_{t-s}(\lambda)) (1-\e^{-\delta s})\d s\bigg\} \cr
 \ar\ar\qquad
=\, \bP\bigg[\exp\bigg\{-\int_{(0,u]}\int_{D[0,\infty)} F(w)N_\beta(\d s,\d w) - \lambda\beta\delta^{-1}(1-\e^{-\delta t})\bigg\}\bigg] \cr
 \ar\ar\qquad\qquad
\cdot\exp\bigg\{-v_{t-r}(\lambda)Y_r^\beta + v_{t-r}(\lambda) \beta\delta^{-1}(1-\e^{-\delta r})\bigg\} \cr
 \ar\ar\qquad\qquad
\cdot\exp\bigg\{- \beta\delta^{-1}\int_0^{t-r} \psi_m(v_{t-r-s}(\lambda)) (1-\e^{-\delta(r+s)})\d s\bigg\} \cr
 \ar\ar\qquad
=\, \bP\bigg[\exp\bigg\{-\int_{(0,u]}\int_{D[0,\infty)} F(w)N_\beta(\d s,\d w) - \lambda\beta\delta^{-1}(1-\e^{-\delta t})\bigg\}\bigg] \cr
 \ar\ar\qquad\qquad
\cdot\exp\bigg\{-v_{t-r}(\lambda)Y_r^\beta + \lambda\beta\delta^{-1}\e^{-\delta(t-r)} (1-\e^{-\delta r})\bigg\} \cr
 \ar\ar\qquad\qquad
\cdot\exp\bigg\{\beta\delta^{-1}(1-\e^{-\delta r})\int_0^{t-r}\e^{-\delta s} \psi_m(v_{t-r-s}(\lambda))\d s\bigg\} \cr
 \ar\ar\qquad\qquad
\cdot\exp\bigg\{- \beta\delta^{-1}\int_0^{t-r} \psi_m(v_{t-r-s}(\lambda)) (1-\e^{-\delta(r+s)})\d s\bigg\} \cr
 \ar\ar\qquad
=\, \bP\bigg[\exp\bigg\{-\int_{(0,u]}\int_{D[0,\infty)} F(w)N_\beta(\d s,\d w) - \lambda\beta\delta^{-1}(1-\e^{-\delta(t-r)})\bigg\}\bigg] \cr
 \ar\ar\qquad\qquad
\cdot\exp\bigg\{-v_{t-r}(\lambda)Y_r^\beta - \beta\delta^{-1}\int_0^{t-r} \psi_m(v_{t-r-s}(\lambda)) (1-\e^{-\delta s})\d s\bigg\} \cr
 \ar\ar\qquad
=\, \bP\bigg[\exp\bigg\{-\int_{(0,u]}\int_{D[0,\infty)} F(w)N_\beta(\d s,\d w) \cr
 \ar\ar\qqquad\qqquad\qqquad\qqquad
-\, v_{t-r}(\lambda)Y_r^\beta - \beta\int_0^{t-r} v_s(\lambda))\d s\bigg\}\bigg],
 \eeqnn
where we used \eqref{s6.13} for the last equality. That gives the desired result. \qed

Since the CB- and CBI-processes have Feller transition semigroups, they have c\`{a}dl\`{a}g realizations. By Proposition~A.7 of Li (2011), any realizations of the processes has a c\`{a}dl\`{a}g modification. In the case of $\phi^\prime(\infty)= \infty$, let $(X_t^z,\mcr{G}_t^z)$, $(Y_t^\beta,\mcr{G}_t^\beta)$ and $(Y_t^\nu,\mcr{G}_t^\nu)$ be defined as in Theorem~\ref{ts6.2}, Corollary~\ref{ts6.5} and Corollary~\ref{ts6.6}, respectively. In the case of $\phi^\prime(\infty)< \infty$, we define those processes as in Theorem~\ref{ts6.7}, Theorem~\ref{ts6.8} and Corollary~\ref{ts6.6}, respectively. In both cases, let $Y_t = X_t^z + Y_t^\beta + Y_t^\nu$ and $\mcr{G}_t = \sigma(\mcr{G}_t^z\cup \mcr{G}_t^\beta\cup \mcr{G}_t^\nu)$. It is not hard to show that $\{(Y_t,\mcr{G}_t): t\ge 0\}$ is a CBI-process with branching mechanism $\phi$ given by \eqref{s2.13} and immigration mechanism $\psi$ given by \eqref{s5.9}.

The existence of the excursion law for the branching mechanism $\phi(z) = bz + cz^2$ was first proved by Pitman and Yor (1982). As a special case of the so-called \index{Kuznetsov measure} \textit{Kuznetsov measure}, the existence of the law for measure-valued branching processes was derived from a general result on Markov processes in Li (2003, 2011), where it was also shown that the law only charges sample paths starting with zero. In the setting of measure-valued processes, El~Karoui and Roelly (1991) used \eqref{s6.3} to construct the excursion law; see also Duquesne and Labb\'{e} (2014). The construction of the CB- or CBI-process based on a excursion law was first given by Pitman and Yor (1982). This type of constructions have also been used in the measure-valued setting by a number of authors; see, e.g., Dawson and Li (2003), El~Karoui and Roelly (1991), Li (1996, 2003, 2011) and Li and Shiga (1995).

\newpage

\section{Martingale problem formulations}

 \setcounter{equation}{0}

In this section we give several formulations of the CBI-process in terms of martingale problems. Let $(\phi,\psi)$ be given by \eqref{s2.13} and \eqref{s5.9}, respectively. We assume \eqref{s5.11} is satisfied and define $\psi^\prime(0)$ by \eqref{s5.13}. Let $C^2[0,\infty)$ denote the set of bounded continuous real functions on $[0,\infty)$ with bounded continuous derivatives up to the second order. For $f\in C^2[0,\infty)$ define
 \beqlb\label{s7.1}
Lf(x) \ar=\ar cxf^{\prime\prime}(x) + x\int_{(0,\infty)} \big[f(x+z) - f(x) -
zf^\prime(x)\big] m(\d z) \cr
 \ar\ar
+\, (\beta - bx)f^\prime(x) + \int_{(0,\infty)} \big[f(x+z) - f(x)\big] \nu(\d
z).
 \eeqlb
We shall identify the operator $L$ as the generator of the CBI-process.

\bgproposition\label{ts7.1} Let $(P_t)_{t\ge 0}$ be the transition semigroup defined by \eqref{s2.19} and \eqref{s5.10}. Then for any $t\ge 0$ and
$\lambda\ge 0$ we have
 \beqlb\label{s7.2}
\int_{[0,\infty)} \e^{-\lambda y} P_t(x,\d y)
 =
\e^{-x\lambda} + \int_0^t \d s\int_{[0,\infty)} [y\phi(\lambda) -
\psi(\lambda)] \e^{-y\lambda} P_s(x,\d y).
 \eeqlb
\edproposition

\proof Recall that $v_t^\prime(\lambda)= (\partial/\partial \lambda) v_t(\lambda)$. By differentiating both sides of \eqref{s5.10} we get
 \beqnn
\int_{[0,\infty)} y\e^{-y\lambda} P_t(x,\d y)
 \ar=\ar
\int_{[0,\infty)} \e^{-y\lambda} P_t(x,\d y)\bigg[xv_t^\prime(\lambda)
+ \int_0^t \psi^\prime(v_s(\lambda))v_s^\prime(\lambda)\d s\bigg].
 \eeqnn
{From} this and \eqref{s3.7} it follows that
 \beqnn
{\partial\over\partial t}\int_{[0,\infty)} \e^{-y\lambda} P_t(x,\d y)
 \ar=\ar
- \bigg[x{\partial\over \partial t}v_t(\lambda) + \psi(v_t(\lambda))\bigg]
\int_{[0,\infty)} \e^{-y\lambda} P_t(x,\d y) \cr
 \ar=\ar
\bigg[x\phi(\lambda) v_t^\prime(\lambda) - \psi(\lambda)\bigg]\int_{[0,\infty)}
\e^{-y\lambda} P_t(x,\d y) \cr
 \ar\ar
- \int_0^t \psi^\prime(v_s(\lambda)){\partial\over \partial
s}v_s(\lambda) \d s \int_{[0,\infty)} \e^{-y\lambda} P_t(x,\d y) \cr
 \ar=\ar
\bigg[x\phi(\lambda) v_t^\prime(\lambda) - \psi(\lambda)\bigg]\int_{[0,\infty)}
\e^{-y\lambda} P_t(x,\d y) \cr
 \ar\ar
+\, \phi(\lambda)\int_0^t \psi^\prime(v_s(\lambda))v_s^\prime(\lambda) \d
s \int_{[0,\infty)} \e^{-y\lambda} P_t(x,\d y) \cr
 \ar=\ar
\int_{[0,\infty)} [y\phi(\lambda) - \psi(\lambda)] \e^{-y\lambda}
P_t(x,\d y).
 \eeqnn
That gives \eqref{s7.2}. \qed

Suppose that $(\Omega, \mcr{G}, \mcr{G}_t, \bP)$ is a filtered probability space satisfying the usual hypotheses and $\{y(t): t\ge 0\}$ is a c\`{a}dl\`{a}g process in $[0,\infty)$ that is adapted to $(\mcr{G}_t)_{t\ge 0}$ and satisfies $\bP[y(0)]< \infty$. Let $C^{1,2}([0,\infty)^2)$ be the set of bounded continuous real functions $(t,x)\mapsto G(t,x)$ on $[0,\infty)^2$ with bounded continuous derivatives up to the first order relative to $t\ge 0$ and up to the second order relative to $x\ge 0$. Let us consider the following properties:
 \benumerate

\itm[{\rm(1)}] For every $T\ge 0$ and $\lambda\ge 0$,
 \beqnn
\exp\bigg\{-v_{T-t}(\lambda)y(t) - \int_0^{T-t} \psi(v_s(\lambda))\d
s\bigg\}, \qquad 0\le t\le T,
 \eeqnn
is a martingale.

\itm[{\rm(2)}] For every $\lambda\ge 0$,
 \beqnn
H_t(\lambda)
 :=
\exp\bigg\{-\lambda y(t)+\int_0^t[\psi(\lambda) - y(s)\phi(\lambda)]\d
s\bigg\}, \quad t\ge 0,
 \eeqnn
is a local martingale.

\itm[{\rm(3)}] The process $\{y(t): t\ge 0\}$ has no negative jumps and the optional random measure
 \beqnn
N_0(\d s,\d z) := \sum_{s>0}1_{\{\Delta y(s)\ne 0\}}\delta_{(s,\Delta y(s))}(\d s,\d z),
 \eeqnn
where $\Delta y(s) = y(s)-y(s-)$, has predictable compensator $\hat{N}_0(\d s,\d z) = \d s\nu(\d z) + y(s-)\d sm(\d z)$. Let $\tilde{N}_0(\d s,\d z) = N_0(\d s,\d z) - \hat{N}_0(\d s,\d z)$. We have
 \beqnn
y(t) = y(0) + M^c(t) + M^d(t) - b\int_0^t y(s-)\d s + \psi^\prime(0)t,
 \eeqnn
where $\{M^c(t): t\ge 0\}$ is a continuous local martingale with quadratic variation $2cy(t-)\d t$ and
 \beqnn
M^d(t) = \int_0^t\int_{(0,\infty)} z \tilde{N}_0(\d s,\d z), \qquad t\ge 0,
 \eeqnn
is a purely discontinuous local martingale.

\itm[{\rm(4)}] For every $f\in C^2[0,\infty)$ we have
 \beqlb\label{s7.3}
f(y(t)) = f(y(0)) + \int_0^t Lf(y(s))\d s + \mbox{local mart.}
 \eeqlb

\itm[{\rm(5)}] For any $G\in C^{1,2}([0,\infty)^2)$ we
have
 \beqlb\label{s7.4}
G(t,y(t))= G(0,y(0)) + \int_0^t \big[G^\prime_t(s,y(s)) + LG(s,y(s))\big]\d s + \mbox{local mart.}\quad
 \eeqlb
where $L$ acts on the function $x\mapsto G(s,x)$.

 \eenumerate

\bgtheorem\label{ts7.2} The above properties {\rm(1)}, {\rm(2)}, {\rm(3)}, {\rm(4)} and {\rm(5)} are equivalent to each other. Those properties hold if and only if
$\{(y(t),\mcr{G}_t): t\ge 0\}$ is a CBI-process with branching mechanism $\phi$ and immigration mechanism $\psi$. \edtheorem

\proof Clearly, {\rm(1)} holds if and only if $\{y(t): t\ge 0\}$ is a Markov
process relative to $(\mcr{G}_t)_{t\ge 0}$ with transition semigroup
$(P_t)_{t\ge 0}$ defined by \eqref{s5.10}. Then we only need to
prove the equivalence of the five properties.

{\rm(1)}$\Rightarrow${\rm(2)}: Suppose that {\rm(1)} holds. Then $\{y(t):
t\ge 0\}$ is a CBI-process with transition semigroup $(P_t)_{t\ge
0}$ given by \eqref{s5.10}. By \eqref{s7.2} and the Markov
property it is easy to see that
 \beqlb\label{s7.5}
Y_t(\lambda) := \e^{-\lambda y(t)} + \int_0^t [\psi(\lambda) -
y(s)\phi(\lambda)] \e^{-\lambda y(s)}\d s
 \eeqlb
is a martingale. By integration by parts applied to
 \beqlb\label{s7.6}
Z_t(\lambda) := \e^{-\lambda y(t)}
 ~~\mbox{and}~~
W_t(\lambda) := \exp\bigg\{\int_0^t[\psi(\lambda) - y(s)\phi(\lambda)]\d
s\bigg\}
 \eeqlb
we obtain
 \beqnn
\d H_t(\lambda)
 =
\e^{-\lambda y(t-)}\d W_t(\lambda) + W_t(\lambda)\d\e^{-\lambda y(t)}
 =
W_t(\lambda)\d Y_t(\lambda).
 \eeqnn
Then $\{H_t(\lambda)\}$ is a local martingale.

{\rm(2)}$\Rightarrow${\rm(3)}: For any $\lambda\ge 0$ let $Z_t(\lambda)$ and $W_t(\lambda)$ be defined by \eqref{s7.6}. We have $Z_t(\lambda)= H_t(\lambda)W_t(\lambda)^{-1}$ and so
 \beqlb\label{s7.7}
\d Z_t(\lambda)= W_t(\lambda)^{-1}\d H_t(\lambda) - Z_{t-}(\lambda)[\psi(\lambda) -
\phi(\lambda)y(t-)] \d t
 \eeqlb
by integration by parts. Then the strictly positive process $t\mapsto Z_t(\lambda)$ is a special semi-martingale; see, e.g., Dellacherie and Meyer (1982, p.213). By It\^{o}'s formula we find $t\mapsto y(t)$ is also a special semi-martingale. Now define the optional random measure $N_0(\d s,\d z)$ on $[0,\infty)\times \mbb{R}$ by
 \beqnn
N_0(\d s,\d z) = \sum_{s>0}1_{\{\Delta y(s)\ne 0\}} \delta_{(s,\Delta y(s))}(\d s,\d z),
 \eeqnn
where $\Delta y(s) = y(s) - y(s-)$. Let $\hat{N}_0(\d s,\d z)$ denote the predictable compensator of $N_0(\d s,\d z)$ and let $\tilde{N}_0(\d s,\d z)$ denote the compensated random measure; see, e.g., Dellacherie and Meyer (1982, p.375). It follows that
 \beqlb\label{s7.8}
y(t) = y(0) + U(t) + M^c(t) + M^d(t),
 \eeqlb
where $t\mapsto U(t)$ is a predictable process with locally bounded variations, $t\mapsto M^c(t)$ is a continuous local martingale and
 \beqnn
t\mapsto M^d(t) := \int_0^t\int_{\mbb{R}} z \tilde{N}_0(\d s,\d z)
 \eeqnn
is a purely discontinuous local martingale; see, e.g., Dellacherie and Meyer (1982, p.353 and p.376) or Jacod and Shiryaev (2003, pp.84--85). Let $\{C(t)\}$ denote the quadratic variation process of $\{M^c(t)\}$. By It\^{o}'s formula,
 \beqlb\label{s7.9}
Z_t(\lambda) \ar=\ar Z_0(\lambda) - \lambda\int_0^tZ_{s-}(\lambda)\d y(s)
+ \frac{1}{2}\lambda^2\int_0^t Z_{s-}(\lambda)\d C(s) \cr
 \ar\ar
+ \int_0^t\int_{\mbb{R}} Z_{s-}(\lambda) \big(\e^{-z\lambda} - 1 +
z\lambda\big) N_0(\d s,\d z) \cr
 \ar=\ar
Z_0(\lambda) - \lambda\int_0^tZ_{s-}(\lambda)\d U(s)
+ \frac{1}{2}\lambda^2\int_0^t Z_{s-}(\lambda)\d C(s) \cr
 \ar\ar
+ \int_0^t\int_{\mbb{R}} Z_{s-}(\lambda) \big(\e^{-z\lambda} - 1 +
z\lambda\big) \hat{N}_0(\d s,\d z) + \mbox{local mart.}
 \eeqlb
In view of \eqref{s7.7} and \eqref{s7.9} we get
 \beqnn
[y(s-)\phi(\lambda) - \psi(\lambda)]\d s
 =
- \lambda\d U(s) + \frac{1}{2}\lambda^2\d C(s) + \int_{\mbb{R}} \big(\e^{-z\lambda} - 1 + z\lambda\big) \hat{N}_0(\d s,\d z)
 \eeqnn
by the uniqueness of canonical decompositions of special semi-martingales; see, e.g., Dellacherie and Meyer (1982, p.213). By substituting the representations \eqref{s2.13} and \eqref{s5.9} for $\phi$ and $\psi$ into the above equation and comparing both sides we find
 \beqnn
~\d C(s) = 2cy(s-)\d s, ~ \d U(s) = [\psi'(0)-by(s-)]\d s
 \eeqnn
and
 \beqnn
\hat{N}_0(\d s,\d z) = \d s\nu(\d z) + y(s-)\d sm(\d z).
 \eeqnn
Then the process $t\mapsto y(t)$ has no negative jumps.

{\rm(3)}$\Rightarrow${\rm(4)}: This follows by It\^{o}'s formula.

{\rm(4)}$\Rightarrow${\rm(5)}: For $t\ge 0$ and $k\ge1$ we have
 \beqnn
G(t,y(t))
 \ar=\ar
G(0,y(0)) + \sum_{j=0}^\infty \big[G(t\land j/k,y(t\land (j+1)/k)) \cr
 \ar\ar\qqquad\qqquad\qqquad
-\, G(t\land j/k,y(t\land j/k))\big] \cr
 \ar\ar
+ \sum_{j=0}^\infty \big[G(t\land (j+1)/k,y(t\land (j+1)/k)) \cr
 \ar\ar\qqquad\qqquad
-\, G(t\land j/k,y(t\land (j+1)/k))\big],
 \eeqnn
where the summations only consist of finitely many nontrivial terms. By
applying {\rm(4)} term by term we obtain
 \beqnn
G(t,y(t))
 \ar=\ar
G(0,y(0)) + \sum_{j=0}^\infty \int_{t\land j/k}^{t\land (j+1)/k}
\bigg\{[\beta - by(s)]G^\prime_x(t\land j/k,y(s)) \cr
 \ar\ar
+\, cy(s)G^{\prime\prime}_{xx}(t\land j/k,y(s)) + y(s)\int_{(0,\infty)}
\Big[G(t\land j/k,y(s) + z) \cr
 \ar\ar
-\, G(t\land j/k,y(s)) - z G^\prime_x(t\land j/k,y(s))\Big] m(\d z) \cr
 \ar\ar
+ \int_{(0,\infty)} \Big[G(t\land j/k,y(s) + z) - G(t\land j/k,y(s))\Big]
\nu(\d z)\bigg\}\d s \cr
 \ar\ar
+ \sum_{j=0}^\infty\int_{t\land j/k}^{t\land (j+1)/k}
G^\prime_t(s,y(t\land (j+1)/k))\d s + M_k(t),
 \eeqnn
where $\{M_k(t)\}$ is a local martingale. Since $\{y(t)\}$ is a
c\`{a}dl\`{a}g process, letting $k\to\infty$ in the equation above gives
 \beqnn
G(t,y(t)) \ar=\ar G(0,y(0)) + \int_0^t \bigg\{G^\prime_t(s,y(s)) + [\beta
- by(s)] G^\prime_x(s,y(s)) \cr
 \ar\ar
+\, cy(s)G^{\prime\prime}_{xx}(s,y(s)) + y(s)\int_{(0,\infty)}
\Big[G(s,y(s)+z) \cr
 \ar\ar
-\, G(s,y(s)) - z G^\prime_x(s,y(s))\Big]m(\d z) \cr
 \ar\ar
+ \int_{(0,\infty)}\Big[G(s,y(s)+z) - G(s,y(s))\Big]\nu(\d z)\bigg\}\d s + M(t),
 \eeqnn
where $\{M(t)\}$ is a local martingale. Then we have \eqref{s7.4}.

{\rm(5)}$\Rightarrow${\rm(1)}: For fixed $T\ge 0$ and $\lambda\ge 0$ we define the function
 \beqnn
G_T(t,x) = \exp\bigg\{-v_{T-t}(\lambda)x-\int_0^{T-t}\psi(v_s(\lambda))\d s\bigg\}, \quad 0\le t\le T, x\ge0,
 \eeqnn
which can be extended to a function in $C^{1,2}([0,\infty)^2)$. Using \eqref{s3.6} we see
 \beqnn
{\d\over\d t}G_T(t,x) + LG_T(t,x)= 0, \quad 0\le t\le T, x\ge0,
 \eeqnn
Then \eqref{s7.4} implies that $t\mapsto G(t\land T,y(t\land T))$ is a local martingale, and hence a martingale by the boundedness. \qed

\bgcorollary\label{ts7.3} Let $\{(y(t), \mcr{G}_t): t\ge 0\}$ be a c\`{a}dl\`{a}g realization of the CBI-process satisfying $\bP[y(0)]< \infty$. Then for every $T\ge 0$ there is a constant $C_T\ge 0$ such
that
 \beqnn
\bP\Big[\sup_{0\le t\le T}y(t)\Big]
 \le
C_T\big\{\bP[y(0)] + \psi^\prime(0) + \sqrt{\bP[y(0)]} + \sqrt{\psi^\prime(0)}\big\}.
 \eeqnn
\edcorollary

\proof By the above property (3) and Doob's martingale inequality we have
 \beqnn
\ar\ar\bP\Big[\sup_{0\le t\le T}|y(t)-y(0)|\Big] \cr
 \ar\ar\qquad
\le T\psi^\prime(0) + \bP\bigg[|b|\int_0^T y(s)\d s\bigg] + \bP\Big[\sup_{0\le t\le T}|M_t^c|\Big] \cr
 \ar\ar\qquad\quad
+\, \bP\bigg[\int_0^T\int_{(1,\infty)} zN_0(\d s,\d z)\bigg] + \bP\bigg[\int_0^T\int_{(1,\infty)} z\hat{N}_0(\d s,\d z)\bigg] \cr
 \ar\ar\qquad\quad
+\, \bP\bigg[\sup_{0\le t\le T} \bigg|\int_{(0,1]} \int_0^1 z\tilde{N}_0(\d s,\d z)\bigg|\bigg] \cr
 \ar\ar\qquad
\le T\psi^\prime(0) + \bP\bigg[|b|\int_0^T y(s)\d s\bigg] + 2\bigg\{\bP\bigg[c\int_0^T y(s)\d s\bigg]\bigg\}^{1/2} \cr
 \ar\ar\qquad\quad
+\, 2T\int_{(1,\infty)} z \nu(\d z) + 2\bP\bigg[\int_0^Ty(s)\d s\int_{(1,\infty)} z m(\d z)\bigg] \cr
 \ar\ar\qquad\quad
+\, 2\bigg\{\bP\bigg[T\int_{(0,1]} z^2 \nu(\d z) + \int_0^Ty(s)\d s\int_{(0,1]} z^2 m(\d z)\bigg]\bigg\}^{1/2}.
 \eeqnn
Then the desired inequality follows by simple estimates based on \eqref{s5.12}. \qed

\bgcorollary\label{ts7.4} Let $\{(y(t), \mcr{G}_t): t\ge 0\}$ be a c\`{a}dl\`{a}g realization of the CBI-process satisfying $\bP[y(0)]< \infty$. Then the above properties {\rm(3)}, {\rm(4)} and {\rm(5)} hold with the local martingales being martingales.
\edcorollary

\proof Since the arguments are similar, we only give those for {\rm(4)}. Let $f\in C^2[0,\infty)$ and let
 \beqnn
M(t) = f(y(t)) - f(y(0)) - \int_0^t Lf(y(s))\d s, \qquad t\ge 0.
 \eeqnn
By property (4) we know $\{M(t)\}$ is a local martingale. Let $\{\tau_n\}$ be a localization sequence of stopping times for $\{M(t)\}$. For any $t\ge r\ge 0$ and any bounded $\mcr{G}_r$-measurable random variable $F$, we have
 \beqnn
\ar\ar\bP\bigg\{\bigg[f(y(t\land\tau_n)) - f(y(0)) - \int_0^{t\land\tau_n} Lf(y(s))\d s\bigg]F\bigg\} \cr
 \ar\ar\qqquad
=\,\bP\bigg\{\bigg[f(y(r\land\tau_n)) - f(y(0)) - \int_0^{r\land\tau_n} Lf(y(s))\d s\bigg]F\bigg\}.
 \eeqnn
In view of \eqref{s7.1}, there is a constant $C\ge 0$ so that $|Lf(x)|\le C(1+x)$. By Corollary~\ref{ts7.3} we can let $n\to \infty$ and use dominated convergence in the above equality to see $\{M(t)\}$ is a martingale. \qed

Note that property {\rm(4)} implies that the generator of the CBI-process is the closure of the operator $L$ in the sense of Ethier and Kurtz (1986). This explicit form of the generator was first given in Kawazu and Watanabe (1971). The results of Theorem~\ref{ts7.2} were presented for measure-valued processes in El~Karoui and Roelly (1991) and Li (2011).

\newpage

\section{Stochastic equations for CBI-processes}

 \setcounter{equation}{0}

In this section we establish some stochastic equations for the CBI-processes. Suppose that $(\phi,\psi)$ are branching and immigration mechanisms given respectively by \eqref{s2.13} and \eqref{s5.9} with $\nu(\d u)$ satisfying condition \eqref{s5.11}. Let $(P_t)_{t\ge 0}$ be the transition semigroup defined by \eqref{s2.19} and \eqref{s5.10}. In this and the next section, for any $b\ge a\ge 0$ we understand
 \beqnn
\int_a^b = \int_{(a,b]} ~\mbox{and}~ \int_a^\infty = \int_{(a,\infty)}.
 \eeqnn

Let $\{B(t)\}$ be a standard Brownian motion and $\{M(\d s,\d z,\d u)\}$ a Poisson time-space random measure on $(0,\infty)^3$ with intensity $\d sm(\d z)\d u$. Let $\{\eta(t)\}$ be an increasing L\'{e}vy process with $\eta(0)=0$ and with Laplace exponent $\psi(z)= -\log \bP\exp\{-z\eta(1)\}$. We assume that $\{B(t)\}$, $\{M(\d s,\d z,\d u)\}$ and $\{\eta(t)\}$ are defined on a complete probability space and are independent of each other. Consider the stochastic integral equation
 \beqlb\label{s8.1}
y(t)
 \ar=\ar
y(0) + \int_0^t \sqrt{2cy(s-)}\d B(s) - b\int_0^t y(s-)\d s \cr
 \ar\ar\quad
+ \int_0^t\int_0^\infty \int_0^{y(s-)} z \tilde{M}(\d s,\d z,\d u) +
\eta(t),
 \eeqlb
where $\tilde{M}(\d s,\d z,\d u) = M(\d s,\d z,\d u) - \d sm(\d z)\d u$ denotes the compensated measure. We understand the forth term on the right-hand side of \eqref{s8.1} as an integral over the set $\{(s,z,u): 0<s\le t, 0<z<\infty, 0<u\le y(s-)\}$ and give similar interpretations for other stochastic integrals in this section. The reader is referred to Ikeda and Watanabe (1989) and Situ (2005) for the basic theory of stochastic equations.

\bgtheorem\label{ts8.1} A positive c\`{a}dl\`{a}g process $\{y(t): t\ge 0\}$ is a CBI-process with branching and immigration mechanisms $(\phi,\psi)$ given respectively by \eqref{s2.13} and \eqref{s5.9} if and only if it is a weak solution to \eqref{s8.1}. \edtheorem

\proof Suppose that the positive c\`{a}dl\`{a}g process $\{y(t)\}$ is a weak solution to \eqref{s8.1}. By It\^{o}'s formula one can see $\{y(t)\}$ solves the martingale problem \eqref{s7.3}. By Theorem~\ref{ts7.2} we infer that $\{y(t)\}$ is a CBI-process with branching and immigration mechanisms given respectively by \eqref{s2.13} and \eqref{s5.9}. Conversely, suppose that $\{y(t)\}$ is a c\`{a}dl\`{a}g realization of the CBI-process with branching and immigration mechanisms given respectively by \eqref{s2.13} and \eqref{s5.9}. By Theorem~\ref{ts7.2} the process has no negative jumps and the random measure
 \beqnn
N_0(\d s,\d z):= \sum_{s>0}1_{\{\Delta y(s)>0\}}\delta_{(s,\Delta y(s))}(\d s,\d z)
 \eeqnn
has predictable compensator
 \beqnn
\hat{N}_0(\d s,\d z) = y(s-)\d sm(\d z) + \d s\nu(\d z).
 \eeqnn
Moreover, we have
 \beqnn
y(t) \ar=\ar y(0) + t\bigg[\beta + \int_0^\infty u\nu(\d u)\bigg] - \int_0^t
by(s-)\d s \cr
 \ar\ar\qquad
+\, M^c(t) + \int_0^t \int_0^\infty z \tilde{N}_0(\d s,\d z),
 \eeqnn
where $\tilde{N}_0(\d s,\d z) = N_0(\d s,\d z) - \hat{N}_0(\d s,\d z)$ and $t\mapsto M^c(t)$ is a continuous local martingale with quadratic variation $2cy(t-)\d t$. By Theorem~III.7.1$'$ in Ikeda and Watanabe (1989, p.90), on an extension of the original probability space there is a standard Brownian motion $\{B(t)\}$ so that
 \beqnn
M^c(t)= \int_0^t\sqrt{2c y(s-)} \d B(s).
 \eeqnn
By Theorem~III.7.4 in Ikeda and Watanabe (1989, p.93), on a further extension of the probability space we can define independent Poisson random measures $M(\d s,\d z,\d u)$ and $N(\d s,\d z)$ with intensities $\d sm(\d z)\d u$ and $\d s\nu(\d z)$, respectively, so that
 \beqnn
\int_0^t \int_0^\infty z \tilde{N}_0(\d s,\d z)
 =
\int_0^t\int_0^\infty \int_0^{y(s-)} z \tilde{M}(\d s,\d z,\d u) +
\int_0^t\int_0^\infty z \tilde{N}(\d s,\d z).
 \eeqnn
Then $\{y(t)\}$ is a weak solution to \eqref{s8.1}. \qed

\bgtheorem\label{ts8.2} For any initial value $y(0)=x\ge 0$, there is a pathwise unique positive strong solution to \eqref{s8.1}. \edtheorem

\proof By Theorem~\ref{ts8.1} there is a weak solution to \eqref{s8.1}. Then we only need to prove the pathwise uniqueness of the solution; see, e.g., Situ (2005, p.76 and p.104). Suppose that $\{x(t): t\ge 0\}$ and $\{y(t): t\ge 0\}$ are two positive solutions of \eqref{s8.1} with deterministic initial states. By Theorem~\ref{ts8.1}, both of them are CBI-processes. We may assume $x(0)$ and $y(0)$ are deterministic upon taking a conditional probability. In view of \eqref{s5.12}, the processes have locally bounded first moments. Let $\zeta(t)= x(t)-y(t)$ for $t\ge 0$. For each integer $n\ge 0$ define $a_n= \exp\{-n(n+1)/2\}$. Then $a_n\to 0$ decreasingly as $n\to \infty$ and
 \beqnn
\int_{a_n}^{a_{n-1}}z^{-1} \d z = n, \qquad n\ge 1.
 \eeqnn
Let $x\mapsto g_n(x)$ be a positive continuous function supported by $(a_n,a_{n-1})$ so that
 \beqnn
\int_{a_n}^{a_{n-1}}g_n(x)\d x=1
 \eeqnn
and $g_n(x)\le 2(nx)^{-1}$ for every $x>0$. For $n\ge 0$ and $z\in \mbb{R}$ let
 \beqnn
f_n(z)=\int_0^{|z|}\d y\int_0^yg_n(x)\d x.
 \eeqnn
Then $f_n(z)\to |z|$ increasingly as $n\to \infty$. Moreover, we have $|f_n^\prime(z)|\le 1$ and
 \beqlb\label{s8.2}
0\le |z|f_n^{\prime\prime}(z) = |z|g_n(|z|)\le 2/n.
 \eeqlb
For $z,\zeta\in\mbb{R}$ it is easy to see that
 \beqnn
|f_n(\zeta+z) - f_n(\zeta) - zf_n^\prime(\zeta)|
 \le
|f_n(\zeta+z) - f_n(\zeta)| + |zf_n^\prime(\zeta)|\le 2|z|.
 \eeqnn
By Taylor¡¯s expansion, when $z\zeta\ge 0$, there is $\eta$ between $\zeta$ and $\zeta+z$ so that
 \beqnn
|\zeta||f_n(\zeta+z) - f_n(\zeta) - zf_n^\prime(\zeta)|
 \le
|\zeta||f_n^{\prime\prime}(\eta)|z^2/2
 \le
|\eta||f_n^{\prime\prime}(\eta)|z^2/2
 \le
z^2/n.
 \eeqnn
where we used \eqref{s8.2} for the last inequality. It follows that, when $z\zeta\ge 0$,
 \beqlb\label{s8.3}
|\zeta||f_n(\zeta+z) - f_n(\zeta) - zf_n^\prime(\zeta)|
 \le
(2|z\zeta|)\land (z^2/n)
 \le
(1+2|\zeta|)[|z|\land (z^2/n)].
 \eeqlb
{From} \eqref{s8.1} we have
 \beqnn
\zeta(t) \ar=\ar \zeta(0) - b\int_0^t \zeta(s-) \d s +
\sqrt{2c}\int_0^t \big(\sqrt{x(s)}-\sqrt{y(s)}\big)\d B(s) \cr
 \ar\ar\qquad
+ \int_0^t\int_0^\infty \int_{y(s-)}^{x(s-)} z1_{\{\zeta(s-)> 0\}} \tilde{M}(\d s,\d z,\d u) \cr
 \ar\ar\qquad
- \int_0^t\int_0^\infty \int_{x(s-)}^{y(s-)} z1_{\{\zeta(s-)< 0\}} \tilde{M}(\d s,\d z,\d u).
 \eeqnn
By this and It\^o's formula,
 \beqnn
f_n(\zeta(t)) \ar=\ar f_n(\zeta(0)) - b\int_0^t f_n^\prime(\zeta(s)) \zeta(s)\d s
+ c\int_0^t f_n^{\prime\prime}(\zeta(s)) \big[\sqrt{x(s)} - \sqrt{y(s)}\,\big]^2\d s \cr
 \ar\ar
+ \int_0^t\zeta(s)1_{\{\zeta(s)>0\}}\d s \int_0^\infty [f_n(\zeta(s)+z) -
f_n(\zeta(s)) - zf_n^\prime(\zeta(s))] m(\d z) \cr
 \ar\ar
- \int_0^t\zeta(s)1_{\{\zeta(s)<0\}}\d s \int_0^\infty [f_n(\zeta(s)-z) -
f_n(\zeta(s)) + zf_n^\prime(\zeta(s))] m(\d z) \ccr
 \ar\ar
+\, \mbox{mart.}
 \eeqnn
Taking the expectation in both sides of the above and using \eqref{s8.2} and \eqref{s8.3} we see
 \beqlb\label{s8.4}
\mbf{P}[f_n(\zeta(t))]
 \le
f_n(\zeta(0)) + |b|\int_0^t \mbf{P}[|\zeta(s)|] \d s + \varepsilon_n(t),
 \eeqlb
where
 \beqnn
\varepsilon_n(t)= 2cn^{-1}t + \int_0^t(1+2\mbf{P}[|\zeta(s)|])\d s\int_0^\infty [z\land (n^{-1}z^2)] m(\d z).
 \eeqnn
Clearly, we have $\lim_{n\to \infty}\varepsilon_n(t)=0$. Then letting $n\to \infty$ in \eqref{s8.4} we get
 \beqnn
\mbf{P}[|x(t)-y(t)|]
 \le
|x(0)-y(0)| + |b|\int_0^t \mbf{P}[|x(s)-y(s)|] \d s.
 \eeqnn
If $x(0) = y(0)$, we have $\mbf{P}[|x(t)-y(t)|] = 0$ by Gronwall's inequality, and so $\mbf{P}\{x(t)= y(t)\}= 1$ for $t\ge 0$. Then $\mbf{P}\{x(t)= y(t)$ for $t\ge 0\}= 1$ by the right continuity of the processes. That gives the pathwise uniqueness for \eqref{s8.1}. \qed

We can give a formulation of the CBI-process in terms another stochastic integral equation weakly equivalent to \eqref{s8.1}. Let $\{M(\d s,\d z,\d u)\}$ and $\{\eta(s)\}$ be as in \eqref{s8.1}. Let $\{W(\d s,\d u)\}$ be a Gaussian time-space white noise on $(0,\infty)^2$ with intensity $\d s\d u$. We assume $\{W(\d s,\d u)\}$, $\{M(\d s,\d z,\d u)\}$ and $\{\eta(s)\}$ are defined on a complete probability space and are independent of each other. Consider the stochastic integral equation
 \beqlb\label{s8.5}
y(t)
 \ar=\ar
y(0) + \sqrt{2c}\int_0^t\int_0^{y(s-)} W(\d s,\d u) - b\int_0^t y(s-)\d s
\cr
 \ar\ar\qquad
+ \int_0^t\int_0^\infty \int_0^{y(s-)} z \tilde{M}(\d s,\d z,\d u) +
\eta(t).
 \eeqlb
The reader may refer to Li (2011, Section~7.3) and Walsh (1986, Chapter~2) for discussions of stochastic integration with respect to Gaussian time-space white noises.

\bgtheorem\label{ts8.3} A positive c\`{a}dl\`{a}g process $\{y(t): t\ge 0\}$ is a CBI-process with branching and immigration mechanisms $(\phi,\psi)$ given respectively by \eqref{s2.13} and \eqref{s5.9} if and only if it is a weak solution to \eqref{s8.5}. \edtheorem

\proof Suppose that $\{y(t)\}$ is a CBI-process with branching and immigration mechanisms given respectively by \eqref{s2.13} and \eqref{s5.9}. By Theorem~\ref{ts8.1}, the process is a weak solution to \eqref{s8.1}. By El~Karoui and M\'{e}l\'{e}ard (1990, Theorem~III.6), on an extension of the probability space we can define a Gaussian time-space white noise $W(\d s,\d u)$ with intensity $\d s\d u$ so that
 \beqnn
\int_0^t\sqrt{y(s-)} \d B(s)
 =
\int_0^t\int_0^{y(s-)} W(\d s,\d u).
 \eeqnn
Then $\{y(t)\}$ is a weak solution to \eqref{s8.5}. Conversely, suppose that $\{y(t)\}$ is a weak solution to \eqref{s8.5}. By It\^{o}'s formula one can see $\{y(t)\}$ solves the martingale problem \eqref{s7.3}. By Theorem~\ref{ts7.2} we infer that $\{y(t)\}$ is a CBI-process with branching and immigration mechanisms given respectively by \eqref{s2.13} and \eqref{s5.9}. \qed

\bgtheorem\label{ts8.4} Suppose that $\{y_1(t): t\ge 0\}$ and $\{y_2(t): t\ge 0\}$ are two positive solutions to \eqref{s8.5} with $\bP\{y_1(0)\le y_2(0)\}= 1$. Then we have $\bP\{y_1(t)\le y_2(t)$ for all $t\ge 0\}= 1$. \edtheorem

\proof Let $\zeta(t) = y_1(t)-y_2(t)$ for $t\ge 0$. For $n\ge 0$ let $f_n$ be the function defined as in the proof of Theorem~\ref{ts8.2}. Let $h_n(z)= f_n(z\vee 0)$ for $z\in \mbb{R}$. Then $h_n(z)\to z_+:= z\vee 0$ increasingly as $n\to \infty$. {From} \eqref{s8.5} it follows that
 \beqnn
\zeta(t) \ar=\ar \zeta(0) - b\int_0^t \zeta(s-) \d s + \sqrt{2c}\int_0^t \int_{y_2(s-)}^{y_1(s-)} 1_{\{\zeta(s-)>0\}} W(\d s,\d u) \cr
 \ar\ar
-\, \sqrt{2c}\int_0^t \int_{y_1(s-)}^{y_2(s-)} 1_{\{\zeta(s-)<0\}} W(\d s,\d u) \cr
 \ar\ar
+ \int_0^t\int_0^\infty \int_{y_2(s-)}^{y_1(s-)} 1_{\{\zeta(s-)>0\}}z \tilde{M}(\d s,\d z,\d u) \cr
 \ar\ar
- \int_0^t\int_0^\infty \int_{y_1(s-)}^{y_2(s-)} 1_{\{\zeta(s-)<0\}}z \tilde{M}(\d s,\d z,\d u).
 \eeqnn
Since $h_n(z)=0$ for $z\le 0$, by It\^o's formula we have
 \beqnn
h_n(\zeta(t)) \ar=\ar -\, b\int_0^t h_n^\prime(\zeta(s-))\zeta(s-) \d s + c\int_0^t h_n^{\prime\prime}(\zeta(s-))|\zeta(s-)| \d s \cr
 \ar\ar
+ \int_0^t\zeta(s-)1_{\{\zeta(s-)>0\}}\d s \int_0^\infty \Big[h_n(\zeta(s-)+z) - h_n(\zeta(s-)) \cr
 \ar\ar
-\, zh_n^\prime(\zeta(s-))\Big] m(\d z) - \int_0^t\zeta(s-)1_{\{\zeta(s-)<0\}}\d s \int_0^\infty \Big[h_n(\zeta(s-)-z) \cr
 \ar\ar
-\, h_n(\zeta(s-)) + zh_n^\prime(\zeta(s-))\Big] m(\d z) + \mbox{local mart.} \cr
 \ar=\ar
-\, b\int_0^t h_n^\prime(\zeta(s-))\zeta(s-)_+ \d s + c\int_0^t h_n^{\prime\prime}(\zeta(s-))\zeta(s-)_+ \d s \cr
 \ar\ar
+ \int_0^t\zeta(s-)_+\d s \int_0^\infty \Big[h_n(\zeta(s-)+z) - h_n(\zeta(s-)) \cr
 \ar\ar
-\, zh_n^\prime(\zeta(s-))\Big] m(\d z) + \mbox{local mart.}
 \eeqnn
For any $k\ge 1$ define $\tau_k= \inf\{t\ge 0: \zeta(t)_+\ge k\}$. Taking the expectation in the above equality at time $t\land\tau_k$ and using \eqref{s8.2} and \eqref{s8.3} we have
 \beqnn
\bP[h_n(\zeta(t\land\tau_k))]
 \le
|b|\bP\bigg[\int_0^{t\land\tau_k} \zeta(s-)_+ \d s\bigg] + \varepsilon_n(t),
 \eeqnn
where
 \beqnn
\varepsilon_n(t)
 =
2cn^{-1}t + \bP\bigg[\int_0^{t\land\tau_k}(1+2\zeta(s-)_+)\d s\bigg]\int_0^\infty (z\land n^{-1}z^2) m(\d z).
 \eeqnn
Then we let $n\to \infty$ to obtain
 \beqnn
\bP[\zeta(t\land\tau_k)_+]
 \le
|b|\bP\bigg[\int_0^{t\land\tau_k}\zeta(s-)_+ \d s\bigg]
 \le
|b|\int_0^t \bP[\zeta(s\land\tau_k)_+] \d s.
 \eeqnn
By Gronwall's inequality, for each $t\ge 0$ we have
 \beqnn
\bP[(y_1(t\land\tau_k) - y_2(t\land\tau_k))_+]= \bP[\zeta(t\land\tau_k)_+]= 0.
 \eeqnn
By letting $k\to \infty$ and using Fatou's lemma we see $\bP[(y_1(t)-y_2(t))_+]= 0$ for $t\ge 0$, and so $\bP\{y_1(t)\le y_2(t)$ for all $t\ge 0\}= 1$ by the right continuity of the processes. \qed

\bgtheorem\label{ts8.5} For any initial value $y(0)=x\ge 0$, there is a pathwise unique positive strong solution to \eqref{s8.5}. \edtheorem

\proof By Theorem~\ref{ts8.3} there is a weak solution $\{y(t)\}$ to \eqref{s8.5}. The pathwise uniqueness of the solution follows from Theorem~\ref{ts8.4}. Then $\{y(t)\}$ is a strong solution to \eqref{s8.5}. See, e.g., Situ (2005, p.76 and p.104). \qed

{From} \eqref{s8.1} or \eqref{s8.5} we see that the immigration of the CBI-process $\{y(t)\}$ is represented by the increasing L\'{e}vy process $\{\eta(t)\}$. By the L\'{e}vy--It\^{o} decomposition, there is a Poisson time-space random measure $\{N(\d s,\d z)\}$ with intensity $\d s\nu(\d z)$ such that
 \beqnn
\eta(t) = \beta t + \int_0^t\int_0^\infty z N(\d s,\d z), \qquad t\ge 0.
 \eeqnn
Then the immigration of $\{y(t)\}$ involves two parts: the \index{continuous part of immigration} \textit{continuous part} determined by the drift coefficient $\beta$ and the \index{discontinuous part of immigration} \textit{discontinuous part} given by the Poisson random measure $\{N(\d s,\d z)\}$.

Now let us consider a special CBI-process. Let $c,q\ge 0$, $b\in \mbb{R}$ and $1<\alpha <2$ be given constants. Let $\{B(t)\}$ be a standard Brownian motion. Let $\{z(t)\}$ be a spectrally positive $\alpha$-stable L\'{e}vy process with L\'{e}vy measure
 \beqnn
\gamma(\d z):= (\alpha-1) \Gamma(2-\alpha)^{-1} z^{-1-\alpha}\d z, \qquad z>0
 \eeqnn
and $\{\eta(t)\}$ an increasing L\'{e}vy process with $\eta(0)=0$ and with Laplace exponent $\psi$. We assume that $\{B(t)\}$, $\{z(t)\}$ and $\{\eta(t)\}$ are defined on a complete probability space and are independent of each other. Consider the stochastic differential equation
 \beqlb\label{s8.6}
\d y(t) = \sqrt{2cy(t-)}\d B(t) + \sqrt[\alpha]{\alpha qy(t-)}\d z(t) - by(t-)\d t
+ \d \eta(t),
 \eeqlb

\bgtheorem\label{ts8.6} A positive c\`{a}dl\`{a}g process $\{y(t): t\ge 0\}$ is a CBI-process with branching mechanism $\phi(z) = bz + cz^2 + qz^\alpha$ and immigration mechanism $\psi$ given by \eqref{s5.9} if and only if it is a weak solution to \eqref{s8.6}. \edtheorem

\proof Suppose that $\{y(t)\}$ is a weak solution to \eqref{s8.6}. By It\^{o}'s formula one can see that $\{y(t)\}$ solves the martingale problem \eqref{s7.3} associated with the generator $L$ defined by \eqref{s7.1} with $m(\d z)= \alpha q\gamma(\d z)$. Then $\{y(t)\}$ is a CBI-process with branching mechanism $\phi(z) = bz + cz^2 + qz^\alpha$ and immigration mechanism $\psi$ given by \eqref{s5.9}. Conversely, suppose that $\{y(t)\}$ is a CBI-process with branching mechanism $\phi(z) = bz + cz^2 + qz^\alpha$ and immigration mechanism $\psi$ given by \eqref{s5.9}. Then $\{y(t)\}$ is a weak solution to \eqref{s8.1} with $\{M(\d s,\d z,\d u)\}$ being a Poisson random measure on $(0,\infty)^3$ with intensity $\alpha q\d s\gamma(\d z)\d u$. Let us assume $q>0$, for otherwise the proof is easier. Define the random measure $\{N_0(\d s,\d z)\}$ on $(0,\infty)^2$ by
 \beqnn
N_0((0,t]\times B)
 \ar=\ar
\int_0^t\int_0^\infty\int_0^{y(s-)} 1_{\{y(s-)>0\}}1_B\Big(\frac{z}
{\sqrt[\alpha]{\alpha qy(s-)}}\Big) M(\d s,\d z,\d u) \cr
 \ar\ar
+ \int_0^t\int_0^\infty\int_0^{1/\alpha q} 1_{\{y(s-)=0\}}1_B(z) M(\d s,\d z,\d
u).
 \eeqnn
It is easy to compute that $\{N_0(\d s,\d z)\}$ has predictable
compensator
 \beqnn
\hat{N}_0((0,t]\times B)
 \ar=\ar
\int_0^t\int_0^\infty 1_{\{y(s-)>0\}}1_B\bigg(\frac{z}
{\sqrt[\alpha]{\alpha qy(s-)}}\bigg)\frac{\alpha qy(s-)(\alpha-1)\d s\d z}{\Gamma(2-\alpha)z^{1+\alpha}} \cr
 \ar\ar
+ \int_0^t\int_0^\infty 1_{\{y(s-)=0\}}1_B(z)\frac{(\alpha-1)\d s\d z}
{\Gamma(2-\alpha)z^{1+\alpha}} \cr
 \ar=\ar
\int_0^t\int_0^\infty 1_B(z)\frac{(\alpha-1)\d s\d z} {\Gamma(2-\alpha)z^{1+\alpha}}.
 \eeqnn
Thus $\{N_0(\d s,\d z)\}$ is a Poisson random measure with intensity $\d s\gamma(\d z)$; see, e.g., Theorem~III.7.4 in Ikeda and Watanabe (1989, p.93). Now define the L\'{e}vy processes
 \beqnn
z(t) = \int_0^t\int_0^\infty z \tilde{N}_0(\d s,\d z)
 ~~\mbox{and}~~
\eta(t) = \beta t + \int_0^t\int_0^\infty z N(\d s,\d z),
 \eeqnn
where $\tilde{N}_0(\d s,\d z) = N_0(\d s,\d z) - \hat{N}_0(\d s,\d z)$.
It is easy to see that
 \beqnn
\int_0^t\sqrt[\alpha]{\alpha qy(s-)}\d z(s)
 \ar=\ar
\int_0^t\int_0^\infty \sqrt[\alpha]{\alpha qy(s-)}\,z \tilde{N}_0(\d s,\d z) \cr
 \ar=\ar
\int_0^t\int_0^\infty \int_0^{y(s-)} z \tilde{M}(\d s,\d z,\d u).
 \eeqnn
Then $\{y(t)\}$ is a weak solution to \eqref{s8.6}. \qed

\bgtheorem\label{ts8.7} For any initial value $y(0)=x\ge 0$, there is a pathwise unique positive strong solution to \eqref{s8.6}. \edtheorem

\proof By Theorem~\ref{ts8.6} there is a weak solution $\{y(t)\}$ to \eqref{s8.6}, so it suffices to prove the pathwise uniqueness of the solution. We first recall that the one-sided $\alpha$-stable process $\{z(t)\}$ can be represented as
 \beqnn
z(t) = \int_0^t\int_0^\infty z\tilde{M}(\d s,\d z),
 \eeqnn
where $M(\d s,\d z)$ is a Poisson random measure on $(0,\infty)^2$ with intensity $\d s\gamma(\d z)$. Let
 \beqnn
z_1(t) = \int_0^t\int_0^1 z\tilde{M}(\d s,\d z)
 \quad\mbox{and}\quad
z_2(t) = \int_0^t\int_1^\infty z M(\d s,\d z).
 \eeqnn
Since $t\mapsto z_2(t)$ has at most finitely many jumps in each bounded interval, we only need to prove the pathwise uniqueness of
 \beqlb\label{s8.7}
\d y(t) \ar=\ar \sqrt{2cy(t-)}\d B(t) + \sqrt[\alpha]{\alpha qy(t-)}\d z_1(t) - by(t-)\d t \ccr
 \ar\ar\qquad
-\, \alpha^{-1}(\alpha-1) \Gamma(2-\alpha)^{-1} \sqrt[\alpha]{\alpha qy(t-)}\d t + \d \eta(t),
 \eeqlb
Suppose that $\{x(t)\}$ and $\{y(t)\}$ are two positive solutions to \eqref{s8.7} with deterministic initial values. Let $\zeta_\theta(t) = \sqrt[\theta]{x(t)} - \sqrt[\theta]{y(t)}$ for $0<\theta\le 2$ and $t\ge 0$. Then we have
 \beqnn
\d\zeta_1(t) \ar=\ar \sqrt{2c}\zeta_2(t-)\d B(t) + \sqrt[\alpha]{\alpha q} \zeta_\alpha(t-)\d z_1(t) - b\zeta_1(t-) \d t \ccr
 \ar\ar\qquad
-\, \alpha^{-1}(\alpha-1) \Gamma(2-\alpha)^{-1}\sqrt[\alpha]{\alpha q} \zeta_\alpha(t-)\d t.
 \eeqnn
For $n\ge 0$ let $f_n$ be the function defined as in the proof of Theorem~\ref{ts8.2}. By It\^o's formula,
 \beqlb\label{s8.8}
f_n(\zeta_1(t)) \ar=\ar f_n(\zeta_1(0)) + c\int_0^t f_n^{\prime\prime}(\zeta_1(s-))\zeta_1(s-)^2\d s - b\int_0^t
f_n^\prime(\zeta_1(s-)) \zeta_1(s-) \d s \cr
 \ar\ar
-\, \alpha^{-1}(\alpha-1) \Gamma(2-\alpha)^{-1}\sqrt[\alpha]{\alpha q}\int_0^tf_n^\prime(\zeta_1(s-))\zeta_\alpha(s-)\d s \cr
 \ar\ar
+\, \int_0^t ds \int_0^1 \Big[f_n(\zeta_1(s-)+\sqrt[\alpha]{\alpha q}\zeta_\alpha(s-)z) - f_n(\zeta_1(s-)) \cr
 \ar\ar
-\, \sqrt[\alpha]{\alpha q}\zeta_\alpha(s-) zf_n^\prime(\zeta_1(s-))\Big] \gamma(\d z) + \mbox{local mart.}
 \eeqlb
For any $k\ge 1+x(0)\vee y(0)$ let $\tau_k = \inf\{s\ge 0: x(s)\ge k$ or $y(s)\ge k\}$. For $0\le t\le \tau_k$ we have $|\zeta_1(t-)|\le k$, $|\zeta_\alpha(t-)|\le \sqrt[\alpha]{k}$ and
 \beqnn
|\zeta_1(t)|\le |\zeta_1(t-)| + |\zeta_1(t) - \zeta_1(t-)|\le k + \sqrt[\alpha]{\alpha qk}.
 \eeqnn
By Taylor's expansion, there exists $0<\xi<z$ so that
 \beqnn
\ar\ar [f_n(\zeta_1(s-) + \sqrt[\alpha]{\alpha q}\zeta_\alpha(s-)z) - f_n(\zeta_1(s-)) - \sqrt[\alpha]{\alpha q}\zeta_\alpha(s-) zf_n^\prime(\zeta_1(s-))] \ccr
 \ar\ar\qquad
=
2^{-1}(\alpha q)^{2/\alpha}f_n^{\prime\prime}(\zeta_1(s-) + \sqrt[\alpha]{\alpha q}\zeta_\alpha(s-)\xi)\zeta_\alpha(s-)^2z^2
\ccr
 \ar\ar\qquad
\le
2^{-1}(\alpha q)^{2/\alpha}\sqrt[\alpha]{k}f_n^{\prime\prime}(\zeta_1(s-) + \sqrt[\alpha]{\alpha q}\zeta_\alpha(s-)\xi)|\zeta_1(s-) + \sqrt[\alpha]{\alpha q}\zeta_\alpha(s-)\xi|z^2 \ccr
 \ar\ar\qquad
\le n^{-1}(\alpha q)^{2/\alpha}\sqrt[\alpha]{k}z^2,
 \eeqnn
where we have used \eqref{s8.2} and the fact $\zeta_1(s-)\zeta_\alpha(s-)\ge 0$. Taking the expectation in both sides of \eqref{s8.8} gives
 \beqnn
\bP[f_n(\zeta_1(t\land\tau_k))]
 \ar\le\ar
\bP[f_n(\zeta_1(0))] + |b|\int_0^t \bP[|\zeta_1(s\land\tau_k)|] \d s +
2cn^{-1}kt \cr
 \ar\ar
+\, n^{-1}(\alpha q)^{2/\alpha}\sqrt[\alpha]{k}\int_0^t\d s\int_0^1z^2 \gamma(\d z) .
 \eeqnn
Now, if $x(0)=y(0)$, we can let $n\to \infty$ in the inequality above to get
 \beqnn
\bP[|x(t\land\tau_k)-y(t\land\tau_k)|]
 \le
|b|\int_0^t \bP[|x(s\land\tau_k)-y(s\land\tau_k)|] \d s.
 \eeqnn
Then $\bP[|x(t\land\tau_k) - y(t\land\tau_k)|] = 0$ for $t\ge 0$ by Gronwall's inequality. By letting $k\to \infty$ and using Fatou's lemma we obtain the pathwise uniqueness for \eqref{s8.6}. \qed

\bgexample\label{es8.1} The stochastic integral equation \eqref{s8.5} can be thought as a continuous time-space counterpart of the definition \eqref{s5.1} of the GWI-process. In fact, assuming $\mu = \bE(\xi_{1,1})< \infty$, from \eqref{s5.1} we have
 \beqlb\label{s8.9xx}
y(n)-y(n-1)
 =
\sum_{i=1}^{y(n-1)}(\xi_{n,i}-\mu) - (1-\mu)y(n-1) + \eta_n.
 \eeqlb
It follows that
 \beqlb\label{s8.9}
y(n) - y(0) = \sum_{k=1}^n\sum_{i=1}^{y(k-1)}(\xi_{k,i}-\mu) - (1-\mu)\sum_{k=1}^ny(k-1) +
\sum_{k=1}^n\eta_k.
 \eeqlb
The exact continuous time-state counterpart of \eqref{s8.9} would be the stochastic integral equation
 \beqlb\label{s8.10}
y(t) = y(0) + \int_0^t\int_0^\infty \int_0^{y(s-)} \xi \tilde{M}(\d s,\d\xi,\d
u) - \int_0^t by(s)\d s + \eta(t),
 \eeqlb
which is a typical special form of \eqref{s8.5}; see Bertoin and Le~Gall (2006) and Dawson and Li (2006). Here the $\xi$'s selected by the Poisson random measure $M(\d s,\d\xi,\d u)$ are distributed in a i.i.d.\ fashion and the compensation of the measure corresponds to the centralization in \eqref{s8.9}. The increasing L\'{e}vy process $t\mapsto \eta(t)$ in \eqref{s8.10} corresponds to the increasing random walk $n\mapsto \sum_{k=1}^n\eta_k$ in \eqref{s8.9}. The additional term in \eqref{s8.5} involving the stochastic integral with respect to the Gaussian white noise is just a continuous time-space parallel of that with respect to the compensated Poisson random measure. \edexample

\bgexample\label{es8.2} The stochastic differential equation \eqref{s8.6} captures the structure of the CBI-process in a typical special case. Let $1<\alpha\le 2$. Under the condition $\mu := \bE(\xi_{1,1})< \infty$, from \eqref{s8.9xx} we have
 \beqnn
y(n) - y(n-1) = \sqrt[\alpha]{y(n-1)} \sum_{i=1}^{y(n-1)} \frac{\xi_{n,i}-\mu} {\sqrt[\alpha]{y(n-1)}} - (1-\mu)y(n-1) + \eta_n.
 \eeqnn
Observe that the partial sum on the right-hand side corresponds to a one-sided $\alpha$-stable type central limit theorem. Then a continuous time-state counterpart of the above equation would be
 \beqlb\label{s8.11}
\d y(t)
 =
\sqrt[\alpha]{\alpha qy(t-)}\d z(t) - by(t)\d t + \beta\d t, \qquad t\ge 0,
 \eeqlb
where $\{z(t): t\ge 0\}$ is a standard Brownian motion if $\alpha=2$ and a spectrally positive $\alpha$-stable L\'{e}vy process with L\'{e}vy measure $(\alpha-1) \Gamma(2-\alpha)^{-1} z^{-1-\alpha}\d z$ if $1<\alpha<2$. This is a typical special form of \eqref{s8.6}. \edexample

\bgexample\label{es8.3} When $\alpha=2$ and $\beta=0$, the CB-process defined by \eqref{s8.11} is a diffusion process, which is known as \index{Feller's branching diffusion} \textit{Feller's branching diffusion}. This process was first studied by Feller (1951). \edexample

\bgexample\label{es8.4} In the special case of $\alpha=2$, the CBI-process defined by \eqref{s8.11} is known in mathematical finance as the \index{Cox--Ingersoll--Ross model} \textit{Cox--Ingersoll--Ross model} \index{CIR-model} (CIR-model), which was used by Cox et al.\ (1985) to describe the evolution of interest rates. The asymptotic behavior of the estimators of the parameters in the CIR-model was studied by Overbeck and Ryd\'en (1997). In the general case, the solution to \eqref{s8.11} is called a \index{$\alpha$-stable Cox--Ingersoll--Ross model} \textit{$\alpha$-stable Cox--Ingersoll--Ross model} \index{$\alpha$-stable CIR-model} ($\alpha$-stable CIR-model); see, e.g., Jiao et al.\ (2017) and Li and Ma (2015). \edexample

As a simple application of the stochastic equation \eqref{s8.1} or \eqref{s8.5}, we can give a simple derivation of the joint Laplace transform of the CBI-process and its positive integral functional. The next theorem extends the results in Section~4.

\bgtheorem\label{ts8.8} Let $Y = (\Omega, \mcr{F}, \mcr{F}_t, y(t), \bP_x)$ be a Hunt realization of the CBI-process. Then for $t, \lambda, \theta\ge 0$ we have
 \beqnn
\bP_x\exp\bigg\{-\lambda y(t) - \theta\int_0^t y(s)\d s\bigg\}
 =
\exp\bigg\{-xv(t) - \int_0^t \psi(v(s))\d s\bigg\},
 \eeqnn
where $t\mapsto v(t) = v(t,\lambda,\theta)$ is the unique positive
solution to \eqref{s4.14}. \edtheorem

\proof We can construct the process $\{y(t): t\ge 0\}$ as the solution to \eqref{s8.1} or \eqref{s8.5} with $y(0)=x\ge 0$. Let
 \beqnn
z(t) = \int_0^t y(s)\d s, \qquad t\ge 0.
 \eeqnn
Consider a function $G=G(t,y,z)$ on $[0,\infty)^3$ with bounded continuous derivatives up to the first order relative to $t\ge 0$ and $z\ge 0$ and up to the second order relative to $x\ge 0$. By It\^o's formula,
 \beqnn
G(t,y(t),z(t)) \ar=\ar G(0,y(0),0) + \mbox{local mart.} + \int_0^t
\bigg\{G^\prime_t(s,y(s),z(s)) \ccr
 \ar\ar
+\, y(s)G^\prime_z(s,y(s),z(s)) + [\beta - by(s)] G^\prime_y(s,y(s),z(s)) \ccr
 \ar\ar
+\, cy(s)G^{\prime\prime}_{yy}(s,y(s),z(s)) + y(s)\int_0^\infty \Big[G(s,y(s)+z,z(s)) \cr
 \ar\ar
-\, G(s,y(s),z(s)) - z G^\prime_y(s,y(s),z(s))\Big]m(\d z)\bigg\}\d s
\cr
 \ar\ar
+ \int_0^t\d s\int_0^\infty\Big[G(s,y(s)+z,z(s)) -
G(s,y(s),z(s))\Big]\nu(\d z).
 \eeqnn
We can apply the above formula to the function
 \beqnn
G_T(t,y,z) = \exp\bigg\{-v(T-t)x - \theta z - \int_0^{T-t}\psi(v(s))\d s\bigg\}.
 \eeqnn
Using \eqref{s4.14} we see $t\mapsto G_T(t\land T,y(t\land T),z(t\land T))$ is a local martingale, and hence a martingale by the boundedness. {From} the relation $\bP_x [G_T(t,y(t),z(t))]= G_T(0,x,0)$ with $T=t$ we get the desired result. \qed

The existence and uniqueness of strong solution to \eqref{s8.1} were first established in Dawson and Li (2006). The moment condition \eqref{s5.11} was removed in Fu and Li (2010). A stochastic flow of discontinuous CB-processes with critical branching mechanism was constructed in Bertoin and Le~Gall (2006) by using weak solutions of a special case of \eqref{s8.1}. The existence and uniqueness of strong solution to \eqref{s8.6} were proved in Fu and Li (2010) and those for \eqref{s8.5} were given in Dawson and Li (2012) and Li and Ma (2008). The results of Bertoin and Le~Gall (2006) were extended to flows of CBI-processes in Dawson and Li (2012) and Li (2014) using strong solutions. Although the study of branching processes has a long history, the stochastic equations \eqref{s8.1}, \eqref{s8.5} and \eqref{s8.6} were not established until the works mentioned above.

A natural generalization of the CBI-process is the so-called \index{affine Markov process} \textit{affine Markov process}; see Duffie et al.\ (2003) and the references therein. Those authors defined the regularity property of affine processes and gave a number of characterizations of those processes under the regularity assumption. By a result of Kawazu and Watanabe (1971), a stochastically continuous CBI-process is automatically regular. Under the first moment assumption, the regularity of affine processes was proved in Dawson and Li (2006). The regularity problem was settled in Keller-Ressel et al.\ (2011), where it was proved that any stochastically continuous affine process is regular. This problem is related to Hilbert's fifth problem; see Keller-Ressel et al.\ (2011) for details.

\newpage

\section{Local and global maximal jumps}

 \setcounter{equation}{0}

In this section, we use stochastic equations of the CB- and CBI-processes to derive several characterizations of the distributions of their local and global maximal jumps. Let us consider a branching mechanism $\phi$ given by \eqref{s2.13}. Let $\{B(t)\}$ be a standard Brownian motion and $\{M(\d s,\d z,\d u)\}$ a Poisson time-space random measure on $(0,\infty)^3$ with intensity $\d sm(\d z)\d u$. By Theorem~\ref{ts8.2}, for any $x\ge 0$, there is a pathwise unique positive strong solution to
 \beqlb\label{s9.1}
x(t) \ar=\ar x + \int_0^t \sqrt{2cx(s-)}\d B(s) - b\int_0^t x(s-)\d s \cr
 \ar\ar\qquad
+ \int_0^t\int_0^\infty \int_0^{x(s-)} z \tilde{M}(\d s,\d z,\d u).
 \eeqlb
By Theorem~\ref{ts8.1}, the solution $\{x(t): t\ge 0\}$ is a CB-process with branching mechanism $\phi$. For $t\ge 0$ and $r>0$ let
 \beqnn
N_r(t)= \int_0^t\int_r^\infty \int_0^{x(s-)} M(\d s,\d z,\d u),
 \eeqnn
which denotes the number of jumps with sizes in $(r,\infty)$ of the trajectory $t\mapsto x(t)$ on the interval $(0,t]$. By \eqref{s3.5} we have
 \beqnn
\bP[N_r(t)]
 =
m(r,\infty)\bP\bigg[\int_0^t x(s)\d s\bigg]
 =
xb^{-1}(1-\e^{-bt})m(r,\infty),
 \eeqnn
where $b^{-1}(1-\e^{-bt}) = t$ for $b=0$ by convention. In particular, we have $\bP\{N_r(t)< \infty\} = 1$. For $r>0$ we can define another branching mechanism by
 \beqlb\label{s9.2}
\phi_r(z) = b_rz + cz^2 + \int_0^r (\e^{-zu}-1+zu) m(\d u),
 \eeqlb
where
 \beqnn
b_r = b + \int_r^\infty u m(\d u).
 \eeqnn
For $\theta\ge 0$ let $t\mapsto u(t,\theta)$ be the unique positive solution to \eqref{s4.17}. Let $t\mapsto u_r(t,\theta)$ be the unique positive solution to
 \beqlb\label{s9.3}
\frac{\partial}{\partial t}u(t,\theta)
 =
\theta - \phi_r(u(t,\theta)),
 \quad
u(0,\theta) = 0.
 \eeqlb
The following theorem gives a characterization of the distribution of the local maximal jump of the CB-process:

\bgtheorem\label{ts9.1} Let $\Delta x(t)= x(t)-x(t-)$ for $t\ge 0$. Then for any $r>0$ we have
 \beqnn
\bP_x\bigg\{\max_{0<s\le t} \Delta x(s)\leq r\bigg\}
 =
\exp\{-x u_r(t)\},
 \eeqnn
where $u_r(t)= u_r(t,m(r,\infty))$. \edtheorem

\proof Let $M_r(\d s,\d z,\d u)$ and $M^r(\d s,\d z,\d u)$ denote the restrictions of $M(\d s,\d z,\d u)$ to $(0,\infty)\times (0,r]\times(0,\infty)$ and $(0,\infty)\times (r,\infty) \times(0,\infty)$, respectively. We can rewrite \eqref{s9.1} into
 \beqnn
x(t) \ar=\ar x + \int_0^t \sqrt{2cx(s-)}\d B(s) + \int_0^t\int_0^r \int_0^{x(s-)} z \tilde{M}_r(\d s,\d z,\d u) \cr
 \ar\ar\quad
- \int_0^t b_rx(s-)\d s + \int_0^t\int_r^\infty \int_0^{x(s-)} z M^r(\d s,\d z,\d u),
 \eeqnn
where the last term collects the jumps with sizes in $(r,\infty)$ of $\{x(t)\}$. Let $\{x_r(t)\}$ be the unique positive strong solution to
 \beqnn
x_r(t) \ar=\ar x - \int_0^t b_rx_r(s-)\d s + \int_0^t \sqrt{2cx_r(s-)}\d B(s) \cr
 \ar\ar\qquad
+ \int_0^t\int_0^r \int_0^{x_r(s-)} z \tilde{M}_r(\d s,\d z,\d u).
 \eeqnn
Then $\{x_r(t)\}$ is a CB-process with branching mechanism $\phi_r$. Let $\tau_r = \inf\{s\ge 0: \Delta x(s)> r\}$. We have $x_r(s) = x(s)$ for $0\le s< \tau_r$ and
 \beqnn
\bigg\{\max_{0<s\le t} \Delta x(s)\leq r\bigg\}
 \ar=\ar\
\bigg\{\int_0^t\int_r^\infty\int_0^{x(s-)} M^r(\d s,\d z,\d u) = 0\bigg\} \cr
 \ar=\ar
\bigg\{\int_0^t\int_r^\infty\int_0^{x_r(s-)} M^r(\d s,\d z,\d u) = 0\bigg\}.
 \eeqnn
Since the strong solution $\{x_r(t)\}$ is progressively measurable with respect to the filtration generated by $\{B(t)\}$ and $\{M_r(\d s,\d z,\d u)\}$, it is independent of $\{M^r(\d s,\d z,\d u)\}$. Then $\{M^r(\d s,\d z,\d u)\}$ is still a Poisson random measure conditionally upon $\{x_r(t)\}$. It follows that
 \beqnn
\bP_x\bigg\{\max_{0<s\le t} \Delta x(s)\leq r\bigg\}
 =
\bP_x\bigg[\exp\bigg\{-m(r,\infty)\int_0^t x_r(s) \d s\bigg\}\bigg].
 \eeqnn
Then the desired result follows by Corollary~\ref{ts4.4}. \qed

\bgcorollary\label{ts9.2} Suppose that the measure $m(\d u)$ has unbounded support. Then we have,
as $r\to\infty$,
 \beqnn
\bP_x\bigg\{\max_{0<s\le t} \Delta x(s)>r\bigg\}
 \sim
xb^{-1}(1-\e^{-bt})m(r,\infty).
 \eeqnn
\edcorollary

\proof Recall that $t\mapsto u(t,\theta)$ is defined by \eqref{s4.17} and $t\mapsto u_r(t,\theta)$ is defined by \eqref{s9.3}. It is easy to see that $u(t,0) = u_r(t,0) = 0$.
Moreover, by \eqref{s4.17} we have
 \beqnn
\frac{\partial}{\partial t}\frac{\partial}{\partial\theta}u(t,0)
 =
1 - b\frac{\partial}{\partial\theta}u(t,0),
 \quad
\frac{\partial}{\partial\theta}u(0,0) = 0.
 \eeqnn
We can solve the above equation to get
 \beqlb\label{s9.4}
\frac{\partial}{\partial\theta}u(t,0)
 =
b^{-1}(1-\e^{-bt}),
 \eeqlb
where $b^{-1}(1-\e^{-bt})= t$ for $b=0$ by convention. Similarly we have
 \beqlb\label{s9.5}
\frac{\partial}{\partial\theta}u_r(t,0)
 =
b_r^{-1}(1-\e^{-b_r t}).
 \eeqlb
By Theorem~\ref{ts9.1} it follows that
 \beqnn
\bP_x\bigg\{\max_{0<s\le t} \Delta x(s)> r\bigg\}
 =
1-\exp\{-x u_r(t,m(r,\infty))\},
 \eeqnn
For $r>q>0$, we have obviously $\phi\leq \phi_r\leq \phi_q$. By Corollary~\ref{ts4.6} we see \beqnn
u_q(t,m(r,\infty))\leq u_r(t,m(r,\infty))\leq u(t,m(r,\infty)).
 \eeqnn
It follows that
 \beqnn
1-\exp\{-xu_q(t,m(r,\infty))\}
 \ar\leq\ar
\bP_x\bigg\{\max_{0<s\le t} \Delta x(s)> r\bigg\} \cr
 \ar\leq\ar
1-\exp\{-x u(t,m(r,\infty))\}.
 \eeqnn
By \eqref{s9.4} and \eqref{s9.5}, as $r\to\infty$ we have
 \beqnn
1-\exp\{-xu(t,m(r,\infty))\}
 \ar\sim\ar
xu(t,m(r,\infty)) \ccr
 \ar\sim\ar
xb^{-1}(1-\e^{-bt})m(r,\infty),
 \eeqnn
and
 \beqnn
1-\exp\{-xu_q(t,m(r,\infty))\}
 \ar\sim\ar
x u_q(t,m(r,\infty)) \ccr
 \ar\sim\ar
xb_q^{-1}(1-\e^{-b_qt})m(r,\infty).
 \eeqnn
The proof is completed as we notice $\lim_{q\to \infty}b_q = b$. \qed

We can also give some characterizations of the global maximal jump of the CB-process. Let $\phi_r^{-1}(\theta):= \inf\{z\ge 0: \phi_r(z)> \theta\}$ for $\theta\ge 0$. It is easy to see that $\phi_r^{-1}(m(r,\infty))\to 0$ as $r\to \infty$ if and only if $b\ge 0$. Let $\phi^\prime(\infty)$ be given by \eqref{s3.2}. By Theorems~\ref{ts4.8} and~\ref{ts9.1} we have:

\bgcorollary\label{ts9.3} Suppose that $\phi^\prime(\infty)> 0$. Then for any $r>0$ with $m(r,\infty)>0$ we have
 \beqnn
\bP_x\Big\{\sup_{s>0} \Delta x(s)\leq r\Big\}
 =
\exp\{-x \phi_r^{-1}(m(r,\infty))\}.
 \eeqnn
\edcorollary

\bgcorollary\label{ts9.4} Suppose that $b>0$ and the measure $m(\d u)$ has unbounded support. Then as $r\to\infty$ we have
 \beqnn
\bP_x\Big\{\sup_{s>0} \Delta x(s)> r\Big\}
 \sim
xb^{-1}m(r,\infty).
 \eeqnn
\edcorollary

The results on local maximal jumps obtained above can be generalized to the case of a CBI-process. Let $(\phi,\psi)$ be the branching and immigration mechanisms given respectively by \eqref{s2.13} and \eqref{s5.9} with $\nu(\d u)$ satisfying \eqref{s5.11}. Let $\{y(t): t\ge 0\}$ be the CBI-process defined by \eqref{s8.1} with $y(0)= x\ge 0$. For $r>0$ let
 \beqnn
\psi_r(z) = \beta z + \int_0^r (1-\e^{-zu}) \nu(\d u).
 \eeqnn
Based on Theorem~\ref{ts8.8}, the following theorem can be proved by modifying the arguments in the proof of Theorem~\ref{ts9.1}:

\bgtheorem\label{ts9.5} Let $\Delta y(t)= y(t)-y(t-)$ for $t\ge 0$. Then for any $r>0$ we have
 \beqnn
\ar\ar\bP_x\bigg\{\max_{0<s\le t} \Delta y(s)\leq r\bigg\} \cr
 \ar\ar\qquad
= \exp\bigg\{-xu_r(t) - \nu(r,\infty)t -\int_0^t \psi_r(u_r(s))\d s\bigg\},
 \eeqnn
where $u_r(t) = u_r(t,m(r,\infty))$. \edtheorem

The results given in this section were adopted from He and Li (2016). We refer the reader to Bernis and Scotti (2018+) and Jiao et al.\ (2017) for more careful analysis of the jumps of CBI-processes. In particular, the distributions of the numbers of large jumps in intervals were characterized in Jiao et al.\ (2017). The analysis is important for the study in mathematical finance as it allows one to describe in a unified way several recent observations on the bond markets such as the persistency of low interest rates together with the presence of large jumps.

\newpage

\section{A coupling of CBI-processes}

 \setcounter{equation}{0}

In this section, we give some characterizations of a coupling of CBI-processes constructed by the stochastic equation \eqref{s8.5}. Using this coupling we prove the strong Feller property and the exponential ergodicity of the CBI-process under suitable conditions. We shall follow the arguments of Li and Ma (2015). Suppose that $(\phi,\psi)$ are the branching and immigration mechanisms given respectively by \eqref{s2.13} and \eqref{s5.9} with $\nu(\d u)$ satisfying \eqref{s5.11}. Let $(P_t)_{t\ge 0}$ be the transition semigroup of the corresponding CBI-process defined by \eqref{s2.19} and \eqref{s5.10}.

\bgtheorem\label{ts10.1} If $\{x(t): t\ge 0\}$ and $\{y(t): t\ge 0\}$ are positive solutions to \eqref{s8.5} with $\bP\{x(0)\le y(0)\}= 1$, then $\{y(t)-x(t): t\ge 0\}$ is a CB-process with branching mechanism $\phi$. \edtheorem

\proof By Theorem~\ref{ts8.4} we have $\bP\{x(t)\le y(t)$ for all $t\ge 0\}= 1$. Let $z(t)= y(t)-x(t)$. {From} \eqref{s8.5} we have
 \beqnn
z(t) \ar=\ar z(0) + \sqrt{2c}\int_0^t\int_{x(s-)}^{y(s-)} W(\d s,\d u) - b\int_0^t z(s-)\d s
\cr
 \ar\ar\qquad
+ \int_0^t\int_0^\infty \int_{x(s-)}^{y(s-)} z \tilde{M}(\d s,\d z,\d u) \cr
 \ar=\ar
z(0) + \sqrt{2c}\int_0^t\int_0^{z(s-)} W(\d s,x(s-)+\d u) - b\int_0^t z(s-)\d s
\cr
 \ar\ar\qquad
+ \int_0^t\int_0^\infty \int_0^{z(s-)} z \tilde{M}(\d s,\d z,x(s-)+\d u), \cr
 \eeqnn
where $W(\d s,x(s-)+\d u)$ is a Gaussian time-space white noise with intensity $\d s\d u$ and $N(\d s,\d z, x(s-)+\d u)$ is a Poisson time-space random measure with intensity $\d sm(\d z)\d u$. That shows $\{z(t)\}$ is a weak solution to \eqref{s8.5} with $\eta(t)\equiv 0$. Then it is a CB-process with branching mechanism $\phi$. \qed

For $x\ge 0$ and $y\ge 0$, let $\{x(t): t\ge 0\}$ and $\{y(t): t\ge 0\}$ be the positive strong solutions to \eqref{s8.5} with $x(0)=x$ and $y(0)=y$. This construction gives a natural \index{coupling} \textit{coupling} of the CBI-processes. Let $\tau(x,y)= \inf\{t\ge 0: x(t)= y(t)\}$ be the \index{coalescence time} \textit{coalescence time} of the coupling. The distribution of this stopping time is given in the following theorem.

\bgtheorem\label{ts10.2} Suppose that Condition~\ref{ts3.5} holds. Then for any $t\ge 0$ we have
 \beqlb\label{s10.1}
\bP\{\tau(x,y)\le t\}= \bP\{y(t)= x(t)\}= \exp\{-|x-y|\bar{v}_t\},
 \eeqlb
where $t\mapsto \bar{v}_t$ is the unique solution to \eqref{s3.19} with singular initial condition $\bar{v}_{0+}= \infty$. \edtheorem

\proof It suffices to consider the case of $y\ge x\ge 0$. By Theorem~\ref{ts10.1} the difference $\{y(t) - x(t): t\ge 0\}$ is a CB-process with branching mechanism $\phi$. By Theorem~\ref{ts3.4} the probability $\bP\{\tau(x,y)\le t\}= \bP\{y(t)= x(t)\}$ is given by \eqref{s10.1}. \qed

\bgtheorem\label{ts10.3} Suppose that Condition~\ref{ts3.5} holds. Then for $t>0$ and $x,y\ge 0$ we have
 \beqlb\label{s10.2}
\|P_t(x,\cdot)-P_t(y,\cdot)\|_{\rm var}
 \le
2(1-\e^{-\bar{v}_t|x-y|})
 \le
2\bar{v}_t|x-y|,
 \eeqlb
where $\|\cdot\|_{\rm var}$ denotes the total variation norm. \edtheorem

\proof Let $\{x(t): t\ge 0\}$ and $\{y(t): t\ge 0\}$ be given as above. Since $\{y(t) - x(t): t\ge 0\}$ is a CB-process with branching mechanism $\phi$, for any bounded Borel function $f$ on $[0,\infty)$, we have
 \beqnn
|P_tf(x)-P_tf(y)|
 \ar=\ar
\big|\bP[f(x(t))] - \bP[f(y(t))]\big| \ccr
 \ar\le\ar
\bP\big[|f(x(t)) - f(y(t))|1_{\{y(t)\neq x(t)\}}\big] \ccr
 \ar\le\ar
\bP\big[(|f(x(t))| + |f(y(t))|)1_{\{y(t)\neq x(t)\}}\big] \ccr
 \ar\le\ar
2\|f\|\bP\{y(t) - x(t)\neq 0\} \ccr
 \ar=\ar
2\|f\|(1-\e^{-\bar{v}_t|x-y|}).
 \eeqnn
where the last equality follows by Theorem~\ref{ts10.2}. Then we get \eqref{s10.2} by taking the supremum over $f$ with $\|f\|\le 1$. \qed

By Theorem~\ref{ts10.3}, for each $t>0$ the operator $P_t$ maps any bounded Borel function on $[0,\infty)$ into a bounded continuous function, that is, the transition semigroup $(P_t)_{t\ge 0}$ satisfies the \index{strong Feller property} \textit{strong Feller property}.

\bgtheorem\label{ts10.4} Suppose that $b>0$. Then the transition semigroup $(P_t)_{t\ge 0}$ has a unique stationary distribution $\eta$ given by
 \beqlb\label{s10.3}
L_\eta(\lambda)
 =
\exp\bigg\{- \int_0^\infty\psi(v_s(\lambda))\d s\bigg\}
 =
\exp\bigg\{-\int_0^\lambda \frac{\psi(z)}{\phi(z)}\d z\bigg\},
 \eeqlb
and $P_t(x,\cdot)\to \eta$ weakly on $[0,\infty)$ as $t\to \infty$ for every $x\ge 0$. Moreover, we have
 \beqlb\label{s10.4}
\int_{[0,\infty)} y\eta(\d y)
 =
b^{-1}\psi^\prime(0),
 \eeqlb
where $\psi^\prime(0)$ is given by \eqref{s5.13}.
\edtheorem

\proof Since $b>0$, we have $\phi(z)\ge 0$ for all $z\ge 0$, so $t\mapsto v_t(\lambda)$ is decreasing. By Corollary~\ref{ts3.2} we have $\lim_{t\to \infty}v_t(\lambda)= 0$. {From} \eqref{s3.6} it follows that
 \beqnn
\int_0^t\psi(v_s(\lambda))\d s
 =
\int_{v_t(\lambda)}^\lambda \frac{\psi(z)}{\phi(z)}\d z.
 \eeqnn
In view of \eqref{s5.10}, we have
 \beqnn
\lim_{t\to \infty}\int_{[0,\infty)} \e^{-\lambda y} P_t(x,\d y)
 \ar=\ar
\exp\bigg\{- \int_0^\infty\psi(v_s(\lambda))\d s\bigg\} \cr
 \ar=\ar
\exp\bigg\{- \int_0^\lambda \frac{\psi(z)}{\phi(z)}\d z\bigg\}.
 \eeqnn
By Theorem~\ref{ts1.2} there is a probability measure $\eta$ on $[0,\infty)$ defined by \eqref{s10.3}. It is easy to show that $\eta$ is the unique stationary distribution for $(P_t)_{t\ge 0}$. The expression \eqref{s10.4} for its first moment follows by differentiating both sides of \eqref{s10.3} at $\lambda=0$. \qed

\bgtheorem\label{ts10.5} Suppose that $b>0$ and Condition~\ref{ts3.5} holds. Then for any $x\ge 0$ and $t\ge r> 0$ we have
 \beqlb\label{s10.5}
\|P_t(x,\cdot)-\eta(\cdot)\|_{\rm var}
 \le
2[x+b^{-1}\psi^\prime(0)]\bar{v}_r\e^{b(r-t)},
 \eeqlb
where $\eta$ is given by \eqref{s10.3}. \edtheorem

\proof Since $\eta$ is a stationary distribution for $(P_t)_{t\ge 0}$, by Theorem~\ref{ts10.3} one can see
 \beqnn
\|P_t(x,\cdot) - \eta(\cdot)\|_{\rm var}
 \ar=\ar
\bigg\|\int_{[0,\infty)} [P_t(x,\cdot)-P_t(y,\cdot)] \eta(\d y)\bigg\|_{\rm var}
\cr
 \ar\le\ar
\int_{[0,\infty)}\|P_t(x,\cdot) - P_t(y,\cdot)\|_{\rm var}\eta(\d y) \cr
 \ar\le\ar
2\bar{v}_t\int_{[0,\infty)} |x-y| \eta(\d y) \cr
 \ar\le\ar
2\bar{v}_t\int_{[0,\infty)} (x+y) \eta(\d y) \ccr
 \ar=\ar
2[x+b^{-1}\psi^\prime(0)]\bar{v}_t,
 \eeqnn
where the last equality follows by \eqref{s10.4}. The semigroup property of $(v_t)_{t\ge 0}$ implies $\bar{v}_t = v_{t-r}(\bar{v}_r)$ for any $t\ge r>0$. By \eqref{s3.4} we see $\bar{v}_t = v_{t-r}(\bar{v}_r)\le \e^{b(r-t)}\bar{v}_r$. Then \eqref{s10.5} holds. \qed

The result of Theorem~\ref{ts10.1} was used to construct flows of CBI-processes in Dawson and Li (2012). Clearly, the right-hand side of \eqref{s10.5} decays exponentially fast as $t\to \infty$. This property is called the \index{exponential ergodicity} \textit{exponential ergodicity} of the transition semigroup $(P_t)_{t\ge 0}$. It has played an important role in the study of asymptotics of the estimators for the $\alpha$-stable CIR-model in Li and Ma (2015).

 \newpage\addcontentsline{toc}{section}{\bf Index}

\printindex

\end{document}